\documentclass[10pt]{article}
\usepackage[latin1]{inputenc} 
\usepackage{latexsym,amsmath,amsfonts,euscript} 
\newtheorem{Theorem}{Theorem}[section] 
\newtheorem{Prop}[Theorem]{Proposition} 
\newtheorem{Def}[Theorem]{Definition} 
\newtheorem{Lemma}[Theorem]{Lemma} 
\newtheorem{Coro}[Theorem]{Corollary} 
\newtheorem{Remark}[Theorem]{Remark} 
 
\begin{document}

\newcommand{\finishproof}{\hfill $\Box$ \vspace{5mm}} 
\newcommand{\cI}{{\cal I}} 
\newcommand{\al}{\alpha} 
\newcommand{\be}{\beta} 
\newcommand{\ga}{\gamma} 
\newcommand{\Ga}{\Gamma} 
\newcommand{\ep}{\epsilon} 
\newcommand{\de}{\delta} 
\newcommand{\De}{\Delta} 
\newcommand{\ka}{\kappa} 
\newcommand{\la}{\lambda} 
\newcommand{\te}{\theta} 
\newcommand{\om}{\omega} 
\newcommand{\si}{\sigma} 
\newcommand{\bg}{{\bar g}} 
\newcommand{\bs}{{\bar s}} 
\newcommand{\bL}{{\bar L}} 
\newcommand{\tL}{{\tilde L}}
\newcommand{\tK}{{\tilde K}}
\newcommand{\tg}{\tilde g} 
\newcommand{\eqdef}{\stackrel{\rm def}{=}} 
\newcommand{\const}{\mathop{\rm const}} 
\newcommand{\sg}{\mathop{\rm sgrad}} 
\newcommand{\id}{{\bf 1}} 
\newcommand{\dds}{\frac{d}{ds}|_{s=0}} 
\newcommand{\p}{{\partial}} 
\newcommand{\C}{\mathbb C}
\newcommand{\LL}{\mathbb L}

\newcommand{\A}{{\mathbb A}}
\newcommand{\Z}{\mathbb Z} 
\newcommand{\N}{\mathbb N} 
\newcommand{\R}{\mathbb R}
\newcommand{\D}{\mathbb D} 
\newcommand{\Q}{\mathbb Q} 
\newcommand{\T}{\mathbb T}
\newcommand{\CC}{{\mathbb C}}
\newcommand{\tI}{\tilde I}

\renewcommand{\mod}{\mathop{\rm mod}} 
\renewcommand{\theequation}{\thesection .\arabic{equation}} 
\renewcommand{\arraystretch}{1.3}

\renewcommand{\theequation}{\thesection .\arabic{equation}}  
\renewcommand{\arraystretch}{1.3}

\title{Invariants of isospectral deformations and spectral rigidity. } 
\author{G.Popov\thanks{Partially supported by Agence Nationale de la Recherche, Projet 
"RESONANCES": GIP ANR-06-BLAN-0063-03.}, P.Topalov\thanks{Partially supported by the NSF grant DMS 0901443.}} 
%\thanks{Subject classification: 58J50, 70H06} 
%\date{12.12.2010} 
\maketitle 

\begin{abstract} 
\noindent
We introduce a notion of weak isospectrality for continuous deformations. 
Consider the Laplace-Beltrami operator on  
a compact Riemannian manifold  with boundary with  Robin
boundary conditions. Given a Kronecker invariant
torus $\Lambda$ of the billiard ball map with a Diophantine vector of rotation 
we prove that certain
integrals on  $\Lambda$ involving the function  in  the Robin boundary 
conditions remain
constant under weak isospectral deformations.  To this end we
construct continuous families of quasimodes associated with
$\Lambda$. We obtain also isospectral invariants of the Laplacian with
a real-valued potential 
on a compact manifold for continuous deformations of the potential. These invariants are obtained from the first  Birkhoff invariant of the  microlocal  monodromy operator associated to $\Lambda$. 
As an application we prove spectral rigidity of the Robin boundary 
conditions
in the case of Liouville billiard 
tables of dimension two in the presence of a $(\Z/2\Z)^2$ group of  symmetries.\\

\noindent
%Preliminary version 16.05.2009\\
\end{abstract}

\section{Introduction} 

This is a part of a project (cf. \cite{PT1}-\cite{PT3})  concerned with the spectral rigidity and integral geometry 
 of compact Liouville billiard tables of dimensions $n\ge 2$.
 The general strategy is   first to  find  a list of spectral invariants and then 
 to prove  for certain manifolds that these invariants imply spectral rigidity. 
 
%The main tool is the wave-trace formula of Selberg,  Duistermaat and Guillemin and of Guillemin and Melrose for manifolds with boundary. 
Substantial progress in the  inverse spectral problem has been done recently due to the wave-trace formula
 \cite{G}, \cite{G-M1}, \cite{IS}, \cite{ISZ}, \cite{MM}, \cite{Z1}-\cite{Z5}, and the semi-classical trace formulas  \cite{IS}, \cite{ISZ}, \cite{GPU}, \cite{GU}. The wave-trace formula, known in  physics as the Balian-Bloch formula and treated rigorously by Y. Colin de Verdi\`ere \cite{CV1} and by J. Duistermaat and V. Guillemin \cite{D-G} (see also \cite{G-M2}, \cite{MM}, \cite{P1}, \cite{S-V}, \cite{SZ}), as well its semi-classical analogue - the Gutzwiller trace formula relate the spectrum of the operator with different invariants of the corresponding closed geodesics  such as their  lengths and the spectrum of the linear Poincaré map. 
It is especially fruitful for deformations of closed manifolds with negative curvature.  V. Guillemin and D. Kazhdan \cite{G-K} and C. Croke and V. Sharafutdinov \cite{C-S} proved that such manifolds are spectrally rigid. There are several positive  results on the  inverse spectral problem for manifolds of non-negative curvature. 
 It has been proved  in \cite{G}, \cite{Z1}, \cite{Z2}, \cite{IS}, \cite{ISZ}, that for certain non-degenerate closed geodesics one can extract the Birkhoff Normal Form (BNF) from the singularity expansions of the wave-trace and even to reconstruct the boundary of analytic planar domains \cite{Z3} - \cite{Z6}. The wave-trace method requires  certain technical conditions such as simplicity of the length spectrum (a non-coincidence condition) and non-degeneracy of the corresponding closed geodesic and its iterates which are difficult to check.

We propose in this paper another approach which avoids the wave-trace formulas and works without assuming any  non-coincidence or non-degeneracy conditions. 
Our aim here is to present a simple idea of how  
quasimodes can be used in inverse spectral problems.  This idea works well for  isospectral 
deformations whenever  {\em continuous} with respect to the parameter 
of the deformation {\em quasimodes} can be constructed 
for the corresponding eigenvalue problem. In the currant paper, 
given a compact billiard table $(X,g)$ with a smooth Riemannian metric $g$  
and the corresponding Laplace-Beltrami operator on it, 
we consider   continuous deformations 
either of the function $K$ in the Robin boundary condition or of a  real-valued 
potential $V$ on $X$ (further applications of the method will be discussed in Sect. \ref{subsec:FurtherRemarks}). 
To construct quasimodes 
we assume that there is an  exponent  $B^m$, $m\ge 1$, 
of the corresponding billiard ball map $B$  which 
admits an invariant Kronecker torus $\Lambda$ with a  Diophantine vector of rotation $\omega$. 
By a Kronecker torus we mean   an embedded  Lagrangian submanifold $\Lambda$ of the coball 
bundle of the boundary, diffeomorphic to the torus $\T^{n-1}:=\R^{n-1}/2\pi\Z^{n-1}$ and 
invariant with respect to $B^m$,  and such that the restriction of 
$B^m$ to $\Lambda$  is smoothly conjugated to the rotation $R_\omega$ 
with a constant vector $\omega$ on $\T^{n-1}$. 
Any regular invariant torus of a completely integrable system,  is a Kronecker torus. Moreover, the  Kolmogorov-Arnold-Moser (KAM) theorem provides  families of Kronecker tori  for  close to completely integrable systems. These tori have  Diophantine vectors of rotation and their union is of a positive Lebesgue measure in the phase space.  
If the deformation is isospectral we prove that  certain integrals of the function 
$K$  (of the potential $V$)  on any Kronecker torus $\Lambda$ with a Diophantine vector of rotation 
remain constant under the deformation. These integrals are related to the first Birkhoff invariant  of the corresponding microlocal monodromy operator associated to $\Lambda$. In comparison with the wave-trace method \cite{G}, \cite{G-M1}, \cite{IS}, \cite{ISZ},  \cite{Z1}-\cite{Z5}, here, instead of looking for all the  coefficients  (or  Birkhoff invariants) of the singular expansion of the wave-trace related to a single  closed geodesic $\gamma$ and its iterates, we make use only of  the first Birkhoff invariant of the microlocal monodromy operator but for a large  family of Kronecker tori $\Lambda$, which can be considered as a Radon transform associated to that family. In particular, our method requires only finite regularity of the function $K$ (the potential $V$).

Typical examples of completely integrable billiard tables are the Liouville billiard tables (L.B.T.s) \cite{PT3}.  
In the case of L.B.T.s   we treat these integral invariants as  values 
of a  suitable Radon transform. 
Then the spectral rigidity follows  from the injectivity of the Radon transform. 
Liouville billiard tables of dimension two have  been studied in \cite{PT1}. 
Liouville billiard tables of  dimension $n\ge 2$ are introduced in \cite{PT3}, where  
the integrability of the corresponding billiard ball map is obtained  using 
a simple variational principle. A typical example of a  L.B.T.  is the  interior of the $n$-axial ellipsoid in $\R^n$  \cite{PT3}. The injectivity of the Radon transform for L.B.T.s in higher
dimensions is investigated  in  \cite{PT2}. 
%An earlier version of this paper is given in 

To summarize, we point out that the method can be applied whenever a continuous quasimode of the problem exists. The notion of  continuous quasimodes is introduced in Sect. \ref{sec:quasimodes-isospectral}. The construction of the continuous quasimodes here is related on Birkhoff normal forms.

A billiard table $(X,g)$ is  a smooth compact Riemannian manifold of dimension    
${\rm dim}\, X = n \ge  2$ equipped with a smooth Riemannian metric $g$ 
and with a $C^\infty$ boundary   
$\Gamma := \partial X\ne\emptyset$.    
Let $\Delta$ be the ``positive'' Laplace-Beltrami operator on $(X,g)$.   
Given  a real-valued function $K\in C^\ell(\Gamma, \R)$,   $\ell\ge 0$, 
we  consider the  operator $\Delta$  with domain     
\begin{equation}  
D:= \left\{u\in H^2(X):\, \frac{\partial u}{\partial \nu}|_\Gamma = K  
u|_\Gamma\right\}\, ,  
\label{domain}
\end{equation}  
where $\nu(x)$, $x\in \Gamma$, is the inward unit normal to  
$\Gamma$ with respect to the metric $g$. We denote this operator by  
$\Delta_{g,K}$.   
It is a selfadjoint operator in $L^2(X)$ with  
discrete spectrum   
\[  
{\rm Spec}\,  \Delta_{g,K}:=\{ \lambda_1 \le \lambda_2 \le \cdots  
\}\, ,  
\]  
where each eigenvalue $\lambda=\lambda_j$  is repeated according to its  
multiplicity, and it solves the spectral problem  
\begin{equation}  
\left\{  
\begin{array}{rcll}  
\Delta\,  u\ &=& \ \la\,  u \, \quad \mbox{in}\ $X$\, , \\  
\displaystyle\frac{\partial u}{\partial \nu}|_\Gamma \  &=& \ K\, u|_\Gamma \, .  
\end{array}  
\right.  
                                         \label{thespectrum}  
\end{equation}

\subsection{Invariants of isospectral families}
The method we use requires only finite regularity of $K$. 
Fix $\ell \ge 0$, and consider a continuous family of  
real-valued functions $K_t$ in the Hölder space $C^\ell(\Gamma,\R)$, $t\in [0,1]$, which means that  
the map $[0,1]\ni t\mapsto K_t$  is continuous
in  $C^\ell(\Gamma, \R)$. To simplify the notations we denote by 
$\Delta_t$ the corresponding operator $\Delta_{g,K_t}$. The family  $\Delta_t$  is said
to be isospectral if
\begin{equation}
\forall\,  t \in [0,1]\, ,\  {\rm Spec}\left(\Delta_t \right)\, = \, 
{\rm Spec}\left(\Delta_0 \right)\, .
\label{isospectral}
\end{equation}
We are going to introduce a weaker notion of isospectrality. Consider the union ${\mathcal I}$ 
of infinitely many disjoint intervals going to infinity,  with lengths tending to zero, and which are polynomially separated. More precisely, fix  two positive constants $d>1/2$ and $c>0$, and suppose that
 
\begin{enumerate}
\item[$(\mbox{H}_1)$] {\em 
 The set $\mathcal{I}\subset (0,\infty) $ 
 is a union of infinitely many disjoint  intervals $[a_k,b_k]$, $k\in\N$,  such that 
\begin{enumerate}
\item[$\bullet$]
$\displaystyle{\lim\,  a_k\, =\, \lim\,  b_k\,  =\,  +\infty}$,
\item[$\bullet$] 
$\lim\,  (b_k -a_k)\,  =\,  0$, 
\mbox{and} 
\item[$\bullet$]
$a_{k+1} - b_{k}\, \ge \, c b_k^{-d}$ \quad for  any 
$k\in \N$. 
\end{enumerate}}
\end{enumerate}
More generally, fix an integer $s\ge 0$ and instead of $(\mbox{H}_1)$ consider the stronger condition $(\mbox{H}_1)_s$, 
where the second assumption is replaced by 
\begin{equation}
\label{strong-iso}
 \lim\, a_k^{s/2}\, (b_k-a_k) \ =\  0 \quad \mbox{as} \quad k\to \infty\,   .
\end{equation}
We impose the following {\em ``weak isospectral assumption''}: 
\begin{enumerate}
\item[$(\mbox{H}_2)$]
{\em There is  $a\ge 1$ such that  
$\forall\,  t \in [0,1]\, ,\  {\rm Spec}\left(\Delta_t
\right)\, \cap [a,+\infty)\  
\subset \   {\mathcal I} \, .$}
\end{enumerate}
Using the asymptotic of the eigenvalues $\lambda_j$ as $j\to \infty$ we will show in Sect. \ref{sec:quasimodes-isospectral} that the conditions (H$_1$)-(H$_2$)
are ``natural'' for any $d> n/2$ and $c>0$. By ``natural'' we mean that for any $d> n/2$ and $c>0$ the isospectral condition \eqref{isospectral} implies that  
there exists a union of infinitely many disjoint intervals $\mathcal I$ such that (H$_1$)-(H$_2$) are satisfied -- see Lemma \ref{lemma:natural}
and Remark \ref{Rem:natural} for  details.
Note that (H$_1$)-(H$_2$) and $(\mbox{H}_1)_s$-(H$_2$) are well adapted for the semi-classical analysis. 

The elastic reflection of geodesics at $\Gamma$ determines
continuous curves on $X$ called {\em billiard trajectories} as well as a discontinuous dynamical system on $T^*X$ --
the {\em billiard flow}. The latter 
 induces  a discrete dynamical system $B$ defined  on an open subset $\widetilde B^\ast \Gamma$ of the open  coball bundle $B^\ast \Gamma$ of 
$\Gamma$ called billiard ball map (see Sect.  2.1 for a definition). The map $B:\widetilde B^\ast \Gamma\to B^\ast \Gamma$  is symplectic.  
We suppose  also that there is an integer $m\ge 1$ such that the map $P=B^m$   
 admits an invariant Kronecker
torus of a  vector of rotation $\omega$ satisfying the following $(\kappa,\tau)$-Diophantine  condition: 
\begin{equation}  
\left\{
\begin{array}{crl}  
\mbox{\em There is $\kappa >0$
and $\tau > n-1$ such that $\forall\, (k,k_{n})\in \Z^{n},\ 
k=(k_1,\ldots,k_{n-1})\neq 0\, $:  }\\ [0.3cm]                      
|\langle\omega,k\rangle + k_{n}|\ \ge \ \kappa\, \left(\sum_{j=1}^{n-1}
|k_j|\right)^{- \tau}\, .
                       \label{sdc}                     
\end{array}  
\right.                     
\end{equation}
For example, if ${\rm dim}\, X=2$, and $\Gamma$ is locally strictly geodesically convex (with respect to the outward normal $-\nu(x)$), then by a result of Lazutkin the union of the invariant circles of $P=B$ with  Diophantine numbers of rotation is of a positive Lebesgue measure in $T^\ast\Gamma$. If $\gamma$ is a suitable elliptic broken geodesic with $m$ vertices in $X$, ${\rm dim}\, X\ge 2$,  applying the KAM theorem,  we get a family of invariant tori of $P=B^m$ with Diophantine vectors  of rotation the union of which has 
a positive measure in $T^\ast\Gamma$.

More precisely, we impose the following dynamical condition. 
\begin{enumerate}
\item[(H$_3$)] {\em There is  an embedded submanifold $\Lambda$ of $ B ^\ast \Gamma$  diffeomorphic
to $\T^{n-1}$  and an integer $m\ge 1$  such that $B^j(\Lambda)$, $0\le j\le m-1$, belong to the domain of definition $\widetilde B^\ast \Gamma$ of $B$, $\Lambda$ is   invariant with respect to $P=B^m$,  and   the restriction of $P$ to $\Lambda$ is
$C^\infty$ conjugated to the  translation  
$R_{2\pi \omega}(\varphi) =\varphi - 2\pi \omega \, \mbox{\rm (mod $2\pi$)}$ in
$\T^{n-1}$, 
where
$\omega$ satisfies the Diophantine condition  (\ref{sdc}). }
\end{enumerate}  
Then 
$\Lambda\subset B ^\ast  \Gamma$ is Lagrangian (see \cite{Her}, Sect. I.3.2).  
Moreover, it follows from (\ref{sdc}) that $P$ is uniquely ergodic on $\Lambda$, i.e. there is a  unique probability measure $d\mu$ on $\Lambda$ (associated with the corresponding Borelian $\sigma$-algebra) which is invariant under $P$. Obviously,  $d\mu$ is  the pull-back of the Lebesgue measure $(2\pi)^{1-n}d\varphi$ on $\T^{n-1}$ via the diffomorphism in (H$_3$). Set $\Lambda_j=B^j(\Lambda)$, $j=0,1,\ldots,m-1$, and $d\mu_j=(B^{-j})^\ast(d\mu)$. Then $\Lambda_j$ is a Kronecker invariant torus of $P$ with a vector of rotation $2\pi\omega$ and $d\mu_j$ is the unique probability measure on it which is invariant with respect to $P$. 
Denote by  $\pi_\Gamma:
T^\ast \Gamma \to \Gamma$ the canonical projection. 
Given $(x,\xi) \in B ^\ast  \Gamma$, we
denote by $\xi^+\in T_x^\ast X$ the corresponding outgoing unit covector which means that the restriction of $\xi^+$ to $T_x\Gamma$ equals $\xi$ and $\langle \xi^+,\nu(x)\rangle\ge 0$, where $\langle \cdot,\cdot\rangle$ stands for the paring between vectors and covectors, and we define 
 \[
 \theta = \theta(x,\xi)\in [0,\pi/2] \quad \mbox{by} \quad \sin \theta =\langle \xi^+,\nu(x)\rangle .
\] 
We require only finite $C^\ell$-Holder regularity on $K$. 
Fix  $d>1/2$ and $\tau > n-1$ and let 
\begin{equation}
\ell > ([2d]+1)(\tau + 2) +2n+(n-1)/2\, ,
\label{regularity-index}
\end{equation}
where $[p]$ stands for the entire part of the real number $p$. In what follows
$d$ will be  the exponent in  (H$_1$), and
$\tau>n-1$  the exponent in the Diophantine  condition
(\ref{sdc}). Our main result is:
 
\begin{Theorem}\label{main}
Let $\Lambda$ be an invariant  Kronecker torus of $P=B^m$ with a vector of rotation 
$2\pi \omega$, where $\omega$ is $(\kappa,\tau)$-Diophantine.   Fix $\ell$ and $d>1/2$ such that  (\ref{regularity-index}) holds. 
Let 
\[
[0,1]\ni t \mapsto K_t\in C^\ell(\Gamma,\R)\ , 
\]
 be a continuous family of real-valued functions 
on $\Gamma$ such that $\Delta_t$  satisfy  the isospectral condition 
 $(H_1)-(H_2)$ with exponent $d$. 
 Then
\begin{equation}
\forall\, t\in [0,1],\quad \displaystyle \sum_{j=0}^{m -1}\int_{\Lambda_j}
\frac{K_t\circ \pi_\Gamma}{\sin \theta}\, d\mu_j \ =\
\sum_{j=0}^{m -1}
\int_{\Lambda_j}
\frac{K_0\circ \pi_\Gamma}{\sin \theta}\, d\mu_j \, . 
                                           \label{theinvariant}
\end{equation}
\end{Theorem}
%Before  giving  applications of the  theorem we would like to make some comments on it.  
The invariant in \eqref{theinvariant} will be  obtained from the first Birkhoff invariant (for a definition see Remark \ref{rem:birkhoff}) of the corresponding microlocal  monodromy operator. To prove that higher order Birkhoff invariants are isospectral invariants, one should impose the isospectral condition  $(H_1)_s-(H_2)$ for some $s\in\N$ (see Sect. \ref{subsec:FurtherRemarks}). 

Theorem \ref{main} is inspired by 
a result of  Guillemin and Melrose \cite{G-M1, G-M2}. They consider a connected clean 
submanifold $\Lambda$ of fixed points of $P= B^m$, $m\ge 2$, satisfying the
so called ``non-coincidence'' condition\footnote{The integer $m$ is not necessarily the minimal exponent so that $P(\varrho)= B^m(\varrho)$ for $\varrho\in\Lambda$}. Let $T_{\Lambda,m}$ be the  common length 
of  the family of closed broken geodesics issuing from $\Lambda$ and having  $m$  reflexions  at $\Gamma$. 
The ``non-coincidence'' condition means that  the broken geodesics of that family  are the only closed
generalized geodesics in $X$ of length $T_{\Lambda,m}$. Under this condition, 
Guillemin and Melrose prove that 
if $K_j$, $j=1,2$, are two real-valued $C^\infty$ functions on
$\Gamma$
such that $\mbox{Spec}\left(\Delta_{g,K_1}\right)
=\mbox{Spec}\left(\Delta_{g,K_2}\right)$, then certain integrals of $K_j/\sin\theta$, $j=1,2$,  on $\Lambda$ are equal. In the case when $X\subset \R^2$ 
is the interior of an ellipse $\Gamma$ they obtain an infinite
sequence of confocal ellipses  $\Gamma_j\subset X$  tending
to $\Gamma$ such that the corresponding invariant circles $\Lambda_j$
of $B$ satisfy the non-coincidence condition. 
 As a consequence they  obtain  in \cite{G-M1} spectral
rigidity of (\ref{thespectrum}) in the case of the ellipse for
$C^\infty$ functions $K$ which are invariant with respect to the
symmetries of the ellipse. 
The main tool in the proof is the trace formula for the wave equation
with Robin boundary conditions in $X$ (see \cite{G-M2}). 
This result was generalized in \cite{PT1} for two-dimensional 
Liouville billiard
tables of classical type. 

There is little  hope to apply the wave-trace formula 
in our situation.
An invariant Kronecker torus $\Lambda$ of the billiard ball map $B$ could  be approximated with periodic points of $P=B^m$ using a variant of 
the Birkhoff-Lewis theorem and a ``Birkhoff normal form'' of $P$ near
$\Lambda$. Unfortunately, we do not know if the corresponding closed broken
geodesics are non-degenerate. Moreover, it will be  difficult to check 
if  the non-coincidence condition holds even in the case of  ellipsoidal billiard tables.

We propose a simple idea  which  relies on a quasimode construction.
It is natural to use quasimodes for this kind of problems since
quasi-eigenvalues  are close to  eigenvalues and they
contain a lot of geometric information. A quasimode with discrepancy $\varepsilon >0$ of a selfadjoit operator $L$ in $L^2(X)$  is a pair $(\lambda, u)$, where $\lambda$ is a real number, $u$ belongs to the domain of definition of $L$, $\|u\|=1$, and $\|Lu - \lambda u\|\le \varepsilon$. By the spectral theorem, this implies
\begin{equation}
\mbox{dist} \left(\mbox{Spec} \left(L\right), \lambda\right) \le \varepsilon\, .
\label{spectral-distance}
\end{equation}
In order to prove (\ref{theinvariant}), we  construct {\em continuous} families of 
quasimodes $(\mu_q(t)^2 , u_q(t,\cdot))$ of $\Delta_t$ with discrepancies $\varepsilon_q=C_M \mu_q^{-M}$, where
 $M = [2d]+1$  is the
entire part of $2d+1$, $d$ is the exponent in (H$_1$),  $q$ belongs to an unbounded index set 
${\mathcal M} \subset \Z^{n}$ and $C_M$ is a positive constant {\em independent} of $q$ and of $t$. By continuous quasimodes we mean that the quasi-eigenvalues $\mu_q(t)^2$ are  {\em continuous functions} 
of  $t\in [0,1]$. More precisely, 
we obtain  by  Theorem 2.2 quasimodes $(\mu_q(t)^2 , u_q(t,\cdot))$ of $\Delta_t$ such that 
$$
\mu_q(t) = \mu_q^0 +c_{q,0} + c_{q,1}(t)(\mu_q^0)^{-1} 
+ \cdots + c_{q,M}(t)(\mu_q^0)^{-M}   \, ,
$$
where $\mu_q^0$ and $c_{q,0}$ are independent of $t$, 
$\lim_{|q|\to \infty}\mu_q^0 = +\infty $,  and  
$c_{q,j}$, $q\in {\mathcal M}$, 
is an uniformly bounded sequence of  continuous functions in $t\in [0,1]$. 
The function $c_{q,1}$ is of the form
$$
c_{q,1}(t)=c'_{q,1} + c''_{1}\sum_{j=0}^{m-1}\,  \int_{\Lambda_j}\, 
\frac{K_t\circ \pi_\Gamma}{\sin \theta}\, d\mu_j \, , 
$$ 
where $c'_{q,1}$ and  $c''_{1}\neq 0$ are independent of $t$ and  $c''_{1}$ does
not depend on $q$ either.  Now \eqref{spectral-distance} reads
$$
\mbox{dist} \left(\mbox{Spec} \left(\Delta_{t}\right), \mu_q(t)^2\right) \le C_M (\mu_q)^{-M}\, .
$$
Since $M>2d$, it follows from   (H$_2$)  that the quasi-eigenvalues  $\mu_q(t)^2$, $|q|\ge q_0 \gg 1$,  belong
to the  union of intervals $[a_k-ca_k^{-d}/4,b_k+cb_k^{-d}/4]$ which do
not intersect each other  in view of (H$_1$). 
 Since $\mu_q(t)^2$  is continuous in $[0,1]$, it
can not jump from one interval to another. Using (H$_1$), this allows us to show that 
$$
 c_{q,1}(t)- c_{q,1}(0) =  
\mu_q^0|\mu_q(t) - \mu_q(0)| + O(\mu_q^0)^{-1}
$$
tend to $0$ as $|q|\to \infty $ for any $t\in [0,1]$, 
which proves (\ref{theinvariant}).  A more general argument involving the condition $(\mbox{H}_1)_s$ is given in Lemma \ref{lemma:invariants}. 

 Theorem \ref{quasimodes} provides a construction of  quasimodes, which depend continuously on $K\in C^\ell(\Gamma,\R)$, 
where $\ell > M(\tau +2)+2n+ (n-1)/2$ and  $\tau$ is the exponent in the Diophantine condition (\ref{sdc}). Choosing $M=[2d]+1$ we obtain that the index of regularity  $\ell$  of $K_t$ in Theorem \ref{main} should satisfy  (\ref{regularity-index}). 

The construction of continuous quasimodes is provided in Sect. \ref{sec:construction}. This section is quite long and we divide it into several subsections. 

In Sect. \ref{subsec:reduction} we reduce the spectral problem (\ref{thespectrum}) microlocally near  $\Lambda$ to an equation with respect to $(\lambda,v)\in \C\times C^\infty(\Gamma)$ of the form
\[
(W_0(\lambda) -\mbox{Id})v\ =\ O_M(\lambda^{-M-1})v,
\]
where  $W_0(\lambda)$  is the corresponding ``microlocal monodromy operator''. The idea of the reduction at the boundary is explained in the beginning of Sect. \ref{subsec:reduction}. 
The monodromy operator has the form 
$$
W_0(\lambda)\ =\ [Q^0(\lambda) + \lambda^{-1}Q^1(\lambda)]S(\lambda),
$$
where $S(\lambda)$ is a ``classical'' Fourier Integral Operator with a large parameter $\lambda$ ($\lambda$-FIO) whose canonical relation is the graph of the Poincar\'e map $P= B^m$, and with a $C^\infty$-smooth compactly supported amplitude in any local chart, $Q^0(\lambda)$ is a ``classical'' pseudodifferential operator with a large parameter $\lambda$ ($\lambda$-PDO) with a $C^\infty$-smooth compactly supported amplitude in any local chart, and $Q^1(\lambda)$ is a  $\lambda$-PDO with a  compactly supported amplitude of finite regularity. Moreover, $S(\lambda)$ and $Q^0(\lambda)$ do not depend on $K$, while $Q^1(\lambda)$ depends continuously on $K$. The symbol of $Q^1(\lambda)$ in any local chart has the form $\sum_{j=0}^{M-1}q_j^0(x,\xi)\lambda^{-j} + O_M(\lambda^{-M})$,  where $q_j^0$ are compactly supported functions. Moreover, $q_j^0$ are polynomials of $K$ and of its derivatives of order $\le 2j$ with smooth coefficients, while the reminder term $O_M(\lambda^{-M})$ depends continuously on $K$ and on its derivatives of order $\le 2M$. 
The symbols of finite regularity that we need are introduced in the Appendix. The reduction to the boundary  is similar to that in \cite{C-P} and it is close also to the construction of a parametrix for the mixed problem for the 
 wave equation in the $C^\infty$ case  which  has been done by Chazarain \cite{Ch} and V. Guillemin and R. Melrose in \cite{G-M2}. 
 A similar construction has been used recently by H. Hezari and S. Zelditch \cite{H-Z}. The reduction to the boundary is based on a construction  of a parametrix for the Helmholtz equation with Dirichlet boundary conditions which is done in Appendix A.1. 
 
Starting from a Birkhoff Normal Form of  $P$ in a neighborhood of $\Lambda$,  we find a Quantum Birkhoff Normal Form (shortly QBNF) of the monodromy operator $W_0(\lambda)$ in Sect. \ref{subsec:QBNF} solving at any step  a suitable homological equation. 
In particular, we obtain a spectral decomposition of $W_0(\lambda)$ modulo $O_M(\lambda^{-M})$, which allows us to construct the quasimodes.  The third term $c_{q,1}(t)$  of the corresponding quasi-eigenvalues is given by the first  Birkhoff invariant of the QBNF of the monodromy operator. It turns out  that higher order  Birkhoff invariants of $W_0(\lambda)$ are also isospectral invariants (see Sect. \ref{subsec:FurtherRemarks}). 
Similar  construction of quasimodes has been  provided in  \cite{C-P} and in  \cite{P}. The method of the ``microlocal  monodromy operator'' has been developed   in the  context of the wave-trace and the semi-classical trace formulas in  \cite{SZ}. It allows also to investigate the contribution of degenerate closed geodesics to the wave-trace (see also \cite{P1}).  

To get a normal form of $W_0(\lambda)$ we make use of Wiener spaces with weights ${\mathcal A}^s(\T^{n-1})$. A function in $\T^{n-1}$ belongs to such a space if the series of its Fourier coefficients multiplied by a suitable polynomial weight is absolutely convergent. These spaces are perfectly adapted for solving the homological equation (see Lemma \ref{lemma:homological}). At any iteration we lose exactly $\tau $ derivatives, which makes $M\tau$ in total. 
On the other hand, $C^q(\T^{n-1})\subset {\mathcal A}^s(\T^{n-1})\subset C^s(\T^{n-1})$ for any $s\ge 0$ and $q>s + (n-1)/2$, and the inclusion maps are continuous by a theorem of Bernstein. We need also $2n$ derivatives on the amplitude of the corresponding $\lambda$-FIO acting on $\T^{n-1}$ to prove $L^2$-continuity. Finally, to construct continuous with respect to $K\in C^l(\Gamma,\R)$ quasimodes we have to  assume that $l > 2M +\tau M+ (n-1)/2 +2n$.

\subsection{Applications. Spectral rigidity in the presence of $(\Z/2\Z)^2$ symmetries. }
Kronecker invariant tori usually 
appear  in Cantor families (with respect to the  
Diophantine vector of rotation $ \omega$ satisfying (\ref{sdc}) with $\tau>n-1$),  the union of which has positive
Lebesgue measure in $T^\ast \Gamma$, and Theorem 1.1 applies to any
single torus $\Lambda$ of that family.  
Consider for example a strictly convex bounded domain $X\subset \R^2$
with $C^\infty$ boundary $\Gamma$, and fix $\tau >1$. 
  It is known from
Lazutkin (cf. \cite{La}, Theorem 14.21) that 
for any  $0<\kappa \le \kappa_0 \ll 1$ 
there is a Cantor set $\Xi_\kappa\subset (0,\varepsilon_0]$, $\varepsilon_0\sim \sqrt{\kappa_0}\ll 1$,  of $(\kappa,\tau)$-Diophantine
numbers $\omega$    such that 
for each $\omega\in \Xi_\kappa$ there is a KAM (Kolmogorov-Arnold-Moser)  invariant circle
$\Lambda_\omega \subset B^\ast \Gamma$ of $B$ satisfying (H$_3$) with
$m=1$  and with a rotation number
$2\pi\omega$. 
Moreover, $\Xi_\kappa$ is of a positive Lebesgue measure in $(0,\varepsilon_0]$,   
the Lebesgue measure of the complement of  $\Xi := \cup\,  \Xi_\kappa$ in $(0,\varepsilon]$  is $O(\varepsilon)$ as $\varepsilon \to 0$, 
and so is the Lebesgue measure
of the complement to the union of the invariant circles in an
$\varepsilon$-neighborhood of $S^\ast \Gamma$ in $B^\ast \Gamma$. 
More generally, the result of Lazutkin holds for any  compact billiard table
$(X,g)$, ${\rm dim}\, X =2$, with  
connected boundary $\Gamma$ which is locally
strictly geodesically convex. The corresponding completely integrable map is obtained by means of the ``approximate interpolating Hamiltonian'' \cite{MM}. Applying Theorem \ref{main} to any KAM circle $\Lambda_\omega$, $\omega\in \Xi$,  we obtain the following 
\begin{Coro}\label{Coro:convex}
Let $(X,g)$, ${\rm dim}\, X =2$,
be a compact locally strictly
geodesically convex billiard  table  with a smooth connected  boundary $\Gamma$. Let
$[0,1]\ni t \mapsto K_t\in C^\ell(\Gamma,\R)$ 
 be a continuous family of functions  such
 that $\Delta_t$  satisfy  $(H_1)-(H_2)$. Let  $\ell >  ([2d]+1)(\tau+2) +9/2$, where  $d>1/2$ is the exponent in $(H_1)$. 
Then 
\begin{equation}
\forall \omega \in \Xi\, ,\ \forall\, t\in [0,1],
\quad \displaystyle \int_{\Lambda_\omega}
\frac{K_t\circ \pi_\Gamma}{\sin \theta}\, d\mu \ =\ 
\int_{\Lambda_\omega}
\frac{K_0\circ \pi_\Gamma}{\sin \theta} \, d\mu \, . 
                                           \label{theinvariant2}
\end{equation}
\end{Coro}
Another example can be 
obtained applying the KAM theorem to the Poincar\'e map $P=B^m$ of a 
periodic broken geodesic $\gamma$ with $m$ vertexes
(in any dimension $n\ge 2$). To this end we suppose that $\gamma$ is {\em elliptic}, that {\em there are no resonances of order $\le 4$}  and that the corresponding Birkhoff normal form of $P$ is {\em non-degenerate} (see Sect. \ref{subsec:bouncing-ball} for definitions). If the dimension of $X$ is two these conditions mean that $P$ is a twist map. Consider in more details the case of a bouncing ball trajectory $\gamma$ ($m=2$) with vertices $x_1$ and $x_2$ in a two-dimensional billiard table.  Assume  that there is a neighborhood $U$ of $\gamma$ in $X$ and two involutions $J_k:U\to U$, $k=1,2$,  acting as isometries on $U$,  such that $J_1(x_j)=x_j$ and $J_2(x_1)=x_2$ and suppose that $J_1$ and $J_2$ commute. Denote by $I(\Gamma\cap U)\cong\Z_2\oplus\Z_2$ the group of local  isometries generated by the restrictions of $J_1$ and $J_2$ to $\Gamma\cap U$. For any $f\in C(\Gamma\cap U)$ we consider the average 
\[
f^{\#}(x)=\frac{1}{4}\sum_{g\in I(\Gamma\cap U)}f(g.x)\, ,\ x\in \Gamma\cap U.
\]
\begin{Coro}\label{Coro:bouncing-ball}
Let $(X,g)$ be a compact billiard table of dimension ${\rm dim}\, X =2$ with a smooth  boundary $\Gamma$.
Let $\gamma$ be an elliptic  bouncing ball trajectory with end points $x_1$ and $x_2$ satisfying the hypothesis above.  Let
$[0,1]\ni t \mapsto K_t\in C^\ell(\Gamma,\R)$ 
 be a continuous family of real-valued functions on $\Gamma$ such
 that $\Delta_t$  satisfy  $(H_1)-(H_2)$ and let $\ell >  3[2d] +15/2$, where $d$ is the exponent in (H$_1$). 
Then for any $t\in (0,1]$ there is an infinite sequence $\{y_k\}_{k\in\N}$ in $\Gamma\setminus \{x_1\}$ such that $\lim y_k=x_1$ and $K_t^{\#}(y_k)=K_0^{\#}(y_k)$ for any $k\in \N$. In particular, 
the Taylor polynomials  of $K_t^{\#}$ at $x_j$, $j=1,2$,  of degree less than $[\ell]+1$ do not depend on $t\in[0,1]$. 
\end{Coro}
In particular, if $\Gamma$ is analytic and $K_0$ and $K_1$ are real-analytic and invariant with respect to $J_1$ and $J_2$, then $K_0=K_1$. Using higher order Birkhoff invariants one could remove the symmetries of $K_t$ (see Sect. \ref{subsec:FurtherRemarks}). 
Corollary \ref{Coro:bouncing-ball} will be proved in Sect.  \ref{subsec:bouncing-ball}. 

Similar results can be obtained from the singular expansion of the wave-trace at $\gamma$ and its iterates by the method developed by Zelditch \cite{Z3,Z4}. As we have mentioned above  the  two methods are of  different nature. 
Here we use only one Birkhoff invariant but for infinitely many invariant circles which requires the condition  $\lim\,  (b_k-a_k)=0$ in $(\mbox{H}_1)$. By the wave-trace method  one obtains   all the terms of the singular expansion of the wave-trace related to only one  closed geodesic $\gamma$ and its iterates, which would require the stronger condition  $\lim\, a_k^s (b_k-a_k)=0$ for any $s\in\N$.

Theorem \ref{main} can be applied as well in the completely integrable case, for example for  the ellipse 
or the ellipsoid, or more generally for Liouville billiard
tables of classical type \cite{PT1,PT2, PT3} in any dimension $n\ge 2$. 
We are going to prove spectral rigidity of the Robin boundary problem 
for two dimensional  Liouville billiard tables of classical type 
(see Sect. \ref{subsec:L.B.T.} for definitions).
Such billiard tables have a group of isometries 
$I(X)\cong\Z_2\oplus\Z_2$ which induces a group of isometries 
$I(\Gamma)\cong\Z_2\oplus\Z_2$ on the boundary. Given $f\in C(\Gamma)$ we consider  its average $f^{\#}$ with respect to $I(\Gamma)$. 
We denote by ${\rm Symm}^\ell (\Gamma)$ 
the space of all $C^\ell$ real-valued functions which are invariant with respect to 
$I(\Gamma)$, i.e. a smooth function $f$ belongs to ${\rm Symm}^\ell (\Gamma)$ if $f=f^{\#}$ . Applying Theorem 1.1 for Diophantine numbers of rotation $\omega$, we show  that {\em any continuous weakly isospectral deformation} of $K$ in 
${\rm Symm}^\ell (\Gamma)$ is {\em trivial}. 
More precisely, we have 
\begin{Coro}\label{Coro:liouville}
Let $(X,g)$, ${\rm dim}\, X =2$,
be a Liouville billiard table of classical type.  Let 
$K_t$, $t\in [0,1]$, 
 be a continuous family of real-valued functions in $C^\ell(\Gamma,\R)$
  such that $\Delta_t$  satisfy   $(H_1)-(H_2)$. Suppose that $\ell >  3[2d] +15/2$, where $d$ is the exponent in (H$_1$). 
Then $K_t^{\#}\equiv K_0^{\#}$ for any $t\in [0,1]$. In particular, if $K_0, K_1\in {\rm Symm}^\ell (\Gamma)$
then $K_0=K_1$.
\end{Coro}

\noindent
It seems that even for the ellipse this result has not been known.
%Using Lemma 2.1 and Corollary \ref{Coro:liouville} 
%we obtain that {\em any continuous  isospectral deformation} of $K$ 
%in the sense of (\ref{isospectral}) in 
%${\rm Symm}(\Gamma)$ is  {\em trivial}. 
We point out that the Liouville 
billiard tables that we consider are not analytic in general 
and  the methods used in \cite{G-M1} and \cite{PT1} can not be applied. 
Corollary \ref{Coro:liouville} will be proved in Sect. \ref{subsec:L.B.T.}.

In the same way we treat the operator 
$\Delta_t = \Delta + V_t$ in $X$ with fixed Dirichlet or Robin (Neumann) 
boundary conditions on $\Gamma$,
where $V_t\in C^\ell(X)$, $t\in [0,1]$, 
is a continuous family of real-valued potentials in $X$. 
The corresponding results are 
proved in Sect. \ref{subsec:potential}.  Injectivity of the Radon transform and 
spectral rigidity of Liouville billiard tables in higher dimensions is investigated in 
\cite{PT2}.

\section{Quasimodes and spectral invariants}\label{sec:quasimodes}
\setcounter{equation}{0}

\subsection{Quasimodes and isospectral deformations}\label{sec:quasimodes-isospectral}

First we shall show that the isospectral condition (H$_1$)-(H$_2$) is natural for any $d > n/2$ and $c>0$. 
Given $d>n/2$, $c>0$, and $\al\gg 1$ we consider the set
$$
\tilde{\mathcal I}\equiv \tilde{\mathcal I}(\Delta_0)  :=\  \left\{ \lambda \ge \al:\, 
\left|\, \mbox{Spec}\left(\Delta_0 \right)\, -\,  \lambda\right|\, \le\, 
2c\lambda^{-d}\right\} \, .
$$
Let us write $\tilde{\mathcal I}$ as a disjoint union of connected intervals $[\overline a_k,\overline b_k]$,
$1\le k\le m\le \infty$.\footnote{We do not exclude a priori  the case when $m$ is finite and $\overline b_m=\infty$.}
\begin{Lemma}\label{Lem:infinite_intervals}
The intervals $[\overline a_k,\overline b_k]$ are infinitely many and $\lim_{k\to\infty}\overline a_k=\lim_{k\to\infty}\overline b_k=\infty$.
\end{Lemma}
{\em Proof of Lemma \ref{Lem:infinite_intervals}}. 
Assume that $m$ is finite and $\overline b_m=\infty$. Let $\la_s,\la_{s+1},...$, $\lim_{j\to\infty}\la_j=\infty$, be the eigenvalues of $\Delta_0$ (listed with multiplicities)
that belong to the interval $[\overline a_m,\infty)$. Clearly,
\begin{equation}\label{e:divergence}
\sum_{j\ge s}|\la_{j+1}-\la_j|=\infty\,.
\end{equation}
By Weyl's formula there exists $v>0$ such that
\begin{equation}\label{e:weyl}
\la_j=2v j^{2/n}(1+o(1))\quad\quad\mbox{as}\quad\quad j\to\infty\,.
\end{equation}
On the other hand, for any $j\ge s$, there exists $x_*\in [\la_j,\la_{j+1}]$ such that $|\la_j-x_*|\le 2cx_*^{-d}$ and $|\la_{j+1}-x_*|\le 2cx_*^{-d}$.
Hence, $\forall j\ge s$,
\begin{equation}\label{e:lambda_interval}
|\la_{j+1}-\la_j|\le 4c x_*^{-d}\le 4c \la_j^{-d}\,.
\end{equation}
By \eqref{e:weyl}, \eqref{e:lambda_interval}, and $d>n/2$ we get,
\[
\sum_{j\ge s}|\la_{j+1}-\la_j|\le C\sum_{j\ge s}\la_j^{-d}\le C\sum_{j\ge s} j^{-\frac{2d}{n}}\le C\int_s^\infty t^{-\frac{2d}{n}}\;dt<\infty\,,
\]
where $C$ stands for different positive constants. The last inequality contradicts \eqref{e:divergence}. Hence, any of the intervals $[\overline a_k,\overline b_k]$,
$1\le k\le m$, is finite. Finally, as $\lim_{j\to\infty}\la_j=\infty$ and  $\la_j\in \tilde{\mathcal I}$ we get  $m=\infty$ and 
$\lim_{k\to\infty}\overline a_k=\lim_{k\to\infty}\overline b_k=\infty$.
\finishproof

Now, for any $k\in\N$ we set $a_k = \overline a_k + \frac{3}{2}c\overline a_k^{-d}$ and $b_k = \overline b_k - \frac{3}{2}c\overline b_k^{-d}$. 
By definition, there is $\lambda,\, \lambda'\in \mbox{Spec}\left(\Delta_0 \right)$ such that $\lambda - \overline a_k= 2c\overline a_k^{-d}$ and 
$ \overline b_k - \lambda' = 2c\overline b_k^{-d}$, hence, 
 $\overline b_k - \overline a_k \ge 2c(\overline a_k^{-d} + \overline b_k^{-d})$, and we obtain  
$b_k - a_k \ge \frac{1}{2}c(\overline a_k^{-d} + \overline b_k^{-d}) >0$. 
By construction 
\begin{equation}\label{e:third_condition}
a_{k+1}-b_k > \frac{3}{2} c \overline b_k^{-d}
\end{equation}
since the intervals $[\overline a_k,\overline b_k]$ are disjoint. Hence, the intervals $[a_k, b_k]$, $k\in\N$, are disjoint as well.
Denote by ${\mathcal I}={\mathcal I}(\Delta_0)$ the union of the disjoint intervals $[a_k,b_k]$, $k\ge 1$. Note that ${\mathcal I}$ depends also on
the choice of the constant $\al\gg 1$.

\begin{Lemma}\label{lemma:natural}
For any $d>n/2$ and $c>0$ there exists $\al\gg 1$ such that the set  ${\mathcal I}\left(\Delta_0\right)$ satisfies {\rm (H$_1$)}.
In particular, the isospectral condition  \eqref{isospectral} implies  {\rm (H$_1$) - (H$_2$)} with ${\mathcal I} = {\mathcal I} \left(\Delta_0\right)$.
\end{Lemma}
{\em Proof of Lemma \ref{lemma:natural}}. 
The first condition in {\rm (H$_1$)} follows from Lemma \ref{Lem:infinite_intervals}.
Let us prove the second one. Fix $k\in\N$ and consider the interval $[\overline a_k, \overline b_k]$. Let $\lambda_s\le\cdots \le \lambda_{r}$ be the eigenvalues of $\Delta_0$ in $[\overline a_k, \overline b_k]$. Then, by \eqref{e:lambda_interval},
\begin{equation}\label{e:lambda_interval_1}
0\, \le\,  \lambda_{j+1}\, -\, \lambda_{j} \, \le\,  
4 c\lambda_{j}^{-d}\, 
\end{equation}
for $s\le j\le r$. On the other hand, by Weyl's formula \eqref{e:weyl} and Lemma \ref{Lem:infinite_intervals} there exists $k_0\in\N$ and a constant $v>0$
such that for any $k\ge k_0$
\begin{equation}\label{e:weyl_1}
v j^{2/n}\le\la_j\le 4 v j^{2/n}\,.
\end{equation}
Then by choosing $k\ge k_0$,  we get from \eqref{e:lambda_interval_1} and \eqref{e:weyl_1} that
\begin{eqnarray}
\displaystyle b_k - a_k &=& \frac{1}{2}c \overline a_k^{-d}+\sum_{s\le j\le r-1}(\la_{j+1}-\la_j)+\frac{1}{2}c \overline b_k^{-d}\le
c \overline a_k^{-d}+4c\sum_{s\le j\le r-1}\la_j^{-d}\le\nonumber\\
&=& c \overline a_k^{-d}+C \sum_{s\le j\le r-1} j^{-2d/n}\le c \overline a_k^{-d}+C \int_s^\infty t^{-2d/n}\;dt\le 
c \overline a_k^{-d}+C s^{1-\frac{2d}{n}}\le\nonumber\\
&=&c \overline a_k^{-d}+C \la_s^{\frac{n}{2}\left(1-\frac{2d}{n}\right)}\le c \overline a_k^{-d}+C \overline a_k^{\frac{n}{2}-d}\le
2C \overline a_k^{\frac{n}{2}-d}\label{e:ab-estimate}
\end{eqnarray}
where $C$ stands for different positive constants. Hence, $b_k -a_k = o(1)$ as $d>n/2$. This proves the second condition in {\rm (H$_1$)}.

We have,
\begin{equation}\label{e:b_equiv_1}
b_k = \overline b_k - \frac{3}{2}c\overline b_k^{-d}=\overline b_k(1+o(1))\quad\quad\mbox{as}\quad\quad k\to\infty.
\end{equation}
Combining this with \eqref{e:third_condition} and taking $k_0\ge 1$ sufficiently large we see that for any $k\ge k_0$,
\[
a_{k+1}-b_k > c b_k^{-d}\,.
\]
Hence, the set $\mathcal I$ satisfies {\rm (H$_1$)} for $\al\ge a_{k_0}$. This completes the proof of the first statement of the lemma.
It is clear from  construction that $\mathcal I$ satisfies {\rm (H$_2$)} for any $a\ge\al$.
\finishproof

\begin{Remark}\label{Rem:natural}
It follows from Lemma \ref{lemma:natural} and inequality \eqref{e:ab-estimate} that for any $d>n/2$ and $c>0$ the isospectral condition (\ref{isospectral}) implies $(\mbox{H}_1)_s$-(H$_2$)
with ${\mathcal I}={\mathcal I}(\De_0)$ and $0\le s<2d-n$.
\end{Remark}
Indeed, as in \eqref{e:b_equiv_1} we get
\begin{equation}\label{e:a_equiv_1}
a_k = \overline a_k + \frac{3}{2}c\overline a_k^{-d}=\overline a_k(1+o(1))\quad\quad\mbox{as}\quad\quad k\to\infty.
\end{equation}
Then,  taking $k_0\gg 1$  we get form \eqref{e:ab-estimate} that 
\[
\displaystyle 0<b_k - a_k\le 2C a_k^{\frac{n}{2}-d}\le c a_k^{s/2}\,.
\]
We are going to give now a formal  argument which allows one to get spectral invariants from continuous quasimodes.
\begin{Def}\label{def:continuous-quasimode}
Let $A_t$, $t\in [0,1]$, be a family of selfadjoint operators in a Hilbert space ${\mathcal H}$. By a  continuous  quasimode  of order $M\ge 1$ of $A_t$ we mean  a sequence $(\mu_q(t)^2,u_{q}(t))$, $q\in \N$, such that 
\begin{itemize}
\item   
$u_{q}(t)$ belongs to the domain of definition $D(A_t)$ of $A_t$ and $\|u_{q}(t)\|=1$,  
\item there is  a constant $C_M>0$ independent of $q$ and of $t$  such that
\begin{equation}
\left\|\left(A_t  - \mu_q(t)^2 \right)u_{q}(t)\right\| \le C_M\,  \mu_q(t)^{-M}
\label{discrepancy}
\end{equation}
for any $q\in \N$ and $t\in [0,1]$, 
\item $\mu_q(t)\in \R$ and 
\[
\mu_q(t) = \mu_q^0 +c_{q,0}(t) + c_{q,1}(t) (\mu_q^0)^{-1}  + \cdots + 
c_{q,M}(t) (\mu_q^0)^{-M} ,
\]
where
\begin{enumerate}
\item[(i)]
$\mu_q^0$  are independent of $t$ and $\lim\mu_q^0= +\infty$  as $q\to\infty$, 
\item [(ii)]
$\{c_{q,j}\in C([0,1],\R):\, q\in \N,\, 0\le j\le M\}$ is a bounded family of  continuous functions. 
\end{enumerate}
\end{itemize}
\end{Def}
We have the following
\begin{Lemma}\label{lemma:invariants} Suppose that the family of selfadjoint operators  $A_t$, $t\in [0,1]$,   in the Hilbert space  ${\mathcal H}$ is isospectral in the sense of $(\mbox{H}_1)_s$-$(\mbox{H}_2)$, where $s\ge 0$ is an integer. Let  $(\mu_q(t)^2,u_{q}(t))$, $q\in \N$, be a continuous family of quasimodes of order
$M > \max\{2d,s\}$  such that 
\begin{enumerate}
\item[{\rm (iii)}]$c_{q,j}$ is independent of $t$ for $j\le s$, 
\item [{\rm (iv)}] $c_{q,s+1}(t) = c'_{q,s+1} +  c(t)$, where $c'_{q,s+1}$ does not depend on $t$ and $c(t)$ does not depend on $q$. 
\end{enumerate}
Then $ c(t)=c(0)$ for any  $t\in [0,1]$. 
\end{Lemma}
{\em Proof.}
It is easy to see that for any $q\in \N$ and $t\in [0,1]$, the distance from $\mu_q(t)^2$ to the spectrum of $A_t$ can be estimated above by
\[
d_{t,q}:=\left|\, \mbox{Spec}\left(A_t \right)\, -\,  
\mu_q(t)^2\right|\,  \le\,  C_M\, \mu_q(t)^{-M}\,   .
\]
Indeed, if $d_{t,q}\neq 0$  the spectral theorem and \eqref{discrepancy} yield
\[
\frac{1}{d_{t,q}} \ge \|(A_t - \mu_q(t)^2)^{-1}\|  \ge 
\|(A_t - \mu_q(t)^2)u_q(t)\|^{-1} \ge \frac{\mu_q(t)^M}{C_M}. 
\]
For any $q\in\N$ and  $ t\in [0,1]$ let us take $\lambda_{t,q}\in \mbox{Spec}\left(A_t \right)$ such that
\[
\left|\, \lambda_{t,q} \, -\,  
\mu_q(t)^2\right|\, \le\, 
C_M \mu_{q}(t)^{-M}\, . 
\]
The properties (i) and (ii) of the quasimode imply that  $\lim\mu_q(t)= \infty$ as $q\to\infty$ uniformly with respect to $t\in [0,1]$. 
There is $q_0\gg 1$ such that 
\begin{equation}
\forall\, q\ge q_0,\ \forall\,  t\in [0,1]\, :\quad \frac{1}{4} (\mu_q^0)^2 \le \frac{1}{2} \mu_q(t)^2 \le\lambda_{t,q} \le 2\mu_q(t)^2 \le  4 (\mu_q^0)^2  ,
\label{quasi-inequalities}
\end{equation}
and  for any $q\ge q_0$ and $t\in [0,1]$ we get
\begin{equation}
\left|\, \lambda_{t,q} \, -\,  
\mu_q(t)^2\right|\, \le\, 
C' \lambda_{t,q}^{-M/2}\,   
\label{quasi-to-spectrum}
\end{equation}
where $C'=2^{M/2}C_M $. Moreover,  \eqref{quasi-inequalities} and (i) imply that  
\begin{equation}
\mbox{$\lim \lambda_{t,q} = \infty$  as  $q\to\infty$ uniformly with respect to $t\in [0,1]$.} 
\label{uniform-limit-lambda}
\end{equation}
In particular,   
there is   $q_0\ge 1$ such that $\lambda_{t,q}\ge a$  for any $q\ge q_0$ and   $t\in [0,1]$, where $a\ge 1$ is the constant in (H$_2$). Now  using (H$_2$) we get 
  for any $q\ge q_0$ and   $t\in [0,1]$   an integer $k=k(q,t)\ge 1$ such that 
\begin{equation}\label{e:la_in}
\lambda_{t,q}\in [a_k,b_k].
\end{equation}
Fix $\beta$ so that $M >\beta > \max\{2d,s\}$. Then choosing  $q_0\gg 1$ we obtain from \eqref{quasi-to-spectrum}, \eqref{uniform-limit-lambda} and \eqref{e:la_in}   that for any $q\ge q_0$ and  $t\in [0,1]$ the quasi-eigenvalue $\mu_q(t)^2$ belongs to  the interval
$[a_k-\frac{c}{2} \la_{q,t}^{-\beta/2},b_k+ \frac{c}{2} \la_{q,t}^{-\beta/2}]$, where  $c>0$ is the constant of the third assumption of  $(\mbox{H}_1)$. Since $\lambda_{t,q} > a_k$, we get
\begin{equation}\label{e:the_intervals}
\forall\, q\ge q_0, \ \forall\, t\in [0,1]\, :\quad \mu_q(t)^2\in I_k:=\left[a_k-\frac{c}{2} a_k^{-\beta/2},b_k+ \frac{c}{2} a_k^{-\beta/2}\right]\, , 
\end{equation}
where $ k=k(q,t)$.
In particular, 
\[
b_{k(q,t)}\ge  \mu_q(t)^2 - \frac{c}{2} a_{k(q,t)}^{-\beta/2} \ge \frac{1}{2}(\mu_q^0)^2 - \frac{c}{2} a_{1}^{-\beta/2},
\]
 which implies that $k(q,t) = \infty$  as $q\to\infty$ uniformly with respect to $t\in [0,1]$. 
On the other hand, using the third assumption of $(\mbox{H}_1)$, the relation $b_k=a_k(1+o(1))$ as $k\to\infty$, which  follows from the first two assumptions in {(\rm H$_1$)}, 
and the inequality $\beta>2d$, we get 
\[
(a_{k+1}-\frac{c}{2}a_{k+1}^{-\be/2})-(b_k+\frac{c}{2}a_k^{-\be/2})=(a_{k+1}-b_k)-\frac{c}{2}a_{k+1}^{-\be/2}-\frac{c}{2}a_k^{-\be/2}\ge
c b_k^{-d}-c a_k^{-\be/2}>0
\]
for any $k\ge k_0$  and $k_0\gg 1$. This shows that the intervals  $I_k$ in \eqref{e:the_intervals} do not intersect each other for  $k\ge k_0$.
Choose $q_0\gg 1$ so that $k(q,t)\ge k_0$ for any $q\ge q_0$ and $t\in [0,1]$. 
The function  $\mu_q(t)^2$ is continuous on $[0,1]$, hence,  it can not jump from one interval to another for $q\ge q_0$. Consequently,    $k(q,t)$ does not depend on $t$ for $q\ge q_0$. We have proved that for any $q\ge q_0$ there is $k=k(q)\in\N$ independent
of $t$ such that
\begin{equation}\label{e:mu_in}
\forall\, t\in [0,1]\, : \quad \mu_q(t)^2 \in \left[a_k-\frac{c}{2}a_k^{-\beta/2},b_k+ \frac{c}{2} a_k^{-\beta/2}\right]\,.
\end{equation}
Moreover,  $k(q)\to\infty$ as $q\to\infty$  and we obtain  
\[
\mu_q(t)^2=a_{k(q)}+O(b_{k(q)}-a_{k(q)}) + O(a_{k(q)}^{-\be/2})=a_{k(q)}(1+o(1)),\quad \mbox{as}\ q\to\infty\,,
\]
uniformly in $t\in[0,1]$. On the other hand  {\rm (i)} and {\rm (ii)} imply
\[
\mu_q(t)=\mu_q^0(1+o(1)),\quad \mbox{as}\  q\to\infty\, ,
\]
uniformly in $t\in[0,1]$. 
 The last two formulas combined with \eqref{e:mu_in} show that for any $q\ge q_0\gg 1$ 
$$
(\mu_q^0)^{s+1}|\mu_q(t) - \mu_q(0)|\le C\mu_q(t)^{s}|\mu_q(t)^2 - \mu_q(0)^2|\le C \left(a_{k(q)}^{s/2}\left(b_{k(q)} -a_{k(q)}\right) + ca_{k(q)}^{(s-\beta)/2}\right):=\epsilon_q
$$
where $C>0$ stands for different constants which are independent of $q$ and of $t$.  The inequality $\beta>s$ and \eqref{strong-iso} imply that
$\epsilon_q \to 0$ as $q\to\infty$. Moreover,  (iii) and (iv) yield 
\[
\mu_q(t)-\mu_q(0)=(c(t)-c(0))(\mu_q^0)^{-(s+1)}+O((\mu_q^0)^{-(s+2)}) \quad \mbox{as}\  q\to\infty\, .
\] 
Thus
\begin{equation}
\label{mu}
|c(t)-c(0)| \le (\mu_q^0)^{s+1}|\mu_q(t) - \mu_q(0)| + O((\mu_q^0)^{-1})= \epsilon_q + O((\mu_q^0)^{-1}) \to 0
\end{equation}
as $q\to\infty$, which proves the claim since $c(t)$ does not depend on $q$. 
\finishproof

\noindent
We point out that there is no restriction on the nature of the spectrum of the operators $A_t$. An analogue of the lemma can be proved in the case of resonances close to the real axis replacing the intervals in condition $(\mbox{H}_1)_s$-$(\mbox{H}_2)$ by boxes $[b_k,c_k] + i[0,d_k]$, where $d_k>0$, $\lim d_k=0$, and $a_k<b_k$ satisfy  $(\mbox{H}_1)_s$,  and using  results of Stefanov \cite{Ste}, and  Tang and Zworski \cite{TZ} about the localization of resonances near quasi-eigenvalues. 
%Lemma \ref{lemma:invariants} will be useful to obtain  invariants for continuous isospectral deformations (see Sect. \ref{subsec:FurtherRemarks}).

\subsection{\bf Billiard ball map}\label{subsec:billiard-ball}

We recall from Birkhoff \cite{Birkhoff} (see also \cite{Tabach}) the definition of the 
billiard ball map $B$ associated to a billiard table
$(X,g)$ with a smooth boundary $\Gamma$. We are interested here only in the ``broken geodesic flow'' given by the elastic reflection of geodesics hitting transversally the boundary. It induces a discrete dynamical system at the boundary which can be described as follows. 

Denote by $h$  the Hamiltonian on $T^\ast X$ 
corresponding to the Riemannian metric $g$ on $X$ via the Legendre
transformation and by $h_0$ the Hamiltonian on $T^\ast\Gamma$ corresponding to the induced
Riemannian metric on $\Gamma$. 
The billiard ball map $B$ lives in an open  subset of  the open 
coball bundle 
\[
B^\ast \Gamma = \{(x,\xi)\in T^\ast \Gamma:\, h_0(x,\xi) < 1\}\,  
\] 
and it is
defined as follows. Let $\overline{B^\ast \Gamma}$ be the closure  of $B^\ast\Gamma$ in $T^\ast \Gamma$. Denote by   $S^\ast X$  the cosphere bundle of $X$ which consists of all $(x,\xi)\in T^\ast X$ such that $h(x,\xi) =  1$,  and set
\[
\Sigma = S^\ast X|_{\Gamma}:= \{(x,\xi)\in S^\ast X:\, x\in
\Gamma\}\ \mbox{and}\ \Sigma^{\pm}:=  \{(x,\xi)\in \Sigma :\,
\pm \langle \xi,\nu(x)\rangle >0\}\, .
\] 
Consider the natural projection $\pi_\Sigma:\, 
\Sigma \, \rightarrow \, \overline{B^\ast \Gamma}$ assigning to each $(x,\eta)\in
\Sigma$ the covector $(x,\eta|_{T_x\Gamma})$. Its restriction to $B^\ast \Gamma$ admits two smooth inverses 
\begin{equation}\label{outgoing-vector}
\pi_\Sigma^{\pm}:\,  B ^\ast \Gamma\,  \rightarrow\, 
\Sigma^{\pm}\, ,\ \pi_\Sigma^{\pm}(x,\xi) = (x,\xi^\pm) \, ,
\end{equation}
which  can be extended continuously on the  closed coball bundle  $ \overline{B ^\ast \Gamma}$.  The covector $\xi^+$ is called  ``outgoing'' and  $\xi^-$ ``incoming''. 
Take $(x,\xi)$ in  ${B} ^\ast \Gamma$  and  
consider the integral  
curve  $\exp(tX_{h})(x,\xi^+)$  of the Hamiltonian vector field   
$X_{h}$ starting at $(x,\xi^+)\in \Sigma^+$. 
If it intersects transversally    
$\Sigma$ at a time $t_1>0$ and lies entirely in the   
interior  of $S^\ast  X$   
for $t\in (0,t_1)$, we set   
\[
(y,\eta^-)=J(x,\xi^+) =\exp(t_1X_{h})(x,\xi_+)\in \Sigma^-\, ,
\] 
and define $B(x,\xi):=(y,\eta)$,  
where $\eta := \eta_-|_{T_y \Gamma}$. Notice that by construction $(y,\eta)\in B^\ast\Gamma$.       
We  denote by $\widetilde B^\ast \Gamma$ the set  
of all such points $(x,\xi)$. 
In this way we obtain a smooth exact symplectic map   
$B:\widetilde B^\ast \Gamma \rightarrow B^\ast \Gamma $, given by 
$B=\pi_\Sigma\circ J \circ \pi_\Sigma^+$.

\subsection{ Quasimodes}

Fix a positive  integer $M$ and let  $\mathcal M\subset \Z^n$ be  an infinite  set. Recall that the domain of definition of $\Delta_{g,K}$ is given by \eqref{domain}.  Then by a quasimode  $\cal Q$ of
$\Delta_{g,K}$ of order $M$ and with a index set $\mathcal M$ we mean   an infinite
sequence   $(\mu_q,u_q)_{q\in \mathcal M}$, such that  $\mu_q$ are
positive, $\lim \mu_q = +\infty$,   
$u_q$ belongs to the Sobolev space  $H^2(X)$,   $ \|u_q\|_{L^2(X)} =1$,  and there is a constant $C_M>0$ such that
\begin{equation}  
\left\{  
\begin{array}{lcr}  
\left\|\Delta\,  u_q\ - \ \mu_q^2\,  u_q\right\|\ \le \   
C_M\,  \mu_q^{-M}  \, \quad \mbox{in}\ L^2(X)\, , \\ [0.3cm] 
\displaystyle \partial u_q/\partial \nu|_\Gamma \  - \ K\,
u_q|_\Gamma  =\ 0 \, .  
\end{array}  
\right.  
                                         \label{thequasimode}
\end{equation}  
Denote by $A(\varrho)$ the action along the  broken bicharacteristic
starting at $\varrho\in \Lambda$ and with endpoint  $P(\varrho)\in \Lambda$. 
Note that $2A(\varrho)>0$ is just the length of the corresponding geodesic arc. Denote by 
$\pi_\Gamma:
T^\ast \Gamma \to \Gamma$  the inclusion map.  Recall that for any  $(x,\xi) \in B ^\ast  \Gamma$ the angle  $\theta = \theta(x,\xi)\in (0,\pi/2]$ is  determined by   $\sin \theta =\langle \xi^+,\nu(x)\rangle$, where  $\langle \cdot,\cdot\rangle$ stands for the pairing between vectors and covectors.
\begin{Theorem}\label{quasimodes}
Let $\Lambda$ be a Kronecker torus satisfying (H$_3$).  
Fix a positive integer  $M\ge 2$ and $\ell > M(\tau +2)+2n+ (n-1)/2$, where $\tau > n-1$ is the exponent in (\ref{sdc}). 
Then for any $K\in C^\ell(\Gamma,\R)$ there is a quasimode $(\mu_q(K),u_q(K))_{q\in \mathcal M}$ 
  of
$\Delta_{g,K}$ of order $M$ satisfying (\ref{thequasimode})  
and with an infinite index set $\mathcal M \subset \Z^{n}$ independent of $K$, such that 
$$
\mu_q = \mu_q^0 +c_{q,0} + c_{q,1} (\mu_q^0)^{-1}  + \cdots + 
c_{q,M} (\mu_q^0)^{-M}   
$$
where 
\begin{enumerate}
\item[(i)]
$\mu_q^0$  is independent of $K$ and  there is $C^0>0$ such that $\mu_q^0 \ge C^0|q|$ for any 
$q\in \mathcal M$,
\item [(ii)]
the function  $c_{q,j}:C^\ell(\Gamma,\R)\to \R$ assigning to each 
  $ K\in C^\ell(\Gamma,\R)$ the corresponding coefficient $c_{q,j}(K)$ of $\mu_q(K)$ is continuous.
\item[(iii)]  the coefficients $c_{q,j}(K)$ and the positive constant $C_M(K)$ in  (\ref{thequasimode}) are uniformly bounded on any bounded subset ${\mathcal B}$  of $C^\ell(\Gamma,\R)$, which means that 
   there is $C= C({\mathcal B})>0$ 
such that $C_M(K) \le C$ and $|c_{q,j}(K)| \le C$ for any  $q\in {\mathcal M}$, $0\le j\le M$,  and any $K\in {\cal B}$,
\item [(iv)]
$c_{q,0}$ is independent of $K$ and  
\[
c_{q,1}(K)=c'_{q,1} + c''_{1} \sum_{j=0}^{m-1}\int_{\Lambda_j}
\frac{K\circ \pi_\Gamma}{\sin \theta}\, d\mu_j \, , 
\] 
where $c'_{q,1}$ is  independent of $K$, and 
$ c''_{1} =     \displaystyle{\frac{2 }{ 
\int_{\Lambda} A(\varrho)\,  d\mu }. }$
\end{enumerate}
\end{Theorem}
To prove  Theorem \ref{main} we denote by 
 ${\mathcal B}$ the set  $\{K_t\, :\, t\in [0,1]\}$. 
Take  $M = [2d] +1$, the smallest positive integer bigger than $2d$,  and 
consider the continuous quasimode  $(\mu_q(K_t)^2, u_q(K_t))$, $t\in [0,1]$, given by
Theorem \ref{quasimodes}.  Then  apply Lemma \ref{lemma:invariants} with $s=0$.

\section{Construction of continuous quasimodes}\label{sec:construction}
The proof of Theorem \ref{quasimodes}   is quite long and we divide it  in several steps as follows. In Sect. \ref{subsec:reduction} we reduce the problem to the boundary and introduce the corresponding microlocal monodromy operator $W(\lambda)$. This is a $\lambda$-FIO  whose canonical relation is just the graph of $P$. Our goal is to ''separate the variables''
microlocally near $\Lambda$  and then to obtain  quasimodes  of $W(\lambda)$  concentrated at $\Lambda$.  For this reason we first find a symplectic Birkhoff normal form of $P$ in 
Sect. \ref{subsec:BNF}. Then conjugating  $W(\lambda)$ with a suitable  $\lambda$-FIO we get a $\lambda$-FIO  $W_1(\lambda)$ with a simple phase function in Sect. \ref{subsec:quantization}. In Sect. \ref{subsec:QBNF} we obtain a microlocal  Birkhoff normal form  $W^0(\lambda)$ of $W_1(\lambda)$ by conjugating it  with a suitable  $\lambda$-PDO  and  solving at any step the corresponding homological equation.  
In this way we separate microlocally the variables near $\Lambda$, which means that the amplitude of $W^0(\lambda)$ does not depend on the angular variables but only on the action variables.  A microlocal ''spectral decomposition''   of $W(\lambda)$ near $\Lambda$ is obtained in  Sect. \ref{subsec:spectral}. In Sect.  \ref{subsec:construction} we get quasi-eigenvalues $z_q(\lambda)$ of $W(\lambda)$, where $q$ belongs to a unbounded subset of $\Z^n$, and then we solve  the equations $z_q(\lambda)=1$, which gives  the desired quasimodes of the problem. 

\subsection{ Reduction to the boundary.}\label{subsec:reduction}
\setcounter{equation}{0}
We are going to give a brief idea of the reduction to the boundary, which is a variant of the reflection method for the wave equation.  Consider an ``outgoing'' solution of 
 the Helmlotz equation $(\Delta -\lambda^2)u = O_N(|\lambda|^{-N})$ at high frequencies ($|\lambda|\to \infty$) with ``initial data''   $\Psi_0(\lambda) f$ at $\Gamma$ microlocalized near $\Lambda_0=\Lambda$, where $\Psi_0(\lambda)$ is a  $\lambda$-PDO having a frequency set (semi-classical wave front set) in a neighborhood of  $\Lambda$. The solution $u$ is given by   $u= H_0(\lambda)f$,  where $H_0(\lambda)$  is a $\lambda$-FIO.  The construction of the operator $H_0(\lambda)$ is provided in Appendix A.1. 
The restriction of $u$ at $\Gamma$ is given modulo  $O_N(|\lambda|^{-N})$ by $u|_{\Gamma} =  \Psi_0(\lambda) f + G_0(\lambda)f$, where $G_0(\lambda)$ is 
a $\lambda$-FIO of order 0 and its canonical relation is the graph of the billiard ball map near $\Lambda_0$ (see Lemma \ref{Lemma:restriction} and \eqref{boundary-trace3}). To satisfy the boundary conditions near $\Lambda_1=B(\Lambda_0)$ in the case when $m\ge 2$ we use the reflexion method. Namely, we add to $H_0(\lambda)$ an operator  $\tilde{G}_1(\lambda) = H_1(\lambda)Q_1(\lambda)G_0(\lambda)$, where $H_1(\lambda)$ is an ``outgoing'' parametrix of the Helmlotz equation in $X$ with Dirichlet boundary conditions  microlocalized near $\Lambda_1$ and $Q_1(\lambda)$ is a suitable $\lambda$-PDO supported near $\Lambda_1$, and then consider $u= H_0(\lambda)f + \tilde{G}_1(\lambda)  f$. In the case of  Dirichlet boundary conditions one can take $Q_1(\lambda)$ to be $-{\rm Id}$ microlocally near $\Lambda_1$. In the case of Robin boundary conditions one uses the observation  that the trace of the normal derivative of  $H_0(\lambda)f$ at $\Gamma$ is given by the action of a suitable $\lambda$-PDO of order one  on the trace of $H_0(\lambda)f$  at $\Gamma$. The same is valid for the trace of the normal derivative of $\tilde{G}_1(\lambda)$  at $\Gamma$. In this way one reduces  the Robin boundary condition  microlocally near $\Lambda_1$ to  an equation with respect to $Q_1(\lambda)$ (see \eqref{pdo-equation1}). The solution of this equation is given by \eqref{pde-solution} and \eqref{principal-symbols}. 
Similarly one can treat  the boundary conditions near $\Lambda_j=B^j(\Lambda)$ for 
any $0<j<m$ which leads to a solution $G(\lambda)f= H_0(\lambda)f + \tilde{G}_1(\lambda)  f +\cdots + \tilde{G}_{m-1}(\lambda)  f$ of the Helmholtz equation modulo  $O_N(|\lambda|^{-N})$   satisfying the boundary conditions microlocally in a neighborhood of $\Lambda_j$ for 
any $0<j<m$. To satisfy the boundary conditions in a neighborhood of $\Lambda_m=\Lambda_0$ we obtain in the same way an equation of the form 
$(W(\lambda) -\mbox{Id})f=0$, where $W(\lambda)$ is the monodromy operator \eqref{monodromy}. The construction of the parametrix $G(\lambda)$ is similar to that in \cite{C-P} and it is close  that of the mixed problem for the wave equation which has been done by Chazarain \cite{Ch} and by Guillemin and Melrose \cite{G-M2} 
(see also \cite{P-S} and \cite{S-V}).

We proceed now with the construction of the parametrix. 
If $m\ge 2$, without loss of generality we can suppose that $\Lambda_j=B^j(\Lambda) \neq \Lambda$ for $1\le j<m$. Indeed, 
 $P$ acts transitively on each $\Lambda_j$ since
$\omega$ is Diophantine,
hence,  $\Lambda_i\cap \Lambda_j = \emptyset$ if  $0<|i-j|<m$ and $m\ge 2$. 
Choose neighborhoods 
$U_j\subset  \displaystyle \widetilde{B ^\ast }\Gamma$ of
$\Lambda_j$, $0\le j \le m$, such that  $U_{j+1}$ is a  neighborhood of 
the closure of $B(U_j)$ for $j= 0,\ldots, m-1$ and $U_{0}\subset U_m$. Moreover, if $m\ge 2$ we assume that 
$U_i\cap U_j = \emptyset$ for  $0<|i-j|<m$. Consider a $C^\infty$
extension  $(\widetilde X, \widetilde g)$  of $(X,g)$ across $\Gamma$ 
such that  any integral curve $\gamma$ of the
Hamiltonian vector field $X_{\widetilde h}$ ($\widetilde h$ being the  Hamiltonian corresponding to $\widetilde g$) starting at
$\pi^+_\Sigma(U_j)$, $j= 0,\ldots, m-1$, satisfies 
\begin{equation}
\gamma \cap T^\ast
\widetilde X |_\Gamma \subset \pi^+_\Sigma(U_j)\cup
\pi^-_\Sigma(U_{j+1}) .
                          \label{bicharacteristic}
\end{equation}
Then $\gamma$ intersects transversally $T^\ast X|_{\Gamma}$ and  for each $\varrho \in U_j$ there is a unique 
$T_j(\varrho) >0$ such that 
\[
\exp(T_j(\varrho)X_{\widetilde h})(\pi^+_\Sigma(\varrho))\in \pi^-_\Sigma(B(U_j))\, .
\]
Let $\psi_j(\lambda)$, $j=0,1, \ldots, m$, 
be   classical $\lambda$-pseudodifferential
operators ($\lambda$-PDOs) of order $0$ 
on $\Gamma$ with a large parameter $\lambda$  and
compactly supported amplitudes in $U_j$  such that
$$
\mbox{WF}'(\mbox{Id} - \psi_j) \cap \Lambda_j =
\emptyset \,  , 
$$
and
\begin{equation}
\mbox{WF}'(\psi_{j+1}) \subset B(U_j)\, , \ 
\mbox{WF}'(\mbox{Id} - \psi_{j+1}) \cap 
B(\mbox{WF}'(\psi_{j})) = \emptyset \ \mbox{for}\ j=0,\ldots m-1\, .
                                    \label{psi-pdos}
\end{equation} 
Hereafter 
$\mbox{WF}'(\psi_{j})$ stands for the frequency set of
$\psi_{j}$  (the semi-classical $\hbar$-wave front with $\hbar = 1/\lambda$ \cite{D-S}, \cite{E-Z}, \cite{PP}), and by a  
  ``classical'' $\lambda$-PDO  we mean  that in any local coordinates 
the corresponding distribution kernel is of the form (\ref{pdo})
where the amplitude has an asymptotic expansion $q(x,\xi,\lambda) \sim \sum_{k=0}^\infty q_k(x,\xi)\lambda^{-k}$
and $q_k$ are $C^\infty$ smooth and uniformly compactly supported (see Appendix A.1). In particular the distribution 
kernel ${\rm OP}_\lambda(q)(\cdot,\cdot)$ is smooth for each $\lambda$ fixed.
 We take $\lambda$ in a complex strip 
\begin{equation}
\label{strip}
{\cal D}\,  :=\,  \{z \in {\C}:\,  |{\rm Im}\, z| \le D_0,\,  {\rm Re}\, z
\ge 1\}\, ,
\end{equation}
$D_0>0$ being  fixed. 

We are looking for a  microlocal outgoing parametrix  
$H_j(\lambda):L^2(\Gamma)\rightarrow  C^\infty( \widetilde X )$, 
of the Dirichlet problem for the Helmholtz equation with ``initial data'' concentrated 
in $U_j$ such that 
\begin{equation}
\forall N\in \N\, , \quad (\Delta-\lambda^2)H_j(\lambda) = O_N(|\lambda|^{-N})
\label{parametrix}
\end{equation}
in a neighborhood of $X$ in 
$\widetilde X$. 
Hereafter,
\[
O_N(|\lambda|^{-N})\, :\ L^2(\Gamma)
\longrightarrow  L^2 (\widetilde X)
\] 
stands for any family  of continuous  operators depending
on $\lambda$ with norms  
$\le C_{N}(1 + |\lambda|)^{-N}$, 
$C_{N}>0$.  
We will also denote  by 
$O_N(|\lambda|^{-N})\, :\ L^2(\Gamma)
\longrightarrow  L^2(\Gamma)$ 
any family of continuous  operators depending
on $\lambda$ with norms uniformly bounded by $C_{N}(1 + |\lambda|)^{-N}$, 
$C_{N}>0$. We take $N=M$ to be the order of the quasimode we are going to construct. 

The operator  $H_j(\lambda)$  
is a Fourier integral operator of order $-1/4$  with a large parameter
$\lambda\in {\cal D}$ ($\lambda$-FIO) the distribution kernel of which is an
oscillatory integral in the sense of Duistermaat \cite{Du} and it is constructed in Appendix A.1. 
In any local coordinates its amplitude is $C^\infty$ smooth, 
uniformly compactly supported for $\lambda\in {\cal D}$ and  it 
has an asymptotic expansion in powers of $\lambda$ up to any negative
order. In particular, $H_j(\lambda)u$ is a $C^\infty$ smooth function for any fixed $\lambda$ 
and $u\in L^2(\Gamma)$. 
The corresponding canonical relation lies in $T^\ast \widetilde X \times T^\ast
\Gamma $ and it  is given by 
$$
{\mathcal C}_j \ :=\ \left\{\left(  
\exp(s X_ {\widetilde h})(\pi_\Sigma ^+(\varrho)); \varrho\right)\,
:\ \varrho\in U_j\, ,\ -\varepsilon <s < T_j +\varepsilon \right\}\, ,\ \varepsilon >0\, .
$$
Consider the operator of restriction 
$\imath_\Gamma^\ast: C^\infty(\widetilde X) \to  C^\infty(\Gamma)$, $\imath_\Gamma^\ast(u)=u_{|\Gamma}$,    as a $\lambda$-FIO 
of order $1/4$, the canonical relation ${\cal R}$ of which is just the inverse  of 
the canonical relation given by the  conormal bundle of the graph of
the inclusion map $\imath_\Gamma :\Gamma \to \widetilde X$. The composition 
${\cal R}\circ {\mathcal C}_j$ is transversal for any $j$ and it is a disjoint union of
the diagonal in $U_j\times U_j$ (for $s=0$) and of the graph of the billiard ball map
$B:U_j\to U_{j+1}$ (for $s=T_j$) (see Lemma \ref{Lemma:restriction} ). 
Let $\Psi_j(\lambda)$, $1\le j\le m-1$,  be a $\lambda$-PDO of order $0$ such that 
\[
\mbox{WF}'(\Psi_j)\subset U_j \ \mbox{and}\ 
\mbox{WF}'(\Psi_j-\mbox{Id})\cap \mbox{WF}'(\psi_j) =\emptyset. 
\] 
Taking $\Psi_j(\lambda)$  as initial data at $\Gamma$ for $s=0$ 
 in  Lemma \ref{Lemma:restriction} we obtain an operator $H_j(\lambda)$
satisfying (\ref{parametrix}). Moreover, by \eqref{boundary-trace3} and \eqref{boundary-trace4} we have 
\begin{equation}
\imath_\Gamma^\ast H_j(\lambda) = \Psi_j(\lambda) + 
G_j(\lambda)
+ O_M(|\lambda|^{-M})
\, ,  
                     \label{boundary-trace}
\end{equation}
where $G_j(\lambda)= E(\lambda)^{-1}G_j^0(\lambda)E(\lambda)$,  $E(\lambda)$ is a $\lambda$-PDO of order $0$ with a frequency set contained in $\widetilde B^\ast \Gamma$ and which is elliptic in a neighborhood of 
$\cup_{j=0}^m \overline U_j$. Moreover,  $G_j^0(\lambda)$ 
is a $\lambda$-FIO of order $0$ and its 
 canonical relation 
is the graph  of the 
billiard ball map $B: U_j \to U_{j+1}$.   Parameterizing the graph of $B:U_j\to U_{j+1}$ by $\varrho\in U_j$,    
the principal symbol of  $G_j^0(\lambda)$ becomes 
\begin{equation}
\sigma(G_j^0(\lambda))(\varrho)=\exp\left(i\pi\vartheta_j/4\right)\exp (i\lambda A_j(\varrho)) 
                     \label{principal-symbol-G}
\end{equation}
in a neighborhood of   
$\mbox{WF}'(\psi_j)$, where $\vartheta_j\in \Z$ is a  Maslov's index
and    $A_j(\varrho) =
\int_{\gamma_j(\varrho)}\xi dx $ 
is the action along the integral curve $\gamma_j(\varrho)$ 
of the Hamiltonian vector
field $X_{\widetilde h}$ 
starting at $\varrho\in U_j$ and with endpoint $B(\varrho)\in
U_{j+1}$ (see \eqref{boundary-trace4}). 
Hereafter, to simplify the notations we trivialize the corresponding half-density bundles using for example the  volume form on $\widetilde X$. 
In particular,  the frequency 
set $WF'$  of
$G_j(\lambda)$  
is contained in 
$U_{j+1}\times U_{j}$ for any $j=0,\dots,m-1$.  
Note that 
$2A_j(\varrho)$ is just the length $T_j(\varrho)$ of the corresponding geodesic
$\widetilde\gamma_j(\varrho)$ in $X$ and we have   
$$
\pi_\Sigma \left(\exp(2A_j(\varrho)X_{\widetilde h})(\pi_\Sigma^+(\varrho))\right) 
= B(\varrho)\, 
,\  \varrho\in U_j\, .
$$
Fix a bounded set ${\mathcal B}$ in $C^\ell(\Gamma, \R)$ and take 
 $K\in {\mathcal B}$. Consider the operator 
${\mathcal N} = \partial/\partial {\widetilde \nu} - {\widetilde K}$ in a
neighborhood of $\Gamma$ in $\widetilde X$,
where $\widetilde \nu$ is a normal vector
field to $\Gamma$ and $ \widetilde K$ is a $C^\ell$-smooth extension of $K$
with compact support contained in a small neighborhood of $\Gamma$. 
To construct $\widetilde K$ we extend $K$ as a
constant on the integral curves of $\widetilde \nu$ and then multiply
it with a suitable  cut-off function. In this way we obtain a
continuous  map $K \to \widetilde K$  from 
$C^\ell(\Gamma, \R)$ to $C_0^\ell(\widetilde X, \R)$. 

Suppose first that $m=1$. Then $\Lambda = B(\Lambda)$ and  $U_0\cup B(U_0)\subset U_1$. 
Set $G(\lambda)= H_0(\lambda)\psi_0(\lambda)$. We have 
 $(\Delta-\lambda^2)G(\lambda) = O_M(|\lambda|^{-M})$
in a neighborhood of $X$ in 
$\widetilde X$, in view of (\ref{parametrix}).  
To satisfy the boundary conditions we should have $\imath_\Gamma^\ast\, \,   {\mathcal N}\,
G(\lambda) =O_{M}(|\lambda|^{-M})$.
Using the symbolic calculus, (\ref{psi-pdos}) and  \eqref{boundary-trace}
we obtain
\[
\imath_\Gamma^\ast\, \,   {\mathcal N}\,
G(\lambda) = \psi_1(\lambda)(\lambda R_{0}^+ - K)\, \psi_0(\lambda)  + 
\psi_1(\lambda)(\lambda R_{1}^- - K)\, 
G_0(\lambda)\psi_0(\lambda)
 +O_{M}(|\lambda|^{-M})\, .
\]
Here, $ R_{0}^+ (\lambda)$ is a classical $\lambda$-PDO  of order $0$  on $\Gamma$ 
independent of $K$, with a $C_0^\infty$-symbol in any local coordinates, and with 
principal symbol 
\[
\sigma(R_{0}^+)(\varrho)\, =\, i\sqrt{1-h_0(\varrho)}\, , \ 
\varrho \in U_{0}\, , 
\] 
and 
$ R_{1}^- $ is a classical $\lambda$-PDO of order $0$  on $\Gamma$ 
independent of $K$  with 
principal symbol 
\[
\sigma(R_{1}^-)(\varrho)\, =\, -i\sqrt{1-h_0(\varrho)}\, , \ 
\varrho \in U_{1}\, .
\] 
We consider the following equation with respect to   $Q_{1}$ 
\begin{equation}
 \psi_1(\lambda)\left[\lambda R_{1}^-(\lambda) - K + (\lambda R_{0}^+(\lambda) - K ) Q_{1}(\lambda)\right]=
 O_{\mathcal B}(|\lambda|^{-M}) \ ,
\label{pdo-equation1}
\end{equation}
which we solve using
 the classes ${\rm PDO}_{\ell,2,M}(\Gamma;{\mathcal B};\lambda)$ introduced  in 
 Appendix A.2 (see Definition \ref{Def:finite-smoothness}). Hereafter, 
 \[
 O_{\mathcal B}(|\lambda|^{-M})\, :\ L^2(\Gamma) \to  L^2(\Gamma)
\]
 denotes 
any family of continuous  operators depending
on $K\in {\cal B}$ and on $\lambda\in {\cal D}$ with norms uniformly bounded by 
$ C_{\cal B}(1 + |\lambda|)^{-M}$,  where $C_{\cal B}>0$ is a constant independent
of $K\in {\cal B}$ and $\lambda$. 

Let $\left(R_0^+\right)^{-1}(\lambda)$ be a classical $\lambda$-PDO such that 
$\mbox{WF}'\left(\left(R_0^+\right)^{-1} R_0^+ -\mbox{Id}\right)\cap \mbox{WF}'(\psi_1) = \emptyset $. Then a solution of 
(\ref{pdo-equation1}) is given by
$$
Q_{1}= Q_{1}^0 + \lambda^{-1} Q_{1}^1 
$$
where
\begin{equation}
 Q_{1}^0 =-\left(R_0^+\right)^{-1} R_1^-\quad \mbox{and}\quad  Q_{1}^1 = \sum_{j=0}^{M-2} 
\lambda^{-j} \left(\left(R_0^+\right)^{-1}K\right)^{j+1} \left[\mbox{Id} -\left(R_0^+\right)^{-1}R_1^- \right]\, .
                                \label{pde-solution}
\end{equation}
The  operator $Q_{1}^1(\lambda)$ belongs to ${\rm PDO}_{\ell,2,M}(\Gamma;{\mathcal B};\lambda)$ in view of Remark \ref{Rem:commutator}, and it 
is well defined modulo $O_{\mathcal B}(|\lambda|^{-M})$. 
 The corresponding  
  principal symbols are
\begin{equation}
\sigma_0\left(Q_{1}^0\right)(x,\xi) = 1\ \quad \mbox{and}\quad \sigma_0\left(Q_{1}^1\right)(x,\xi) =
- \frac{2iK(x)}{\sqrt{1-h_0(x,\xi)}} 
= - \frac{2iK(x)}{\sin \theta(x,\xi)}
 \label{principal-symbols}
\end{equation}
in a neighborhood of $\mbox{WF}'(\psi_{1})$ in $U_{1}$. In this way the equation
\[
\imath_\Gamma^\ast\, \,   {\mathcal N}\,
G(\lambda)v = O_{M}(|\lambda|^{-M})v
 \]
 reduces to 
$(  W(\lambda) -  \mbox{Id}\, )\psi_0(\lambda)v \ = \ 
O_{\mathcal B}(|\lambda| ^{-M-1})v$, 
where $W(\lambda):= Q_1(\lambda)G_0(\lambda)$ represents  the microlocal monodromy operator.

Suppose now that $m\ge 2$. In order to satisfy the boundary conditions  at $U_{j+1}$, 
$0\le j\le m-2$, we  
are looking for  a $\lambda$-PDO $Q_{j+1}(\lambda)$  such that 
\begin{equation}
\psi_{j+1}(\lambda) \imath_\Gamma^\ast\,  {\mathcal N}\, 
H_{j+1}(\lambda)Q_{j+1}(\lambda) G_{j}(\lambda)\,
+\,  \psi_{j+1}(\lambda) \imath_\Gamma^\ast\,   {\mathcal N}\, 
H_j(\lambda)  = O_{\mathcal B}(|\lambda|^{-M}) \, .
                      \label{theboundary-identity}
\end{equation}
 Using the symbolic calculus  we write
$$
 \psi_{j+1}(\lambda)\imath_\Gamma^\ast\, \,   {\mathcal N}\,
H_{j+1}(\lambda) Q_{j+1}(\lambda)G_j(\lambda) 
 =\psi_{j+1}(\lambda) (\lambda R_{j+1}^+(\lambda) - K) Q_{j+1}(\lambda)G_j(\lambda) + O_{M}(|\lambda|^{-M})
$$
where $ R_{j+1}^+ (\lambda)$ is a classical $\lambda$-PDO  of order $0$  on $\Gamma$ 
independent of $K$, with a $C_0^\infty$-symbol in any local coordinates, and with 
principal symbol 
\[
\sigma(R_{j+1}^+)(\varrho)\, =\, i\sqrt{1-h_0(\varrho)}\, , \ 
\varrho \in U_{j+1}\, .
\] 
 In the same way  we obtain
$$
\psi_{j+1}(\lambda)\imath_\Gamma^\ast\,    {\mathcal N}\, 
 H_{j}(\lambda)  = \psi_{j+1}(\lambda)(\lambda R_{j+1}^- - K)\, 
G_{j}(\lambda)  +O_M(|\lambda|^{-M}) ,
$$
where $ R_{j+1}^- $ is a classical $\lambda$-PDO of order $0$  on $\Gamma$ 
independent of $K$  with 
principal symbol 
\[
\sigma(R_{j+1}^-)(\varrho)\, =\, -i\sqrt{1-h_0(\varrho)}\, , \ 
\varrho \in U_{j+1}\, .
\] 
Then (\ref{theboundary-identity}) reduces into the equation
\begin{equation}
\psi_{j+1}(\lambda)\left[(\lambda R_{j+1}^+ - K) Q_{j+1} + \lambda R_{j+1}^- - K \right]=  O_{\mathcal B}(|\lambda|^{-M})
                        \label{pdo-equation}
\end{equation}
on $U_{j+1}$, which we solve as above in 
 the classes ${\rm PDO}_{\ell,2,M}(\Gamma;{\mathcal B};\lambda)$.  More precisely, we obtain an operator 
$$
Q_{j+1}= Q_{j+1}^0 + \lambda^{-1} Q_{j+1}^1 
$$
 which
is well defined modulo $O_{\mathcal B}(|\lambda|^{-M-1})$,  
where 
$ Q_{j+1}^0$ is a  classical 
$\lambda$-PDOs of order $0$ independent of $K$ and  with a $C^\infty$ symbol, and 
$ Q_{j+1}^1\in {\rm PDO}_{\ell,2,M}(\Gamma;{\mathcal B};\lambda)$. The corresponding  
  principal symbols are
\[
\sigma_0(Q_{j+1}^0)(x,\xi) = 1\ ,\quad \sigma_0(Q_{j+1}^1)(x,\xi) =
-\frac{2iK(x)}{\sqrt{1-h_0(x,\xi)}} 
= -\frac{2iK(x)}{\sin \theta(x,\xi)}
\] 
in a neighborhood of $\mbox{WF}'(\psi_{j+1})$ in $U_{j+1}$.

Consider the family of operators $G(\lambda) : L^2 (\Gamma) \rightarrow 
C^\infty (\widetilde X)$, $\lambda\in {\cal D}$,  defined by 
\begin{equation}\label{operator-G}
G(\lambda) = H_0(\lambda)\psi_0(\lambda) + 
\sum_{k=2}^{m}
H_{k-1}(\lambda)\left( \Pi_{j=0}^{k-2}
Q_{j+1}(\lambda) G_{j}(\lambda)\right)\psi_0(\lambda) \,  .
\end{equation}
Using (\ref{psi-pdos}) -  
(\ref{boundary-trace}) and  (\ref{theboundary-identity}) 
we obtain
\[
\left\{
\begin{array}{rcll}
(\Delta-\lambda^2)G(\lambda)\ &=&\ 
O_{\mathcal B}(|\lambda|^{-M})\, , \\
\imath^\ast_\Gamma \,{\mathcal N}\,  G(\lambda)  \ &=&\  
\psi_{m}(\lambda)(\lambda R_{0}^+ - K)\, \psi_0(\lambda)  + 
\psi_{m}(\lambda)(\lambda R_{m}^- - K)\, 
\widetilde{W}(\lambda)\psi_0(\lambda)
 +O_{\mathcal B}(|\lambda|^{-M})  \, ,
\end{array}
\right.
\]
where 
\[
\widetilde{W}(\lambda)=\imath^\ast_\Gamma 
H_{m-1}(\lambda)\Pi_{j=0}^{m-2} \left(\psi_{j+1}(\lambda)Q_{j+1}(\lambda)
G_{j}(\lambda)\right)\, , 
\]
and $R_{0}^+$ and $R_{m}^-$ are defined as above. As in 
(\ref{pdo-equation1}) we find $Q_{m} = Q_{m}^0 + \lambda^{-1}  Q_{m}^1$ such that  
$$
 \psi_m(\lambda)\left[\lambda R_{m}^- - K + (\lambda R_{0}^+ - K ) Q_{m}(\lambda)\right]=
 O_{\mathcal B}(|\lambda|^{-M}) \ ,
$$
where $Q_{m}^k$, $k=0,1$,  are as above. Denote
\begin{equation}
W(\lambda):= Q_m(\lambda)\widetilde W(\lambda) =
\Pi_{j=0}^{m-1}  \left(\psi_{j+1}(\lambda)Q_{j+1}(\lambda) G_j(\lambda)\right)\, .  
\label{monodromy}
\end{equation}
Then the boundary problem above becomes
\[
\left\{
\begin{array}{rcll}
(\Delta-\lambda^2)G(\lambda)\ &=&\ 
O_{\mathcal B}(|\lambda|^{-M})\, , \\
\imath^\ast_\Gamma \,{\mathcal N}\,  G(\lambda)  \ &=&\  
\psi_{m}(\lambda)(\lambda R_{0}^+ - K)(  \mbox{Id}\, - \, W(\lambda))\, \psi_0(\lambda)
 +O_{\mathcal B}(|\lambda|^{-M})  \, .
 \end{array}
\right.
\]
In this way we reduce the equation 
$\imath^\ast_\Gamma \,{\mathcal N}\,  G(\lambda)v
=O_{\mathcal B}(|\lambda|^{-M})v$ to the following one 
\begin{equation}
(  W(\lambda) -  \mbox{Id}\, )\psi_0(\lambda)v \ = \ 
O_{\mathcal B}(|\lambda| ^{-M-1})v\, .
                               \label{equation}
\end{equation}
The operator $W(\lambda)$ defined by \eqref{monodromy} will be  called a microlocal monodromy operator. 
We summarize the above construction  by the following 
\begin{Prop} \label{prop:monodromy1}
Suppose that there is a sequence $(\lambda_q,v_q)\in \R_+ \times L^2(\Gamma)$, $q\in {\cal M}$, such that ${\cal M}\subset \Z^n$ is an infinite index set, $\lim_{|q|\to \infty}\lambda_q = +\infty$, $\|v_q\|_{L^2(\Gamma)}=1$, and 
\[
(  W(\lambda_q) -  \mbox{Id}\, )\psi_0(\lambda_q)v_q \ = \ 
O_{\mathcal B}(\lambda_q ^{-M-1})\, .
\]
Set $u_q= G(\lambda_q)v_q$. Then 
\[
\left\{
\begin{array}{rcll}
(\Delta-\lambda_q^2)G(\lambda_q)\ &=&\ 
O_{\mathcal B}(\lambda_q^{-M})\, , \\
\imath^\ast_\Gamma \,{\mathcal N}\,  u_q  \ &=&\  
 O_{\mathcal B}(\lambda_q^{-M})  \, .
 \end{array}
\right.
\]
\end{Prop}
We are interested in the structure of the monodromy operator. Recall from \eqref{boundary-trace}, \eqref{boundary-trace3} and \eqref{boundary-trace4} that  $G_j(\lambda)= E(\lambda)^{-1}G_j^0(\lambda)E(\lambda)$, where $E(\lambda)$ is a $\lambda$-PDO of order $0$ which is elliptic on the union of $\overline U_j$, $1\le j\le m$,  and that the principal symbol of $G_j^0(\lambda)$ is given by \eqref{principal-symbol-G}. 
Set 
$$
S(\lambda):=\Pi_{j=0}^{m-1} G_j^0(\lambda).
$$
By construction $G_j(\lambda)$ is elliptic on $\mbox{WF}'(\psi_{j}Q_{j})$, and 
using Lemma \ref{Lemma:commutator} we can commute $G_j(\lambda)$ with $\psi_{j}Q_{j}$. 
Since $ {\rm PDO}_{\ell,2, M}(\Gamma;{\mathcal B};\lambda)$ 
is closed under multiplication (see Remark \ref{Rem:commutator}), we get 
another $\lambda$-PDO of the same class which we commute with $G_{j+1}(\lambda)$ and so on. Then  setting 
$\widetilde \psi_0(\lambda):=E(\lambda)^{-1}\psi_0(\lambda)E(\lambda)$ and $\widetilde \psi_m(\lambda):=E(\lambda)^{-1}\psi_m(\lambda)E(\lambda)$
 we obtain  the following 
\begin{Prop} \label{prop:monodromy2}
The microlocal monodromy operator can be written as follows
\[
W_0(\lambda):=E(\lambda) W(\lambda)E(\lambda)^{-1}=
\widetilde \psi_m(\lambda) \left(Q^0(\lambda) + \lambda^{-1}Q^1(\lambda)\right)
S(\lambda)\widetilde \psi_0(\lambda)   + O_{\mathcal B}(\lambda^{-M-1}) \, ,
\]
where $Q^0(\lambda)$ is a classical  $\lambda$-PDOs on $\Gamma$ with a $C^\infty$ symbol
independent of $K$ and  with a  principal symbol equal to 
 one in a neighborhood of $\Lambda$, and 
 $ Q^1\in {\rm PDO}_{\ell,2, M}(\Gamma;{\mathcal B};\lambda)$. Moreover,  
 the principal
symbol of $Q^1(\lambda)$ is 
$$
\sigma_0(Q^1)(x,\xi) = - 2i \sum_{j=0}^{m-1} \frac{K(\pi_\Gamma(x^j,\xi^j))}{\sin
\theta (x^j,\xi^j)}\, ,\quad (x^j,\xi^j) = B^{-j}(x,\xi)\, ,
$$
in  $P(U_0)$. 
The operator  $ S(\lambda)$ does not depend on $K$, 
and it is a classical $\lambda$-FIO 
of order 0 with a
large parameter $\lambda\in {\cal D}$. The canonical relation of 
$  S(\lambda)$ is given by
the graph of the symplectic map 
$P= B^m:U_0 \to U_m$. Moreover,  parameterizing ${\rm graph}\, P$ by its projection on $U_0$,  the principal symbol of $S(\lambda)$ becomes 
$\exp\left(i \frac{1}{2}\vartheta\right)\exp (i\lambda A(x,\xi))$, $(x,\xi)\in U_0$, where $\vartheta\in\Z$ and   $A(x,\xi) =
\sum_{j=0}^{m-1}A_j(x^j,\xi^j)$ is the action along the corresponding broken geodesic.
\end{Prop}
The principal symbol $\sigma_0(Q^1)$ of the operator $Q^1$ is obtained by means of the  Egorov's theorem (see Lemma \ref{Lemma:commutator}). 
In what follows we will find a simple microlocal normal form of $W_0(\lambda)$. To this end we will use a suitable symplectic normal form of $P$. 

\subsection{ Birkhoff normal  form of $P$. }\label{subsec:BNF}
The exact symplectic map $P=B^m : U_0\to U_m$ can be put in  a symplectic Birkhoff normal form in a neighborhood of  $\Lambda$ as follows. 
Denote by $f: \T^{n-1}\to B ^\ast \Gamma$ a smooth  embedding given by (H$_3$) such that  $\Lambda= f(\T^{n-1})$ 
and   the    diagram
\begin{equation}
\label{diagram} 
\displaystyle{\begin{array}{cccl} 
\T^{n-1}&\stackrel{R_{2\pi \omega}}{\longrightarrow}&\T^{n-1}\cr
\downarrow\lefteqn{f}& &\downarrow\lefteqn{f} \cr
\Lambda&\stackrel{P}{\longrightarrow}&\Lambda&  
\end{array} }
\end{equation}
is commutative,  
 where  
$\omega$ satisfies the Diophantine condition  (\ref{sdc}). 
Denote by $\gamma_j^0$, $j=1,\ldots, n-1$, the cycles 
$\gamma_j^0:=\{(0,\ldots,0,\varphi_j,0,\ldots, 0):\, \varphi_j\in\T\}$ of $\T^{n-1}$
 and set $\gamma_j = f\circ \gamma_j^0$.   
Then  
$\gamma_j$, $j=1,\ldots, n-1$, is  a basis of cycles of the first homology group
$H_1(\Lambda,\Z)$, and we set $I^0 = (I^0_1,\ldots,I^0_{n-1})$, 
 where 
 $$
 I^0_j = (2\pi )^{-1}\int _{\gamma_j} \xi d x\, .
 $$ 
Denote by $\imath : \T^{n-1}\to T^\ast \T^{n-1}$ the embedding $\imath(\varphi) = (\varphi, I^0)$. 
 Then Proposition 9.13, \cite{La}, implies
\begin{Prop}\label{BirkhoffNormalForm}   
There is a neighborhood $\A=\T^{n-1}\times D'$ of 
$\T^{n-1}\times \{I^0\}$ in $T^\ast \T^{n-1}$ and  an exact symplectic map $\chi : \A \to U_m\subset B ^\ast  \Gamma$ such that 
\begin{enumerate}
\item[(i)]   $\chi\circ \imath= f$,  
\item[(ii)]
the symplectic map $P^0 := \chi^{-1}\circ P \circ
\chi$ is defined by a generating function 
\begin{equation}
\label{generating-function}
\Phi^1(x,I)=\langle x,I \rangle +\Phi(x,I) \, ,\ (x,I)\in \R^{n-1}\times D\, ,  
\end{equation}
where $ \Phi\in C^\infty(\R^{n-1}\times D)$ is 
$2\pi$-periodic in $x$, $D$ is a neighborhood of $I^0$,  and 
\item[(iii)]
$\Phi(x,I)=L(I) + \Phi^0(x,I)$, where 
$\nabla L (I^0) = 2\pi  \omega$ and $\partial_I^\alpha  \Phi^0(x,I^0)=0$, for any $x\in
\R^{n-1}$, and $\alpha \in \N^{n-1}$. 
\end{enumerate}
\end{Prop}
To simplify the notations  we denote the class of $x\in \R^{n-1}$ in $\T^{n-1}$ by the same letter $x$.
Recall that  the function $\Phi^1(x,I) $ in \eqref{generating-function}  is a generating function of the symplectic map $P^0$ in  $\R^{n-1}\times D$ if  
$$P^0\left(\nabla_I\Phi^1(x,I),I\right) =  \left(x, \nabla_x \Phi^1(x,I)\right)\, ,\ (x,I)\in \R^{n-1}\times D\, .$$ 
In particular, we have     $\chi (\T^{n-1}\times \{I^0\})= \Lambda$ and          
\begin{equation}
\forall\,   N\in
\N\, , \quad P^0(\varphi, I) = (\varphi -\nabla L(I), I) + O_N(|I-I_0|^N)\, .
                           \label{Birkhoff}
\end{equation}
Proposition 9.13, \cite{La}, provides a symplectic mapping  $\chi$  with the desired properties. Hereafter we take 
 $D'$ to be a small ball in $\R^{n-1}$ centered at $I^0$. Then  $\chi$ is exact symplectic in $\A$  due to the choice of $I^0$. Indeed, the integrals of the closed one-form $\alpha:=\chi^{\ast} (\xi dx) - I d\varphi$ on the cycles $\gamma_j^0\times \{I^0\}$ are all zeros, hence, the class of $\alpha$ in the first cohomology group $H^1(\A,\R)$ is zero and $\alpha$ is exact on $\A$. 

The function $L$ is defined modulo a constant and in 
what follows we shall give a natural choice of  $L(I^0)$ which is related to the choice of the phase function  in \eqref{operator}. 
Moreover, we shall give a geometric interpretation of $L$.

Consider the ``flow-out'' ${\cal T}$ of $\Lambda$ in $T^\ast X$ 
by means of the broken bicharacteristics of $h$. Actually, the propagation along the broken bicharacteristics of $h$ is not even continuous because of the reflections at the boundary. To work with  smooth objects one can  use the compressed cotangent bundle or equivalently the symplectic gluing technique which we explain below.

Recall that the broken bicharacteristic flow of $h$ can be extended continuously from the interior of  $T^\ast X$ to the compressed cotangent bundle by gluing together $\pi_\Sigma^-(\varrho)$ and $\pi_\Sigma^+(\varrho)$ for $\varrho\in B^\ast\Gamma$. Moreover,  
one can use the method of the symplectic gluing (\cite{La}, Chapter 1.4, see also \cite{PT2}, Appendix) to extend it smoothly across the boundary near $\Lambda_j$, $j=1,\ldots,m-1$. More precisely,  there exists a unique smooth symplectic gluing of the symplectic manifold $T^\ast X$ across  a neighborhood  $U\subset T^\ast X|_\Gamma$  of $\cup_{j=1}^{m-1}\pi_\Sigma^{\pm}\Lambda_j$ into a symplectic manifold $\widetilde{T^\ast X}$ and a unique extension of $h$ to a smooth
Hamiltonian $\tilde h$ on $\widetilde{T^\ast X}$ such that the canonical projection 
$p: T^\ast(X\setminus \Gamma)\cup U \to \widetilde{T^\ast X}$ is a smooth symplectic map on   $T^\ast (X\setminus \Gamma)$, and 
$\tilde h\circ p= h$ on $T^\ast(X\setminus \Gamma)\cup U$. In this way, the broken bicharacteristic flow is represented locally by the smooth Hamiltonian flow  $(t,\varrho)\to \exp(tX_{\tilde h})(\varrho)$, $\varrho\in \widetilde{T^\ast X}$, and the flow-out of $\Lambda$ is a smooth Lagrangian torus ${\mathcal T}$ in $\widetilde{T^\ast X}$, invariant with respect to the flow of $\tilde h$. Now $\Lambda$ can be considered as  a smooth section of ${\mathcal T}$, 
$P$ can be identified with the corresponding Poincaré map,  and $2A(\varrho)$, $\varrho\in\Lambda$, is the first return time. In particular, there is a diffeomorphism $\phi:{\mathcal T}\to\T^n$  conjugating the restriction of the Hamiltonian flow of $\tilde h$ on ${\mathcal T}$ to the linear flow $(\varphi,\varphi_n)\to (\varphi,\varphi_n) -t(\omega,1)$, 
$(\varphi,\varphi_n)\in\T^n$. 
%Recall that  $\omega$ is Diophantine. Hence, the vector $(\omega,1)\in \R^n$ is non-resonant and the restriction of the flow of $\tilde h$ to  ${\mathcal T}$ is uniquely ergodic, i.e.  there exists a unique probability measure $d\tilde\mu$ on ${\mathcal T}$ (associated with the Borelian $\sigma$-algebra) which is invariant with respect to the flow oh $\tilde h$. It is given by $(2\pi)^{-n}\phi^\ast(d\varphi d\varphi_n)$. 

We  choose the constant  $L(I^0)$ as follows. 
Let $\rho^0=\chi(\varphi^0,I^0)\in \Lambda$. 
Denote by $\gamma_{n1}(\rho^0)$ the broken bicharacteristic arc in ${\cal
T}$ issuing from $\rho^0$ and having  endpoint at $P(\rho^0)$, and by 
$$\gamma_{n2}(\rho^0):= \chi(\varphi^0 + (s-1)2\pi \omega,I^0)\, ,\ s\in
[0,1],$$ an  arc 
 connecting $P(\rho^0)$ and
$\rho^0$ in $\Lambda$. 
Let $\gamma_n$ be the union of the two arcs. We denote by $L(I^0)$ the
action along $\gamma_n$, i.e. 
\begin{equation}
L(I^0)\ =\ \, \int_{\gamma_n}\, \xi d x\, . 
                                  \label{action}
\end{equation}
Note that the integral   above   depends only on the
homology class of the loop  $\gamma_n$ in the Lagrangian torus 
${\cal T}$. We can give now a geometric interpretation of  $L$ which
will be needed later.  
The Poincaré identity gives 
\begin{equation}
P^\ast  (\xi dx) = \xi dx + dA,
\label{poincare}
\end{equation}
where $\xi dx$ is the fundamental one-form on $T^\ast\Gamma$ and
$A(\rho)$, $\rho =\chi(\varphi, I)$, $|I-I^0| \ll 1$,  
stands for the action along the broken bicharacteristic
$\gamma_{n1}(\rho)$.  Since $\chi$ is exact symplectic  in $\A\subset T^\ast \T^{n-1}$ 
we have 
\begin{equation}
\chi^\ast(\xi dx) = I d\varphi +
d\Psi
\label{exact-symplectic}
\end{equation}
with a suitable smooth function $\Psi\in C^\infty(\A)$. 
\begin{Lemma}
\label{length}
We have 
\[
 L(I) - \langle I, \nabla L(I)\rangle  = A(\chi(\varphi,I)) + \Psi(\varphi,I) -
\Psi(P^0(\varphi,I))  + O_p(|I-I^0|^p)
\]
for any $p\in \N$..
\end{Lemma}
{\em Proof}. 
Combining \eqref{poincare} and \eqref{exact-symplectic} 
we obtain 
\[
(P^0)^\ast  (I d\varphi) - I d\varphi = d((A\circ \chi) +  \Psi -
\Psi\circ P^0). 
\]   
In view of   (\ref{Birkhoff}) this implies  
\[
 L(I) - \langle I, \nabla L(I)\rangle  = A(\chi(\varphi,I)) + \Psi(\varphi,I) -
\Psi(P^0(\varphi,I)) +C+  O_p(|I-I^0|^p)
\]
for any $p\in \N$, where $C\in \R$. Notice that $C$ should be  zero since    
 for 
 $I= I^0$ and $\omega =\nabla L(I^0)/2\pi $ it follows by  (\ref{action}) and \eqref{exact-symplectic} that
\[
\begin{array}{lcrr}
\displaystyle L(I^0) - \langle I^0, \nabla L(I^0)\rangle =
 L(I^0)  - 2\pi \langle I^0, \omega\rangle = 
\int_{\gamma_{n}} 
\xi d x - \int_{\gamma_{n2}^0} 
I^0 d \varphi \\ [0.3cm]
\displaystyle =
\int_{\gamma_{n1}(\rho^0)} \xi d x
+  \Psi(\varphi^0,I^0) - \Psi(\varphi^0-2\pi \omega,I^0)
= A(\chi(\varphi^0,I^0)) + \Psi(\varphi^0,I^0) -
\Psi(P^0(\varphi^0,I^0)) \, , 
\end{array}
\]
where $\gamma_{n2}^0 := (\varphi^0 + (s-1)2\pi \omega,I^0)=\chi^{-1}(\gamma_{n2}(\rho^0))$.  This proves the Lemma. 
\finishproof
\begin{Lemma}
\label{positiv-action}
We have 
\[
L(I^0) -  2\pi \langle I^0, \omega\rangle =    
\int_{\Lambda} A(\varrho)\,  d\mu \,  > \ 0 .                
\]
\end{Lemma}
{\em Proof}. 
 Set $\varrho^j = P^j(\varrho^0)
= \chi(\varphi^0 -2\pi j\omega,I^0)$. 
The measure  
$d\mu = (2\pi)^{-n+1}\chi_\ast (d\varphi)$   on $\Lambda$ is invariant with respect
to the map $P:\Lambda \to \Lambda$ which is ergodic since $\omega$ satisfies (\ref{sdc}). The Birkhoff ergodicity theorem implies 
\[
L(I^0) -  2\pi \langle I^0, \omega\rangle = 
\lim_{j\to \infty}\frac{1}{j} \sum_{k=0}^{j-1}
A(\varrho^k) =    
\int_{\Lambda} A(\varrho)\,  d\mu \,  .                
\]
\finishproof
\subsection{ Quantization of $\chi$. }\label{subsec:quantization}
Using the embedding 
 $f$ in \eqref{diagram}  
 we identify the first cohomology groups 
$H^1(\Lambda,\Z) =  H^1(\T^{n-1},\Z)=\Z^{n-1}$,
and we denote by $\vartheta_0\in \Z^{n-1}$ 
the Maslov class of the invariant
torus $\Lambda$. 
As in \cite{CV} we consider 
the flat Hermitian line bundle $\LL$ over
$\T^{n-1}$ which is associated to the class $\vartheta_0$. 
The sections $s$
in ${\LL}$ can be identified canonically with functions
$\widetilde{s}:\R^{n-1}
\rightarrow \C$ so that
\begin{equation}
\widetilde{s}(x +2\pi p)\ =\ e^{i\frac{\pi}{2}\langle 
\vartheta_0,p\rangle}
\widetilde{s}(x)
                                              \label{sections}
\end{equation}
for each $x\in \R^{n-1}$ and $p\in \Z^{n-1}$. 
An orthonormal basis of $L^2(\T^{n-1},\LL)$ 
is given by
$e_k,\ k\in \Z^{n-1}$,  where
$$
\widetilde e_k (x)\ =\ 
\exp\left( i \langle k + \vartheta_0/4,x\rangle
\right).
$$ 
We quantize the exact symplectic  transformation $\chi: \A=\T^{n-1}\times D' \to T^\ast \Gamma$ as in \cite{CV}. 
More
precisely,  we find  a classical $\lambda$-FIO 
$$ T(\lambda): C^\infty(\T^{n-1}, \LL) \rightarrow 
C^\infty(\Gamma)$$ of order $0$ the canonical relation of which is just the graph
of  $\chi$ and such
that 
$$\mbox{WF}'(T(\lambda) T(\lambda)^\ast - \mbox{Id}_\Gamma)\cap \A^0 =\emptyset. $$
The last formula means that $T(\lambda)$ is a microlocally unitary operator over $\A^0:=\T^{n-1}\times D^0$, where $D^0$ is a small ball centered at $I^0$ and contained in  $D'$. 
One  can take the principal symbol of $T(\lambda)$ to be  equal to one 
in $\T^{n-1}\times D^0$ modulo a Liouville
factor which is given by $\exp(i\lambda \Psi(\varphi,I))$, 
where  the function $\Psi$ satisfies \eqref{exact-symplectic}  (see \eqref{liouville-factor}). Consider the operator
\[
T(\lambda)^{\ast}  W_0(\lambda) T(\lambda) = \left[T(\lambda)^{\ast}  \left(Q^0(\lambda) + \lambda^{-1}Q^1(\lambda)\right)T(\lambda)\right]\left[
T(\lambda)^{\ast}  S(\lambda) T(\lambda)\right].
\]
Using Lemma \ref{Lemma:commutator} and Remark \ref{Rem:commutator} we write the first factor as a sum $P^0(\lambda) + \lambda^{-1}P^1(\lambda)$, where 
$P^0(\lambda)$ is a classical $\lambda$-PDO of order $0$ acting  on $C^\infty(\T^{n-1}, \LL)$ with a $C^\infty$ symbol
independent of $K$ and  with a  principal symbol equal to 
 one in a neighborhood of $\A^0$, while 
 $ P^1$ is in ${\rm PDO}_{\ell,2, M}(\T^{n-1}, \LL;{\mathcal B};\lambda)$ and its principal symbol is  $\sigma_0(P^1)(x,I) = \sigma_0(Q^1)(\chi(x,I))$. 
 
The second factor is a  composition of three $\lambda$-FIOs,  whose canonical relations are graphs of exact symplectic transformations. Then 
 $S^0(\lambda)$ is a $\lambda$-FIO   of order $0$ and its  canonical relation $\mathcal C$  is  the graph of $P^0= \chi^{-1}\circ P \circ \chi$, i.e. 
 $\mathcal C:= \{(P^0(\varrho),\varrho): \varrho\in\A \}$. Denote by ${\mathcal C}'$ the corresponding Lagrangian submanifold of $T^\ast(\T^{n-1}\times \T^{n-1})$ via the relation \eqref{lagrangian-manifold}. It follows from the    $\lambda$-FIO calculus and Proposition \ref{prop:monodromy2} that parameterizing 
$\mathcal C$ by the variables $\varrho=(\varphi,I)\in \A$  the principal symbol of $S^0(\lambda)$ becomes
\[
\sigma( S^0(\lambda))(\varphi,I)= \exp\left(i \frac{1}{2}\vartheta\right)\exp (i\lambda f(\varphi,I)),
\] 
where the exponent of the Liouville factor  is
\begin{equation}
\label{liouville-exponent}
f(\varphi,I)= A(\chi(\varphi,I)) + \Psi(\varphi,I) -
\Psi(P^0(\varphi,I)).
\end{equation} 
On the other hand, $P^0$ is defined by the generating function $\Phi^1(x,I)=\langle x,I\rangle + \Phi(x,I)$, $(x,I)\in \T^{n-1}\times D'$,  in Proposition \ref{BirkhoffNormalForm}, hence, 
the Lagrangian manifold ${\mathcal C}'$ is defined globally by any phase function $\widetilde \Phi(x,y,\eta):= \Phi^1(x,I)- \langle y,I\rangle + C$, where $C\in \R$ is a constant. 
By the   $\lambda$-FIO calculus  there is $C\in \R$ and a smooth amplitude  $s(x,I,\lambda)= s_0(x,I) + \lambda^{-1}s_1(x,I) + \cdots $, $(x,I)\in \R^{n-1} \times D'$,  of order $0$, which is $2\pi$-periodic with respect to $x$ and 
 uniformly compactly supported in  $I\in D'$, such  that 
\[
\widetilde{S^0(\lambda)u}(x) = 
\left(\frac{\lambda}{2\pi}\right)^{n-1}\, \int_{\R^{2n-2}} 
\, e^{i\lambda (\langle x-y,I\rangle + \Phi(x,I)+C)} \, s(x,I,\lambda)\, J(x,I)\, 
\widetilde u (y) \, dI dy\,  , 
\]
for any $u\in C^\infty(\T^{n-1},\LL)$, where  
\[
J(x,I) \, =\,  
\left|{\rm det}\ \left({\rm Id}_{\R^{2n-2}} + 
\frac{\partial^2 \Phi}{\partial x\partial I}(x,I)\right)\right|^{-1/2} \, .
\]
In particular, $s_0(x,I)= \exp\left(i \frac{1}{2}\vartheta\right)$ in $\R^{n-1} \times D^0$. Calculating the Liouville factor in terms of the phase function 
$\langle x-y,I\rangle + \Phi(x,I)+C$ (see \eqref{liouville-factor})
we shall show that $C=0$. Indeed, 
\[
{\mathcal C}'={\mathcal C}'_{\widetilde \Phi}= \{(x,x+\Phi'_I(x,I),I +\Phi'_x(x,I),-I):\, (x,I)\in \T^{n-1} \times D \},
\]
 and we get $f(x,I) = \Phi(x,I) - \langle I, \Phi'_I(x,I)\rangle + C$. Hereafter, to simplify the notations we denote the class of $x\in\R^{n-1}$ in  $\T^{n-1}$ by $x$ as well. Comparing it with \eqref{liouville-exponent} at $(x,I^0)$ and using Lemma \ref{length} we get $C=0$.

Using  Remark \ref{Rem:commutator} we obtain 
\begin{equation}
\label{conjugation}
T(\lambda)^{\ast}  W_0(\lambda) T(\lambda) = \left( P^0(\lambda) + \lambda^{-1}P^1(\lambda) \right)
S^0(\lambda)= e^{i\pi\vartheta/2}
W_1(\lambda) + O_{\mathcal B}(|\lambda|^{-M}), 
\end{equation}
where 
\[
\widetilde{W_1(\lambda)u}(x) = 
\left(\frac{\lambda}{2\pi}\right)^{n-1}\, \int_{\R^{2n-2}} 
\, e^{i\lambda (\langle x-y,I\rangle + \Phi(x,I))} \, w(x,I,\lambda)\, J(x,I)\, 
\widetilde u (y) \, dI dy\,  , 
\]
for any $u\in C^\infty(\T^{n-1},\LL)$. 
The symbol $w(x,I,\lambda)$, $(x,I)\in \R^{n-1} \times D$,  is $2\pi$-periodic with respect to $x$
and uniformly compactly supported in  $I\in D$,  
and 
$w= w_0 + \lambda^{-1}w^0$, where $w_0\in C^\infty(\R^{n-1}\times D)$ is independent of $K$, 
$w_0(x,I)= 1$ for  any $x\in \R^{n-1}$ and $I\in D^0$, 
and 
\begin{equation}
\label{w-symbol}
w^0=\sum^{M-1}_{j=0}w^0_j(x,I)\lambda^{-j}\in S_{\ell,2,M}( \T^{n-1}\times D; {\mathcal B}; \lambda)\, .
\end{equation} 
Moreover, 
\begin{equation}
\label{first-invariant1}
w^0_0(x,I) = iw_0'(x,I) - 2i\sum_{j=0}^{m-1} \left(\frac{K\circ\pi_\Gamma}{\sin
\theta}\right)\left(B^{-j}\chi(x,I)\right)\, ,
\end{equation}
where $w_0'$ is a $C^\infty$ real-valued function independent of $K$ . 
 Moreover, 
\begin{equation}\label{jacobian}
\forall\,  N\in \N\, ,\quad J(x,I) =   1 + O_N(|I-I^0|^N)\, ,
\end{equation}
since $\Phi^0$ is flat at $I=I^0$. We summarize the above construction by the following
\begin{Prop}\label{operator-T}
There is a classical  $\lambda$-FIO  $T(\lambda): C^\infty(\T^{n-1}, \LL) \rightarrow 
C^\infty(\Gamma)$ of order $0$ with a  canonical relation  given by the graph
of  $\chi: \A \to T^\ast \Gamma$ which is microlocally unitary over $\A^0= \T^{n-1}\times D^0$ and with a principal symbol equal to one in $\A^0$ modulo a Liouville factor,  and such that 
\[
T(\lambda)^{\ast}  W_0(\lambda) T(\lambda)
 = e^{i\pi\vartheta/2}
W_1(\lambda) + O_{\mathcal B}(|\lambda|^{-M}), 
\]
where $W_1(\lambda): C^\infty(\T^{n-1}, \LL)\to C^\infty(\T^{n-1}, \LL)$ is a $\lambda$-FIO of order zero with a canonical relation given by the graph of $P^0$ over $\A$. Moreover, the Schwartz kernel of $W_1(\lambda)$ is of the form 
\begin{equation}
\widetilde{W_1}(x,y,\lambda) = 
\left(\frac{\lambda}{2\pi}\right)^{n-1}\, \int_{\R^{n-1}} 
\, e^{i\lambda (\langle x-y,I\rangle + \Phi(x,I))} \, w(x,I,\lambda)\, J(x,I)\, 
 \, dI \,  , 
                                  \label{operator}
\end{equation}
where the phase function $\Phi$ is given by Proposition \ref{BirkhoffNormalForm}, 
$w= w_0 + \lambda^{-1}w^0$,  $w_0\in C^\infty(\R^{n-1}\times D')$ is independent of $K$ and compactly supported in $D'$, 
$w_0(x,I)= 1$ for  $(x,I)\in \R^{n-1}\times D^0$,  $w^0$ satisfies \eqref{w-symbol} and \eqref{first-invariant1}, and the Jacobian $J$ satisfies \eqref{jacobian}. 
\end{Prop}

\subsection{ Homological equation and Quantum Birkhoff Normal  Form. }\label{subsec:QBNF}

Our goal is to get rid of the angle variable $x$ in $w(x,I,\lambda)$. To explain the main ideas in the construction below consider the ``essential part'' of 
$\widetilde{W_1}(\lambda,x,y)$ which is given by
\[
\widetilde{W_2}(x,y,\lambda) = 
\left(\frac{\lambda}{2\pi}\right)^{n-1}\, \int_{\R^{n-1}} 
\, e^{i\lambda (\langle x-y,I\rangle + L(I))} \, (w_0(I)+ \lambda^{-1}w^0_0(x,I))\, 
 \, dI \,  , 
\]
and denote by $W_2(\lambda)$ the corresponding $\lambda$-FIO. We have neglected the terms $O_N(|I-I^0|^N)$ since  $I-I^0 \sim \lambda^{-1}$ on the  frequency support of the quasimode that we are going to construct.  
We are looking for a $\lambda$-FIO $A(\lambda)$ of order $0$ which is   elliptic at $\T^{n-1}\times\{I^0\}$ and  such that
\begin{equation}
W_2(\lambda)A(\lambda) = A(\lambda) W_2^0(\lambda) + O(|\lambda|^{-2}),
                                  \label{operator-equation}
\end{equation}
where 
\[
\widetilde{W_2^0}(x,y,\lambda) = 
\left(\frac{\lambda}{2\pi}\right)^{n-1}\, \int_{\R^{n-1}} 
\, e^{i\lambda (\langle x-y,I\rangle + L(I))} \, (w_0(I)+ \lambda^{-1}p^0(I)\, 
 \, dI \,  , 
\]
The symbol of $A(\lambda)$ has the form $a_0(x,I) + \lambda^{-1}a^0(x,I) + O(|\lambda|^{-2})$ and we can take  $a_0=1$  in a neighborhood of $\T^{n-1}\times\{I^0\}$. Then the zero order terms of the symbols of the operators in \eqref{operator-equation} equal one in a neighborhood of $\T^{n-1}\times\{I^0\}$. Comparing the corresponding symbols of order $-1$  we obtain  the homological equation 
\[
a^0(x-\nabla L(I), I) - a^0(x, I) = w^0(x,I) - p^0(I)\, ,\quad I \in D^0. 
\]  
We are looking for a pair of functions $a^0(x,I)$ and $p^0(I)$ satisfying the above equation modulo $O(|I-I^0|^N)$ for suitable $N\ge 1$. 
We show below that this equation can be solved when $I=I^0$ since $\nabla L(I^0) = 2\pi\omega$ and $\omega$ is Diophantine.
Then taking the Taylor polynomials of  $w^0(x,I)$ at $I=I^0$ we can find successively the Taylor coefficients of $a^0(x,I)$ and $p^0(I)$ at 
$I=I^0$ up certain  order   which suffices for our purposes.

Consider the  homological equation 
\begin{equation}
{\cal L}_\omega u(\varphi) := u(\varphi - 2\pi \omega) -u(\varphi)=f(\varphi),\ \varphi\in \T^{n-1}.
\label{homological}
\end{equation} 
The solution $u$, if it exists, is less regular than $f$. The lost of regularity depends on the exponent $\tau$ in the 
 the small denominator condition 
(\ref{sdc}).  In order to give sharp estimates on the regularity of $u$ it is convenient to use weighted Wiener spaces ${\mathcal A}^s(\T^{n-1})$, $s\ge 0$, which are defined as follows (see also \cite{M-P}). 
Given $u\in C(\T^{n-1})$, we denote by $u_k$, $k\in\Z^{n-1}$, the corresponding Fourier coefficients, 
and for any $s\in\R_+:=[0,+\infty)$ we define the corresponding weighted Wiener norm of $u$ by
\[
\|u\|_s\, :=\ \sum_{k\in\Z^{n-1}}\,  (1+|k|)^s|u_k|\, ,
\]
where  $|k|=|k_1|+\cdots+|k_{n-1}|$ for any $k\in\Z^{n-1}$. We denote by ${\mathcal A}^s(\T^{n-1})$ the Banach space of all $u\in C(\T^{n-1})$ such that $\|u\|_s<\infty$. We list below some useful properties of these spaces. 
If $u\in {\mathcal A}^s(\T^{n-1})$ and $\alpha\in \N^{n-1}$ with $|\alpha|\le s$, then $\partial^\alpha u \in {\mathcal A}^{s-|\alpha|}(\T^{n-1})$. The space ${\mathcal A}^s(\T^{n-1})$ is  a Banach algebra, if $u,v \in {\mathcal A}^s(\T^{n-1})$ then $\|uv\|_s\le \|u\|_s\|v\|_s$, and if $u\neq 0$, then $1/u \in {\mathcal A}^s(\T^{n-1})$ as well. 
Moreover, 
the  following relations between Wiener spaces and Hölder spaces on the torus $\T^p$, $p\ge 1$, hold
\begin{equation}
\label{bernstein}
C^q(\T^p) \hookrightarrow {\mathcal A}^s(\T^p)\hookrightarrow C^s(\T^p)\, ,
\end{equation}
for any $s\ge 0$ and $q>s+p/2$, and the corresponding inclusion maps are continuous. The first relation   (cf. \cite[Sect. 3.2]{kn:AAP}) is a special case of a theorem of Bernstein ($p=1$) and its generalizations for $p\ge 2$. 
One should be  careful when changing the variables on the torus, since a continuous map $\phi:\T^p\to\T^p$ preserves ${\mathcal A}^s(\T^p)$ if and only if $\phi$ is affine linear. 
Wiener spaces are perfectly adapted for solving \eqref{homological},  provided that  $\omega$ satisfies the Diophantine condition  (\ref{sdc}). 
We have the following  
\begin{Lemma}\label{lemma:homological} Let $\omega$ be $(\kappa,\tau)$-Diophantine  and  let $s\ge\tau$.  Then for any  $f\in {\mathcal A}^s(\T^{n-1})$ satisfying 
\[
\int_{\T^{n-1}}f(\varphi)d\varphi =0
\]
the homological equation  
\[
{\cal L}_\omega u=f\, ,\quad  \int_{\T^{n-1}}u(\varphi)d\varphi =0\, 
\]
has one and only one   solution $u\in {\mathcal A}^{s-\tau}(\T^{n-1})$. Moreover, the solution $u$  satisfies the estimate 
$$\|u\|_{s-\tau}\ \le\ \frac{1}{4\kappa}\, \|f\|_{s}.$$
\end{Lemma}
{\em Proof}. 
Comparing the  Fourier coefficients $u_k$
and  $f_k$, $0\neq k\in \Z^{n-1}$, of $u$ and $f$ respectively,  we get $f_0=0$ and 
\[
u_k\ =\ \frac{f_k}{1 -\exp (2\pi i \langle
k,\omega\rangle )}\,  , \ k\neq 0\, ,
\] 
and we set set $u_0 = 0$. Using 
(\ref{sdc}) we obtain
\[
|1 -\exp (2\pi i \langle
k,\omega\rangle )|= 2|\sin(\pi (k_n-\langle
k,\omega\rangle ))|\ge 4 \kappa |k|^{-\tau}\ge 4\kappa  (1+|k|)^{-\tau} 
\]
for $k\neq 0$, where $k_n-\langle
k,\omega\rangle \in (-1/2,1/2)$. Hence,
\[
|u_k|\ \le \ \frac{1}{4\kappa}(1+|k|)^\tau |f_k|\,  , \ k\in\Z^{n-1}\, .
\] 
Summing up we get the function $u$ and the estimate of $\|u\|_{s-\tau}$. 
In this way we obtain a unique solution $u \in {\mathcal  A}^{s-\tau} (\T^{n-1})$
of (\ref{homological}) 
normalized by $\int_{\T^{n-1}}  u(\varphi)\, d \varphi = 0$.
 \finishproof
%Notice  in the corresponding $C^k$ spaces we loose at least $[\tau]+n$ derivatives. 

The frequencies $I$ of the  quasimode we  are going to construct  will 
satisfy $I-I^0 \sim \lambda^{-1}$, 
where $\lambda^2$ are  the corresponding quasi-eigenvalues. 
For that reason we neglect systematically terms of order $O_N(|I-I^0|^N)$ and 
consider  the Taylor 
polynomials of the 
symbols  at $I=I^0$ up to certain order. The corresponding class of symbols is defined as follows. 

Let $\psi\in C_0^\infty(D)$ and  $\psi = 1$ in a  neighborhood of $I^0$. 
Fix  $\ell,\,  l \ge 2$, $s\ge 2$, $s'\ge 1$,  and a positive integer $N\ge 1$ such that  $s\ge s'$ and 
$l \ge sN + 2n $. For any bounded set
${\cal B}\subset C^\ell(\Gamma, \R)$ we denote by 
${\mathcal  A}^{l,N}_{s,s'}(\T^{n-1}\times D; {\mathcal B}; \lambda)$
the class of symbols depending  on  $K\in  {\mathcal B}$ of the form 
\begin{equation}
\left\{
\begin{array}{lcrr} 
\displaystyle a(\varphi,I, \lambda)\ =\   \sum^{N-1}_{j=0}a_{j,K}(\varphi,I)\lambda^{-j}\,
,\\[0.3cm]
\displaystyle a_{j,K}(\varphi,I)\ = \ \psi(I)\,  \sum_{|\alpha| \le N-j-1} 
(I-I^0)^\alpha a_{j,\alpha,K}(\varphi)\, 
\end{array}
\right.
\label{newsymbols}
\end{equation}
where  
$$a_{j,\alpha,K}= \partial_I^\alpha a_{j,K}(\cdot,I^0)/\alpha !\in {\mathcal  A}^{l-sj -s'|\alpha|}(\T^{n-1}), $$
for $j +|\alpha|\le N-1$, and such that the corresponding maps 
$$C^\ell(\Gamma, \R)\supset {\mathcal B} \ni 
K\to a_{j,\alpha,K}\in {\mathcal A}^{l-sj-s'|\alpha|}(\T^{n-1})$$ 
are continuous. 
We  also say that $r$ belongs to the residual class of symbols 
$\widetilde R_{N}(\T^{n-1}\times D; {\mathcal B}; \lambda)$
if 
\begin{equation}
\left\{
\begin{array}{lcrr} 
\displaystyle r(\varphi,I,\lambda)\ =\  \sum^{N-1}_{j=0}r_{j,K}(\varphi,I)\lambda^{-j}\,
,\\[0.3cm]
\displaystyle r_{j,K}(\varphi,I) = \sum_{|\alpha|=N-j} (I-I^0)^\alpha r_{j,\alpha,K}(\varphi,I)
\end{array}
\right.
\label{newreminders}
\end{equation}
where 
$r_{j,\alpha,K}\in C_0^{n}(\T^{n-1}\times D)$, 
the support of $r_{j,\alpha,K}$ is contained in a 
fixed compact subset
of $\T^{n-1}\times D$ independent of $K$,  and 
$\|r_{j,\alpha,K}\|_{C^{n}}\le C_{\cal B}$, where $C_{\cal B}>0$ does not depend on 
$K\in {\cal B}$. 
Note that the factor space  ${\mathcal  A}^{l,N}_{s,s'}/\widetilde R_{ N}$ 
does not depend   of $\psi$. Moreover, 
\begin{equation}\label{A-to-S}
{\mathcal  A}^{l,N}_{s,s'}(\T^{n-1}\times D; {\mathcal B}; \lambda) \subset S_{l,s,N}(\T^{n-1}\times D; {\mathcal B}; \lambda), 
\end{equation}
where the class of symbols $S_{l,s,N}(\T^{n-1}\times D; {\mathcal B}; \lambda)$ is introduced in Definition \ref{Def:finite-smoothness}. In particular, the corresponding $\lambda$-PDOs are uniformly bounded in $L^2$ by Lemma \ref{Lemma:L2-continuity}. The $\lambda$-PDOs  with symbols $r_{j,\alpha,K}(\varphi,I)$ are also uniformly bounded in $L^2$ by Remark \ref{Rem:reminder}.

 We shall take below $s=\tau +2$ and $s'=\tau$ which is suggested by Lemma \ref{lemma:homological}. 
%The choice of the residual class is motivated by  the proof of Proposition \ref{prop:spectral-decomposition} below. 
From now on we will  drop the index $K$ keeping in mind that the corresponding symbols depend on $K$. The following Proposition provides a normal form of the microlocal  monodromy operator $W_1(\lambda)$ given by  Proposition \ref{operator-T}. 
\begin{Prop}\label{prop:commutator}
Fix $l\ge M(\tau +2) +2n$  and $\ell>l+(n-1)/2$,  and suppose that $K$ belongs 
to a bounded subset $\mathcal B$ of $C^\ell(\Gamma,\R)$. 
Then there exists a $\lambda$-PDO $A(\lambda)$ of order $0$ 
 acting on  $C^\infty(\T^{n-1},{\LL})$ and  a
$\lambda$-FIO $W^0(\lambda)$ of the form (\ref{operator}) 
such that 
\[
W_1(\lambda)A(\lambda)\ =\ A(\lambda)W^0(\lambda) +
R^0(\lambda)\, + \, O_{\mathcal B}(|\lambda|^{-M-1})    \,  ,
\]
where 
\begin{enumerate}
\item[(1)] the full symbols of $A(\lambda)$ and of $W^0(\lambda)$ are
\[
\sigma(A)(\varphi,I,\lambda)=a_0(I) + \lambda^{-1}a^0(\varphi,I,\lambda) \ \mbox{and} \ 
\sigma(W^0)(\varphi,I,\lambda)=  p_0(I) + \lambda^{-1} p^{0}(I,\lambda)\, ,\ \mbox{where}
\] 
\item[(2)] $a_0,p_0,\in C^\infty_0(D)$ do not depend on $K$, $a_0(I)=p_0(I)=1$  in a neighborhood $D^0$ of $I^0$, 
and 
\begin{equation}
\displaystyle{
\left\{
\begin{array}{lcrr} 
\displaystyle \mbox{the symbol}\quad  a^0(\varphi,I,\lambda)=\ \psi(I)\,  \sum_{j+|\alpha| \le M-1}\, 
\lambda^{-j}(I-I^0)^\alpha a^0_{j,\alpha,K}(\varphi)\ \\[0.3cm] 
 \mbox{belongs to}\  {\mathcal  A}^{l -\tau, M}_{\tau+2,\tau}(\T^{n-1}\times D;{\mathcal B}; \lambda)\, ; 
\\[0.3cm] 
\displaystyle  p^0(I,\lambda) =   \sum_{j=0}^{M-1} 
\lambda^{-j}  p^0_{j}(I) \quad \mbox{and}\quad p^0_{j}(I)= \psi(I)\sum_{|\alpha| \le M-j-1} (I-I^0)^\alpha p^0_{j,\alpha}\,  
 \, ,
\end{array}
\right.}
\label{symbols}
\end{equation}
\item[(3)] the maps $C^\ell(\Gamma, \R)\ni K\to p^0_{j,\alpha}\in \C$, $j+|\alpha|\le M-1$,  are continuous, and 
\begin{equation}
\label{first-invariant2}
p^0_{0,0} = \frac{1}{(2\pi)^{n-1}}\, \int_{\T^{n-1}} \, w^0_0(\varphi,I^0) d\varphi\, ,
\end{equation}
where  $w^0_0$ is given by  \eqref{first-invariant1}, 
\item[(4)]
$R^0$ is a  $\lambda$-FIO of the form 
(\ref{operator}) with symbol 
\begin{equation}
r=\sum_{j=0}^{M}r_j\lambda^{-j}\in 
\widetilde R_{M+1}(\T^{n-1}\times D;{\mathcal B}, \lambda)\, ,   
          \label{reminders}
\end{equation}
\end{enumerate}
\end{Prop}
Proposition \ref{prop:commutator} will be proved in the Appendix A.2. 
The main idea of the proof  is given in the beginning of this section. 
\begin{Remark}\label{rem:birkhoff}
The operator $W^0(\lambda)$ is called a Quantum Birkhoff Normal Form (shortly QBNF) of the  monodromy operator $W(\lambda)$. 
The coefficients $p^0_{j,\alpha}$, $j+|\alpha|\le M-1$ are the corresponding Birkhoff invariants. It can be shown that the Birkhoff invariants do not depend on the choice of the operators $T(\lambda)$ and $A(\lambda)$. 
\end{Remark}
Taking into account (\ref{first-invariant1}) and (\ref{first-invariant2}) we obtain 
\begin{equation}
\label{first-invariant3}
p_0^0(I^0) \,  =\, ic \, -\, \frac{2i}{(2\pi)^{n-1}}\,
 \sum_{j=0}^{m-1} \int_{\T^{n-1}}\frac{K\circ \pi_\Gamma}{\sin
\theta}(B^{-j}\chi(\varphi, I^0))\, d\varphi \, 
=\, ic\, -\, 2i  
\sum_{j=0}^{m-1} \int_{\Lambda_j}\frac{K\circ \pi_\Gamma}{\sin \theta}\, d\mu_j\, ,  
\end{equation}
where $c$ is independent of $K$.

\subsection{Spectral decomposition of $W(\lambda)$ near $\Lambda$.}\label{subsec:spectral}

In what follows we shall find solutions  $(\lambda, v(\lambda))$ of the equation \eqref{equation} of the form $v(\lambda)=E(\lambda)^{-1}T(\lambda)A(\lambda)e(\lambda)$. In view of \eqref{conjugation} and  Proposition  \ref{prop:commutator}, $e(\lambda)$ should satisfy the equation 
\begin{equation}
\label{equation1}
e^{i\pi\vartheta/2}W^0(\lambda)e(\lambda) + e^{i\pi\vartheta/2}R^0(\lambda)e(\lambda)  = e(\lambda) + O_{\mathcal B}(|\lambda|^{-M-1})e(\lambda). 
\end{equation} 
Candidates for $e(\lambda)$ are the sections $e_k$, $k\in \Z^{n-1}$. 
Since $\lambda\in {\mathcal D}$ is complex, 
we consider an almost analytic extensions of  order $3M+3$ 
of the phase function  $\Phi$ in $I =\xi + i\eta$, $|\eta|\le C$, which is   
given  by $\Phi(x,\xi + i\eta)=L(\xi + i\eta)+ \Phi^0(x,\xi + i\eta)$, where
\[
L(\xi + i\eta) = \sum_{|\alpha|\le 3M+3} \partial_\xi^\alpha L(\xi)(i\eta)^\alpha (\alpha !)^{-1}\ \mbox{and}\
\Phi^0(x,\xi + i\eta) = \sum_{|\alpha|\le 3M+3} \partial_\xi^\alpha
\Phi^0(x,\xi)(i\eta)^\alpha (\alpha !)^{-1}\,.
\]
It is easy to see that $\overline \partial_I \Phi(x,\xi + i\eta) =
O(|\eta|^{3M+3})$ as $\eta\to 0$. Moreover, 
\begin{equation}
\Phi^0(x,\xi + i\eta) = O\left(|\xi - I^0|^{3M+3}\right)\, ,\quad |\eta|\le C. 
\label{almost-analytic}
\end{equation}
In the same way we construct 
 an almost analytic extension of order $M$ of the function 
  $\psi$,  which was used in (\ref{symbols}). 
We have   $\psi(\xi + i\eta) = 1$ in
a complex neighborhood of $I^0$ and  $\psi(\xi + i\eta) = 0$ for 
$\xi\notin D$.

\begin{Prop}  \label{prop:spectral-decomposition}
Fix  $C> 1$. Then for any  $\lambda \in {\cal D}$, $|\lambda|>1$,  and any $k\in \Z^{n-1}$,  such that $|k|\le C|\lambda|$ we  have 
\begin{equation}
\label{newoperator}
\begin{array}{rcll}
W^0(\lambda)e_k(\varphi)\, &=&\,  
\exp\left[i\lambda \Phi(\varphi,(k + \vartheta_0/4)\lambda^{-1})\right] \\[0.3cm]
&\times&\, \left(p_0 + \lambda^{-1}p^0\right)((k + \vartheta_0/4)\lambda^{-1},\lambda)\,   e_k(\varphi)\,  + \, 
O_{\cal B}(|\lambda|^{-M-1}) \, ,
\end{array}
\end{equation}
and 
\begin{equation}
R^0(\lambda)e_k(\varphi)\,  = \, 
O_{\cal B}\left(|\lambda|^{-M-1}+ |I^0 - (k+\vartheta/4)\lambda^{-1} |^{M+1}\right)\,   .
                                  \label{newrest}
\end{equation}
\end{Prop}
The proof of the proposition is given in the Appendix.

Now choosing $\lambda$ and $k$ so that $(k + \vartheta_0/4)\lambda^{-1}\sim I^0$ we can get rid of $\Phi^0(\varphi,(k + \vartheta_0/4)\lambda^{-1})$ which will allow us  to obtain quasi-eigenvalues of $W^0(\lambda)$.

\subsection{Construction of quasimodes.}\label{subsec:construction}
\subsubsection{Quantization conditions.}\label{subsec:quantization-conditions}
The index set ${\cal M}$ of the
quasimode ${\cal Q}$ we are going to construct is  defined  as
follows. Fix $d_n>0$ such that $2d_n$ is greater than the diameter of $[0,1]^{n-1}\times[0,2\pi]$ in $\R^n$.
We say that the  pair
$q=(k,k_n)\in\Z^{n-1}\times\Z$ belongs to $\cal M$ if there exists  
$\lambda=\mu^0_q\ge 1$ such that the following
quantization conditions hold:
\begin{equation}
\left|\lambda (I^0, L(I^0))\  -\ (k+\vartheta_0/4, 2\pi k_n-\pi \vartheta/2)\right|
\,  \le \, d_n\, ,
                                                     \label{quantization}
\end{equation} 
where $|\cdot|$ is the Euclidean norm in   $\R^n$. 
 We have  $(I^0, L(I^0))\neq (0,0)$ in view of
Lemma \ref{positiv-action}, hence, ${\cal M}$ is an infinite subset of $\Z^{n}$. Moreover, there is $C>0$ such that $\mu^0_q\ge C|q|$. 
Choose a sequence $\left(\mu^0_q\right)_{q\in{\mathcal M}}$ such that $\lambda = \mu^0_q\ge 1$ satisfies \eqref{quantization}. 
Fix  $C_0 > 0$ and set
\[
B(\mu^0_q):= \{\lambda \in \C\, :\ |\lambda
-\mu^0_q| \le C_0\}\, . 
\] 
Then \eqref{quantization} implies
\begin{equation}
\sup_{\lambda\in B(\mu^0_q)}\left|\lambda (I^0, L(I^0))\  -\ (k+\vartheta_0/4, 2\pi k_n-\pi \vartheta/2)\right|
\, =\, O(1)
                                                     \label{quantization1}
\end{equation}
uniformly with respect to $q\in {\mathcal M}$. Using (\ref{almost-analytic}), (\ref{newoperator}) and \eqref{quantization1} for $q\in {\cal
M}$ and $\lambda\in B(\mu^0_q)$
 we obtain
\[
W^0(\lambda)e_k\, =\, Z_q(\lambda)\,   e_k\, + \, 
O_{\cal B}(|\lambda|^{-M-1})  e_k\, ,
\]
where
\[
\begin{array}{lcrr}
Z_q(\lambda)\, =\, 
e^{ i\lambda L((k + \vartheta_0/4)\lambda^{-1}) } 
\left[1+ \lambda^{-1}p^0((k + \pi \vartheta_0/4)\lambda^{-1},\lambda)\right]\\[0.3cm] 
=\ 
\exp\, \left[i\lambda L((k + \vartheta_0/4)\lambda^{-1}) 
+  {\rm Log}\, \left(1 +  \lambda^{-1}p^0((k + \vartheta_0/4)\lambda^{-1},\lambda)\right)\right]\, , 
\end{array}
\]
where $ {\rm Log}\, z = \ln |z| + i\,  {\rm arg}\, z,\ -\pi <{\rm arg}\,
z < \pi$. On the other hand, (\ref{newrest}) and (\ref{quantization1}) imply 
\[
R(\lambda)e_k\, =\, O_{\cal B}(|\lambda|^{-M-1})\, e_k\ \mbox{for}\ \lambda\in B(\mu^0_q)\, . 
\]
Hence,
\begin{equation}
\begin{array}{lcrr}
e^{i\pi \vartheta/2}W_1(\lambda)A(\lambda)e_k\, = \, e^{i\pi \vartheta/2}A(\lambda)\left(W^0(\lambda)e_k\, + R^0(\lambda)e_k\right) \\[0.3cm]
 =\, 
e^{i\pi \vartheta/2}Z_q(\lambda)A(\lambda)e_k + O_{\cal B}(|\lambda|^{-M-1}) e_k\, .
\end{array}
\label{relation}
\end{equation}
and \eqref{equation1} reduces to the equation 
\begin{equation}\label{quasi-equation}
e^{i\pi \vartheta/2}Z_q(\lambda)\ =\ 1\ ,  \quad \lambda\in
B(\mu^0_q)\, , 
\end{equation} 
modulo $O_{\cal B}(|\lambda|^{-M-1}) $. 

\subsubsection{Solving \eqref{quasi-equation}.}
In order to solve \eqref{quasi-equation}  
we are looking for a  perturbation
$\lambda=\mu_q\in B(\mu^0_q)$ of $\mu^0_q$ such that
\[
\mu_q L\left(\frac{k + \vartheta_0/4}{\mu_q}\right) + \pi
\vartheta/2
\]
\[
+ \frac{1}{i}\, {\rm Log}\, \left[1 +\frac{1}{\mu_q}\, p^0\left(\frac{k + \vartheta_0/4}{\mu_q},\mu_q\right) \right]
\
= \  \ 2\pi k_n + 
O_{\cal B}(|\mu_q|^{-M-1})\, .
\]                 
Introduce a small parameter $\varepsilon_q =  (\mu_q^0)^{-1}$.  We
are looking for 
\begin{equation}
\label{mu_q}
\left\{
\begin{array}{rcll}
\mu_q \, &=&\, \mu_q^0 + c_{q, 0} + c_{q, 1}\varepsilon_q +
\cdots  c_{q, M}\varepsilon_q^{M}\ , \\[0.3cm] 
\zeta_q \, &=&\, I^0  +  b_{q, 0}\varepsilon_q+
\cdots  b_{q, M} \varepsilon_q^{M+1} + b_{q, M+1} \varepsilon_q^{M+2}
\end{array}
\right.
\end{equation}
such that
\begin{equation}
\label{system}
\left\{
\begin{array}{rcll}
\mu_q\zeta_q &=& k + \vartheta_0/4 \, \\[0.3cm] 
\displaystyle \mu_q L(\zeta_q) &=&  \displaystyle   2\pi k_n - \pi
\vartheta/2 - \frac{1}{i}\, {\rm Log}\, \left(1 + \mu_q^{-1}p^0(\zeta_q,\mu_q)\right) +  
O_{\cal B}(\varepsilon_q^{M+1})\, .
\end{array}
\right.
\end{equation}
Using \eqref{mu_q} we write 
\[
\mu_q\zeta_q - k - \vartheta_0/4 = \sum_{j=0}^M \varepsilon_q^j \left[ b_{q, j }+ c_{q, j} I^0 - W_{q,j}  \right] 
+ \varepsilon_q^{M+1}\left[(\varepsilon_q \mu_q) b_{q, M+1}- W_{q,M+1}  \right], 
\]
where 
\begin{equation}
\label{W}
\left\{
\begin{array}{lcrr}
\displaystyle W_{q,0} = k + \vartheta_0/4 - \mu_q^0 I^0\,  ,\ W_{q,j} = -\sum_{r+s=j-1}c_{q, r}b_{q, s }\ \mbox{for}\  1\le j \le M,\ \mbox{and}\\ 
\displaystyle  W_{q,M+1} = -\sum_{s=0}^M\sum_{r= M-s}^Mc_{q, r}b_{q, s }.
\end{array}
\right.
\end{equation}
Recall that 
\[
p^0(\zeta_q,\mu_q)= \sum_{m=0}^{M-1} 
\sum_{|\alpha|\le M-m-1}p_{m,\alpha}^0(\zeta_q-I^0)^\alpha \mu_q^{-m} 
\, . 
\]
Then 
$$
\begin{array}{lcrr}
\displaystyle{\rm Log}\, \left(1 +\mu_q^{-1}p^0(\zeta_q,\mu_q)\right) = \sum_{j=1}^{M}\frac{(-1)^{j-1}\varepsilon_q^j }{j}
\left( 1  + \varepsilon_q \sum_{r=0}^{M-1} c_{q, r}\varepsilon_q^r \right)^{-j}\\[0.3cm]
\displaystyle \times \left[ \sum_{m+|\alpha|\le M-1}\varepsilon_q^{m+|\alpha|}p_{m,\alpha}^0\left(\sum_{s=0}^{M-1}b_{q, s}\varepsilon_q^s\right)^\alpha 
\left( 1  + \varepsilon_q \sum_{r=0}^{M-1} c_{q, r}\varepsilon_q^r\right)^{-m}\right]^j  + O_{\cal B}(\varepsilon_q^{M+1}) \\[0.5cm]
\displaystyle  = \sum_{j=1}^{M} u_{q,j}\varepsilon_q^{j} + O_{\cal B}(\varepsilon_q^{M+1})\, ,
\end{array}
$$
where 
\[
u_{q,j} = p^0_{j-1,0} + \sum_{m=0}^{j-2} \sum_{|\alpha|= j-m-1} p^0_{m,\alpha}b_{q, 0}^\alpha + v_{q,j} ,
\]
and  $v_{q,j}$ is a  polynomial of $c_{q, r}$, $b_{q,r'}$ and  of 
$p_{m,\alpha}^0$, where $0 \le r,r' \le j-1$ and  $m+|\alpha| \le  j-2$. 
Using the Taylor expansion of $L(I)$ at $I^0$ up to order $M+1$ we obtain
$$
\begin{array}{lcrr}
\displaystyle
\mu_q L(\zeta_q)=  \sum_{|\alpha|\le M+1}\mu_q \frac{(\zeta_q-I^0)^\alpha }{\alpha !}\partial^\alpha L(I^0) + O_{\cal B}(\varepsilon_q^{M+1}) \\[0.3cm]
\displaystyle = \mu_q L(I^0) + \sum_{j=1}^{M} \varepsilon_q^j \left[L( I^0)c_{q, j} +  2\pi \langle \omega,b_{q, j}\rangle + u'_{q,j}\right] + O_{\cal B}(\varepsilon_q^{M+1}), 
\end{array}
$$
where  $u'_{q,j}$ are polynomials of $c_{q, r}$ and  $b_{q,r'}$,  $0 \le r,r' \le j-1$. 

We obtain from \eqref{system} 
the following linear  systems
\[
\left\{
\begin{array}{rcll}
b_{q, j} + c_{q, j} I^0   &=&\
W_{q,j} \\[0.3cm]
\displaystyle L( I^0)c_{q, j} +  2\pi \langle \omega,b_{q, j}\rangle  &=& \displaystyle   
V_{q,j}\, ,
\end{array} 
\right.
\]
for for $0 \le j\le M$, and we put
$b_{q, M+1}= (\varepsilon_q \mu_q)^{-1} W_{q,M+1}$, 
where $W_{q,j}$ is given by \eqref{W}, 
\begin{equation}
\label{V}
V_{q,j} = -\frac{1}{i}p^0_{j-1,0}  -\frac{1}{i} \sum_{m=0}^{j-2} \sum_{|\alpha|= j-m-1} p^0_{m,\alpha}b_{q, 0}^\alpha + V'_{q,j}\, , 
\end{equation} 
and $V'_{q,j}$ is a  polynomial of $c_{q, r}$, $b_{q,r'}$ and  of 
$p_{m,\alpha}^0$, where $0 \le r,r' \le j-1$ and  $m+|\alpha| \le  j-2$. In particular, 
\begin{equation}
\label{V0}
V_{q,1} = -\frac{1}{i}p^0_{0,0}  -2\pi \langle \omega, b_{q,0}\rangle 
-\frac{1}{2}\langle \nabla^2 L(I^0)b_{q,0},b_{q,0}\rangle. 
\end{equation}
In view of  Lemma \ref{positiv-action},  the corresponding determinant is
$$
D(I^0):= L(I^0) -  2\pi \langle I^0, \omega\rangle =  \frac{1}{(2\pi)^{n-1}}\, 
\int_{\Lambda} A(\varrho)\,  d\mu > 0\, ,
$$
and we obtain a unique
solution $(c_{q, j},b_{q ,j})$, $0\le j \le M-1$. 
More precisely, 
\begin{equation}
\label{quasi-coefficients}
\left\{
\begin{array}{lcrr}
 c_{q, j} = D(I^0)^{-1}\left[ V_{q,j} - 
 2\pi\langle \omega,W_{q, j}\rangle\right]\\[0.3cm]
b_{q, j} = W_{q,j} -c_{q, j} I^0 .
\end{array}
\right.
\end{equation} 
We have  
\[
\left\{
\begin{array}{lcrr}
W_{q,0}=k + \vartheta_0/4 - \mu_q^0 I^0 = O(1), \\[0.3cm]
V_{q,0} = 2\pi k_n -\pi \vartheta/2- \mu_q^0 L(I^0)= O(1)\, ,\ q\in \cal M\, ,
\end{array}
\right.
\]
in view  of (\ref{quantization}). 
Hence, $b_{q,0}$ and $c_{q,0}$ are uniformly bounded and they do not depend on $K$. 
 By recurrence we prove that $b_{q, j}$ and $c_{q, j}$ are
continuous with respect to $K$ and uniformly bounded with respect to $q\in {\cal M}$ and $K\in {\cal B}$. To evaluate $b_{q, M+1}$ observe that $\varepsilon_q \mu_q = 1 + O(\varepsilon_q)$. 
For  $j=1$ we obtain  from  \eqref{W} and \eqref{V0} 
\begin{equation}
\label{system1}
\left\{
\begin{array}{lcrr}
\displaystyle V_{q,1}= -\frac{1}{i}p^0_{0,0}  -2\pi \langle \omega, b_{q,0}\rangle 
-\frac{1}{2}\langle \nabla^2 L(I^0)b_{q,0},b_{q,0}\rangle , \\[0.3cm]
W_{q,1}=-c_{q,0}b_{q,0}
\end{array}
\right.
\end{equation}
and we get from \eqref{first-invariant3} that
\[
\displaystyle c_{q, 1}\, =\,    c_{q,1}' \, + \, \frac{2}
{\int_{\Lambda_1} A(\varrho) d \mu_1}\
 \sum_{j=0}^{m-1} \int_{\Lambda_j}\frac{K\circ \pi_\Gamma}{\sin
\theta}\, d\mu_j\, , 
\] 
where $c_{q,1}'$ does not depend on $K$.  

\subsubsection{Quasimodes}\label{susec:quasi}
For each  $q =(k,k_n)\in {\cal M}$ we set
\[
v_q^0\ :=\ E(\lambda)^{-1}T(\mu_q)A(\mu_q)e_k \quad {\rm and} \quad 
u_q^0\ :=\  G(\mu_q)v_q^0\ =\ 
G(\mu_q)E(\lambda)^{-1}T(\mu_q)A(\mu_q)e_k\, .
\]
Recall that  $T(\lambda)$ is defined by Proposition \ref{operator-T},  $A(\lambda)$ by Proposition \ref{prop:commutator}, and $G(\lambda)$ by \eqref{operator-G}. 
Using (\ref{relation}), we obtain 
\begin{equation}
( W(\mu_q) -\,  {\rm Id}\, )v_q^0 = 
O_{\cal B}(|\mu_q|^{-M-1})\, v_q^0\, ,
\label{theequation}
\end{equation}
and by Proposition \ref{prop:monodromy1}  we get
\[
\left|
\begin{array}{rcll}  
\left(\, \Delta \, -\, \mu_q^2 \, \right)\, u_q^0\ &=&  
O_{\cal B}(|\mu_q|^{-M})\,u_q^0\quad 
{\rm in}\  L^2(X)\, ,\\[0.3cm]
{\mathcal N} u_q^0 |_\Gamma \ &=&\ O_{\cal B}(|\mu_q|^{-M})\, u_q^0 \quad 
{\rm in}\  L^2(\Gamma)\, . 
\end{array} 
\right.
\]
It remains to estimate the $L^2$-norm of $u_q^0$ and to satisfy the boundary conditions in (\ref{thequasimode}) adding to $u_q^0$  a term of magnitude  $O_{\cal B}(|\mu_q|^{-M})$. 
\begin{Lemma}\label{lemma:estimates}
There is $C >1 $ such that 
\[
C^{-1} \le \|u_q^0\|_{L^2(X)} \le C
\]
for any $q\in {\cal M}$.
\end{Lemma}
{\em Proof}. 
The operator $A(\lambda)$ is  of the form $A(\lambda)= A_0(\lambda) + \lambda^{-1}A^0(\lambda)$, where $A_0(\lambda)$ is a classical $\lambda$-PDO of order $0$ of symbol $a_0$ and it is elliptic on $\T^{n-1}\times D^0$, while $A^0(\lambda)$ has a symbol $a^0$ in $S_{l -\tau, \tau+2,M}(\T^{n-1}\times D;{\mathcal B}; \lambda)$
according to Proposition \ref{prop:commutator} and  \eqref{A-to-S}. By Lemma \ref{Lemma:L2-continuity} the family of operators $A^0(\lambda)$ is uniformly bounded in $L^2$ with respect to $\lambda\in {\mathcal D}$ and $K\in {\mathcal B}$ since $l \ge (\tau +2)M + 2n$. 
Moreover,   $E(\lambda)^{-1}T(\lambda)$ is a classical $\lambda$-FIO of order 0 the canonical
relation of which is the graph of a canonical transformation $\chi$. Moreover, it is  elliptic over  a neighborhood of 
the frequency set of $\{e_k:\, q=(k,k_n)\in {\cal M}\}$. Then there is $C>1$ such that 
\[
\forall\,  q\in  {\cal M}\, ,\quad  C^{-1} \ \le \ \|v_q^0\|_{L^2(\Gamma)}\ \le C\, . 
\]
We write $G(\lambda) = \sum_{j=0}^{m-1} \tilde G_j(\lambda)$, where $\tilde G_0(\lambda) = H_0(\lambda)\psi_0(\lambda)$ 
and 
\[
\tilde G_j(\lambda) = H_{j}(\lambda)\Pi_{s=0}^{j-1}
\left( Q_{s+1}(\lambda) G_{s}(\lambda)\right)\psi_0(\lambda) 
\]
for $1\le j\le m-1$ and $m\ge 2$ (see \eqref{operator-G}). Then using Proposition \ref{Lemma:continuity} for $H_j(\lambda): L^2(\Gamma) \to L^2(\widetilde X)$, $0\le j\le m-1$,  we find  $C>0$ such that 
 $C^{-1} \le \|u_q^0\| \le C$ for  any $q\in {\cal M}$, 
which completes the proof of the lemma. \finishproof

To satisfy the boundary conditions in (\ref{thequasimode}) we choose $\kappa\in C^\infty_0((-\varepsilon,\varepsilon))$, $\kappa\ge 0$ such 
that $\kappa = 1$ in a neighborhood of $0$,  and set 
\[
\tilde u_q (x',t)=t\kappa(t)\left(\frac{\partial u_q^0}{\partial t}(x',0) -K(x')u_q^0(x',0)\right)\, ,
\] 
where $(x',t)\in \Gamma\times \R$ are normal coordinates to the boundary $\Gamma$. Then 
$\|\tilde u_q\|_{L^2(X)} \le C'_{\cal B} \mu_q^{-M}$ and $u_q^0 -\tilde u_q$ satisfies the boundary 
condition. We put now
\[ 
w_q = u_q^0 -\tilde u_q 
+(\Delta_{K} - \mu_q^2 - i)^{-1}(\Delta - \mu_q^2)\tilde u_q \, .
\]
Notice that the $L^2$-norm of the second and the third terms can be estimated by $O_{\cal B} (\mu_q^{-M})$
and $$\|(\Delta - \mu_q^2)w_q\|_{L^2(X)} \le C_{\cal B} \mu_q^{-M}.$$  Moreover, $w_q$  belongs to the domain of 
definition of $\Delta_K$. 
Normalizing 
$ u_q = w_q\|w_q\|^{-1}$ and using Lemma \ref{lemma:estimates}, we obtain  a quasimode $(\mu_q, u_q)$ of
order $M$. 
Next  we  show that $\mu_q$ can be chosen
real-valued. 
 Applying Green's
formula we get
\[
|\mu_q^2 - \overline{\mu_q}\, ^2|\ = \ 
|\langle \mu_q^2 u_q, u_q \rangle\,  -\,  
\langle  u_q, \mu_q^2 u_q \rangle|\ = \ 
 O_{\cal B}(|\mu_q|^{-M})\, ,
\]
which allows us to take $\mu_q$ in $\R$. Choosing $|q|\gg 1$ we
can suppose that $\mu_q$ is positive. Notice that $K$ should be in
$C^\ell(\Gamma, \R)$ with $l > M(\tau + 2)+2n+(n-1)/2$. 

\section{Spectral rigidity in the presence of a $(\Z/2\Z)^2$-group of symmetries.}\label{sec:rigidity}
\setcounter{equation}{0}
\subsection{ Spectral rigidity for a bouncing ball geodesic.}\label{subsec:bouncing-ball}
Let $\gamma$ be a closed broken geodesic in $(X,g)$, ${\rm dim}\, X \ge 2$, 
with $m\ge 2$ vertices $x_j$, $0\le j \le m-1$. Denote by $\varrho_j=(x_j,\xi_j) = B^{j}(\rho_0)$, $0\le j\le m-1$,  the corresponding periodic trajectory of $B$. Then $\varrho_0$ is a fixed point of the symplectic map $P=B^m$. Recall that $\gamma$ is called elliptic if $\varrho_0$ is an elliptic fixed point of $P$ which means that the spectrum of $dP(\varrho_0):T_{\varrho_0}\Gamma \to T_{\varrho_0}\Gamma$ lies on the unit circle of the complex plane and it consist of distinct eigenvalues different from one, i.e.
\[
{\rm Spec}\, (DP(\varrho_0))=\{e^{\pm i2\pi\alpha_j}:\,  1\le j\le n-1\}, \ \mbox{where}\ 0 <\alpha_1< \cdots <\alpha_{n-1}\le 1/2\, .
\] 
Set $\alpha=(\alpha_1,\ldots,\alpha_{n-1})$. It is said that $\gamma$ admits no resonances of order less or equal to $4$ if the scalar product $\langle \alpha,k\rangle$ is not  integer for any integer vector $k=(k_1,\ldots,k_{n-1})\in \Z^{n-1}$ different from $0$ and  such that $|k_1|+\cdots+|k_{n-1}|\le 4$. 
In this case $P$ admits a Birkhoff normal form 
\[
P(\theta,r) = (\theta +\nabla B(r) + 0(|r|^{3/2}), r +0(|r|^{2})),\ (\theta,r)\in \T^{n-1}\times\R_+^{n-1},
\]
where $(r,\theta)$ are suitable polar symplectic coordinates in a neighborhood of $\varrho_0$ such that $r(\varrho_0)=0$,
and $B(r) = \langle\alpha,r\rangle +\langle A r,r\rangle/2$, where $A$ is a symmetric matrix. The Birkhoff normal form of $P$ is non-degenerate if $\det A\neq 0$. In this case, applying the KAM theorem we get a family of invariant tori $\Lambda_\omega$, called KAM tori,  satisfying (H$_2$) and having Diophantine vectors of rotations $\omega\in \Xi$ (see \cite{La}, Theorem 13.6). Moreover, for any neighborhood $U$ of $\varrho_0$ in $B^\ast\Gamma$ the union of the KAM tori lying in $U$ has a positive Lebesgue measure. Now Theorem \ref{main} applies to any single torus of the family. 

\vspace{0.3cm}
\noindent
{\em Proof of Corollary \ref{Coro:bouncing-ball}.} 
Consider in more details the case of a bouncing ball trajectory in a two-dimensional billiard table ($n=2$ and  $m=2$). Denote the restrictions of the two involutions to $\Gamma\cap U$ also by $J_1$ and $J_2$ and
 by $\tilde J_j:T^\ast(\Gamma\cap U)\to T^\ast(\Gamma\cap U)$ the corresponding lifts.  In this case $\Gamma\cap U$ has two connected components $\Gamma_j$, $j=1,2$,  and $J_1(\Gamma_j)=\Gamma_j$ while $J_2(\Gamma_1)=\Gamma_2$. Since  $J_1$ and $J_2$ act as isometries and commute with each other, using the definition of $B$ in Sect. 2.1,  we obtain that the involutions $\tilde J_j$, $j=1,2$,  commute with each other and also with $B$. 
 
 For any $\omega\in \Xi$ we set  $\Lambda_\omega^1 =\Lambda_\omega$ 
and $\Lambda_\omega^2 =B(\Lambda_\omega)$. 
Then $\tilde J_1(\Lambda_\omega^j)$, $j=1,2$, are also invariant circles of $P= B^2$ of rotation number $\omega$  and  
 $\Lambda_\omega^j =\tilde J_1(\Lambda_\omega^j)$ for $j=1,2$, while $\Lambda_\omega^2 =\tilde J_2(\Lambda_\omega^1)$.
To prove it we use the following argument. Since $\mbox{dim}\, T^\ast\Gamma_j=2$ the KAM circle    $\Lambda_\omega^j$ divides $T^\ast\Gamma_j$ into two connected components,  and it contains the elliptic fixed point $\varrho_j=(x_j,0)$ of $P$ in its interior $D_j$. Moreover, $\tilde J_j(\varrho_j)=\varrho_j$,  hence, $\tilde J_1(\Lambda_\omega^j)$ contains 
$\varrho_j$ in its interior $\tilde J_1(D_j)$  as well. 
On the other hand, $\tilde J_1$ preserves the volume form of $T^\ast\Gamma_1$, hence, $\Lambda_\omega^j$ intersects  $\tilde J_1(\Lambda_\omega^j)$. This implies  $\Lambda_\omega^j =\tilde J_1(\Lambda_\omega^j)$, since $P$  acts  transitively on both of them. In the same way we prove that  $\Lambda_\omega^2 =\tilde J_2(\Lambda_\omega^1)$.

For any $K\in C(\Gamma)$ we have 
\[
\int_{\Lambda_\omega^1} \frac{K\circ\pi_\Gamma}{\sin \theta}\, d\mu_1 = \int_{\Lambda_\omega^1} (\frac{K\circ\pi_\Gamma}{\sin \theta}\circ\tilde J_1)\, d\mu_1\, .
\]
Since
 $(B\tilde J_2)^\ast d\mu_j=d\mu_j$,  
 %$(BJ_2)^\ast \sin\theta=\sin\theta$,
we have as well 
\[
\int_{\Lambda_\omega^2} \frac{K\circ\pi_\Gamma}{\sin \theta}\, d\mu_2 =\int_{\Lambda_\omega^1} (\frac{K\circ\pi_\Gamma}{\sin \theta}\circ B)\, d\mu_1 = \int_{\Lambda_\omega^1} (\frac{K\circ\pi_\Gamma}{\sin \theta}\circ \tilde J_2)\, d\mu_1\, .
\]
On the other hand, $d\mu_j$ and $\sin\theta$ are invariant with respect to $J_1$, 
and we obtain 
\[
\int_{\Lambda_\omega^1} \frac{K^{\#}\circ\pi_\Gamma}{\sin \theta}\, d\mu_1 = \frac{1}{2}\left(\int_{\Lambda_\omega^1} \frac{K\circ\pi_\Gamma}{\sin \theta}\, d\mu_1 + \int_{\Lambda_\omega^2} \frac{K\circ\pi_\Gamma}{\sin \theta}\, d\mu_2\right).
\]
Now Theorem \ref{main} implies that
\begin{equation}
\int_{\Lambda_\omega} \frac{(K^{\#}_t-K^{\#}_0)\circ\pi_\Gamma}{\sin \theta}\, d\mu = 0
\label{theequality}
\end{equation}
for any $t\in [0,1]$. Parametrize $\Gamma_1$ by the arclength $s\in [-a,a]$ so that $s(x_1)=0$. Then  $J_1(s)=-s$ for any $s$.  Fix $t\in (0,1]$ and set $f=K_t^{\#}-K_0^{\#}$. For any invariant circle  $\Lambda_\omega$, $\omega\in \Xi$, there is $s_\omega >0$ such that    $\pi_\Gamma(\Lambda_\omega) = [-s_\omega,s_\omega]$.  
We are going to show that there exists an infinite sequence $(y_j)_{j\in \N} \subset (0,b)$ such that $\lim y_j=0$ and  $f(y_j)=0$.
Indeed, suppose that   $f(s)\neq 0$ in $(0,b)$ for some $b>0$. Take $\omega\in \Xi$  such that $s_\omega<b$. 
Since $f(s)$ is {\em even} it does not change its sign in the interval 
$[-s_\omega,s_\omega]$. Moreover, $\sin \theta>0$ since $\theta\in (0,\pi/2]$ on the interior of $B^\ast\Gamma$. 
Hence, $(f\circ\pi_\Gamma)/\sin\theta$ does not change its sign and it is non null on  $\Lambda_\omega$, which contradicts (\ref{theequality}). This proves the existence of an infinite sequence $\{y_j\}_{j\in\N}$ such that $f(y_j)=0$ and $y_j\neq x_1$ for any $j\in\N$ and $\lim y_j=x_1$.  Now there exists an infinite sequence $(y_j')_{j\in \N} \subset (0,b)$ such that $y_j\le y_j'\le y_{j+1}$ and $\frac{df}{ds}(y_j')=0$, and so on. This implies that the Taylor polynomials of $f$ of  order less then $[\ell]+1$, $\ell\le +\infty$,  vanish at $s=0$,  which proves the assertion. \finishproof

\subsection{Spectral rigidity for Liouville billiard tables}\label{subsec:L.B.T.}
We recall from \cite{PT1} the definition of Liouville billiard tables of dimension two.
We consider two even functions $f\in C^{\infty}(\R)$, $f(x+2\pi)=f(x)$,  
and $q\in C^{\infty}([-N,N])$, $N>0$, such that  
\begin{itemize}  
\item[(i)]
 $f>0$ if $x\notin  \pi {\Z}$, and        
$f(0)=f(\pi)=0$, $f''(0)>0$;  
\item[(ii)]$q<0$ if $y\ne 0$,  $q(0)=0$ and $q^{''}(0)<0$;  
\item[(iii)] $f^{(2k)}(\pi l)=(-1)^kq^{(2k)}(0)$,  $l=0,1$,    
for every natural $k\in{\N}$.  
\end{itemize}  
 
Consider the quadratic forms  
\[
\begin{array}{ccc}    
dg^2&=&(f(x)-q(y))(dx^2+dy^2)\\  
dI^2&=&(f(x)-q(y))(q(y)dx^2+f(x)dy^2) 
\end{array}
\]
defined on the cylinder $C= \T^1\times[-N,N]$.  
  
The involution $\sigma_0 : (x,y)\mapsto(-x,-y)$ induces an involution  
of the cylinder $C$, that will be denoted by   
$\sigma_0$ as well. We  identify the points $m$ and $\sigma_0(m)$  
on the cylinder and denote by $\widetilde C := C/\sigma_0$ the topological quotient  
space. Let $\sigma : C\to\widetilde C$ be the corresponding projection.  
A point $x\in C$   
is called {\em singular} if $\sigma^{-1}(\sigma(x))=x$, otherwise it is   
a {\em regular} point of $\sigma$. Obviously, the  singular points are   
$F_1=\sigma(0,0)$ and $F_1=\sigma(1/2,0)$.   
It is shown in \cite{PT1} that   
 the quotient space $\widetilde C$ is homeomorphic to  
the unit disk ${\D}^2$ in $\R^2$ and that  
there exist a unique differential structure on $\widetilde C$ such that the  
projection $\sigma: C\to\widetilde C$ is a smooth map, $\sigma$ is a local diffeomorphism  
in the regular points,  and the push-forward $\sigma_*g$ gives a smooth  
Riemannian metric while  $\sigma_*I$ is a smooth integral of the corresponding billiard 
flow on it. We denote by $X$ the space $\widetilde C$ provided with that differentiable 
structure and call $(X, \sigma_*g)$ a Liouville billiard table. 
 Any Liouville billiard table  possesses the string property which means 
that any broken geodesic starting from the singular point $F_1\, (F_2)$ passes through 
$F_2\, (F_1)$ after the first reflection at the boundary and the sum of distances from any
point of $\Gamma$ to $F_1$ and $F_2$ is constant. 

 We impose the following additional  conditions:  
\begin{itemize}  
\item[(iv)] the boundary $\Gamma$ of $X$  is   
locally geodesically convex which amounts  to $q'(N)<0$;   
\item[(v)] $f(x)=f(\pi-x)$ for any $x$ and $f$ is strictly increasing   
on the interval $[0,\pi/2]$;   
\end{itemize}    
Liouville billiard tables   
satisfying  the conditions above    
will be  said to be  {\em of classical type}.  
One of the consequences of the last condition  is that there is a group  
$I(X)\cong{\Z}_2\times{\Z}_2$ acting  on $(X,g)$ by isometries. It is generated by the
involutions $\sigma_1$ and $\sigma_2$ defined by $\sigma_1(x,y)= (x,-y)$ and   
$\sigma_2(x,y) = (\pi - x,y)$. We point out that in contrast to \cite{PT1} we do not assume 
$f$ and $q$ to be analytic. Examples of Liouville billiard tables of classical type on surfaces 
of constant curvature and quadrics are 
provided in \cite{PT1}. The only Liouville billiard table in $\R^2$ is the interior of the 
ellipse because of the string property. 

\vspace{0.3cm}
{\em Proof of \ref{Coro:liouville}}. 
 A first integral of $B$ in $B^\ast \Gamma$
is the function ${\cal I}(x,\xi) = f(x)-\xi^2$ the regular values $h$ of which 
belong to  $(q(N),0)\cup (0,f(\pi/2))$ 
(see \cite{PT1}, Lemma 4.1 and Proposition 4.2).
Each regular level set $L_h$ consists of two connected circles $\Lambda^{\pm}(h)$ 
which are invariant with respect to $B$ for $h\in (q(N),0)$ and to $B^2$ for $h\in (0,f(1/4))$.
The Leray form on $L_h$ is 
\[
 \lambda_h \ =\ 
 \left\{   
 \begin{array}{ccc} 
 \frac{dx}{\sqrt{f(x)-h}},\ \xi > 0\, ,\\  
-\frac{dx}{\sqrt{f(x)-h}}, \ \xi < 0\, .
\end{array}
\right.
\]
Given a continuous  function $K$ on $\Gamma$ we
consider the corresponding ``Radon transform'' assigning to each circle $\Lambda^{\pm}(h)$ the integral
\[
R_K(\Lambda^{\pm}(h)) = \int_{\Lambda^{\pm}(h)}\,  \frac{K\circ \pi_\Gamma}{\sin\theta}\,  \lambda_h\, .
\]
We have
\[
\sin \theta \ =\ \sqrt{\frac{h-q(N)}{f(x)-q(N)}}\ ,
\]
hence,
\[ 
R_{K}(\Lambda^{\pm}(h))  = \pm \frac{1}{\sqrt{h-q(N)}}  
\int\limits_0^{2\pi}\, K(x)\, \sqrt{\frac{f(x)-q(N)}{f(x)-h}}\;dx  \ ,
\ h\in (q(N),0)\, ,
\]  
Fix the exponent $\tau >1$ in the small denominator condition (\ref{sdc}) sufficiently small so that  
$\ell>([2d] +1)(\tau +2) + 4 + 1/2$. For any $\kappa>0$ denote by $\Omega_\kappa^\tau$ the set of all $\omega\in\R^{n-1}$ satisfying  (\ref{sdc}). 
Then applying  
 Theorem \ref{main} for $n=2$ we obtain 
\begin{equation}
\forall\, t\in [0,1]\, ,\quad R_{K_t}(\Lambda^{\pm}(h))\ =\ R_{K_0}(\Lambda^{\pm}(h))
                       \label{radon}
\end{equation}                        
for each regular value $h$ such that the corresponding 
frequency $\omega(h)$ belongs to $\Omega_\kappa^\tau$. 
On the other hand, the union $\mathcal A =\cup_{\kappa>0}\Omega_\kappa^\tau$ 
has full Lebesgue measure in $\R$. 
\begin{Lemma}
There is $\varepsilon >0$ such that the set of all regular values $h\in (q(N), q(N) +\varepsilon)$,
the corresponding frequencies $\omega(h)$ of which belong to $\cal A$, is dense in $(q(N), q(N) +\varepsilon)$.
\end{Lemma}
The proof of the Lemma follows immediately from Proposition 4.4 \cite{PT1}, which claims that
the rotation function $\rho^-(h):=\omega(h)$ is strictly increasing and smooth  in an interval 
$(q(N), q(N) +\varepsilon)$. Then  $\omega^{-1}(\cal A)$ is dense in that interval.

As the function $R_{G_t}(\Lambda^{\pm}(h))$ is analytic in $h\in (q(N),0)$, using the Lemma we obtain 
(\ref{radon}) for any $h\in (q(N), 0)$. 
 Since $K_t$, $t=0,1$, are invariant 
with respect to the action of $I(X)$, this  implies $K_0\equiv K_1$ as in \cite{PT1}. 
\finishproof

Spectral rigidity for higher dimensional Liouville billiard tables will be obtained 
in \cite{PT2}. We point out that we do not need analyticity and 
the billiard tables we consider are supposed to
be smooth only. 

\section{Concluding remarks}\label{sec:remarks}
\subsection{Spectral invariants for continuous deformations of  potentials}\label{subsec:potential}
\setcounter{equation}{0}
Let $V_t$, $t\in [0,1]$, be  a continuous family of $C^\ell$ 
real-valued potentials in $X$, $\ell \in \N$, which means that  
the map $[0,1]\ni t\mapsto V_t$  is continuous
in  $C^\ell(X, \R)$. Denote by $\Delta_t$ the selfadjoint operators
$\Delta +V_t$ in $L^2(X)$ with Dirichlet  boundary conditions on $\Gamma$. 
Consider the corresponding spectral problem
\[ 
\left\{  
\begin{array}{rcll}  
\Delta\,  u\ + \ V_t u\ &=& \ \la\,  u \, \quad \mbox{in}\ $X$\, , \\  
u|_{\Gamma}\ &=&\ 0 \quad \mbox{in}\ \Gamma\ , 
\end{array}  
\right.  
\]          
We suppose as above that there is a Kronecker torus $\Lambda$ 
of $P=B^m$ satisfying $(H_3)$. Without loss of generality we can assume that $B^j(\Lambda)\neq \Lambda$ for any $0<j<m$. Consider the ``flow-out'' ${\mathcal T} = \{\exp(sX_{h})(\pi^+(x,\xi)): (x,\xi)\in \Lambda\}$ of $\Lambda$ with respect to the ``broken bicharacteristic flow'', where  $\pi^+:T^\ast \Gamma \to T^\ast X|_{\Gamma}$ is defined in Sect. \ref{subsec:billiard-ball}.  As in  Sect. \ref{subsec:BNF} we identify ${\mathcal T}$ with a smooth Lagrangian torus in $\widetilde{T^\ast X}$.
We  provide  ${\mathcal T}$ with  ``coordinates'' $(\varrho,t)$, $\varrho\in\Lambda$, $0\le t\le T(\varrho)=2A(\varrho)$,  using the map $(\varrho,t)\to \exp(tX_{\tilde h})(\pi^+(\varrho))$ and introduce a natural measure $d\tilde\mu=d\mu dt$ on it. Recall that $T(x,\xi)= 2A(x,\xi)$ is the corresponding return time function. The measure $d\tilde\mu$ can be related also with the mapping cylinder construction introduced in \cite{F-G}. Moreover, $d\tilde\mu$ is the unique measure on ${\mathcal T}$ which is  invariant with respect to the flow of $\tilde h$ and such that $\mbox{vol}({\mathcal T})= \int_\Lambda T(\varrho)d\mu(\varrho)$. Let $\pi_X: \widetilde{T^\ast X}\to X$ be the canonical projection. 

Fix $\ell > ([2d]+1)(\tau + 2) +2n + (n-1)/2$,  where $\tau > n-1$ is the exponent in the small denominator condition (\ref{sdc}). 
\begin{Theorem}\label{main-potential}
Let $\Lambda$ be a  Kronecker torus of the billiard ball map with a  $(\kappa,\tau)$-Diophantine  
 vector of rotation.   
Let $V_t$, $t\in [0,1]$, be a continuous family of real-valued potentials
in $C^\ell(X,\R)$  
such that $\Delta_t$ satisfy 
the  isospectral condition 
 $(H_1)-(H_2)$ for the Dirichlet problem. 
Then
\[
\forall\, t\in [0,1],\quad \displaystyle  \int_{\mathcal T}\, V_t\circ\pi_X d\tilde\mu \ =\
\int_{\mathcal T}\, V_0\circ\pi_X d\tilde\mu.  
\]
\end{Theorem}
To prove the theorem we construct as in Theorem \ref{quasimodes} 
 a continuous family of quasimodes 
 \[
 (\mu_q(t),u_q(t))_{q\in \mathcal M}\ ,\ \mathcal M \subset \Z^{n}\, ,
 \]  
 of
$\Delta_t$ of order $M$  
 such that 
$$
\mu_q(t) = \mu_q^0 +c_{q,0} + c_{q,1}(t) (\mu_q^0)^{-1}  + \cdots + 
c_{q,N}(t) (\mu_q^0)^{-M}
$$
where $\mu_q^0$  and $c_{q,0}$ are   independent of $t$, $\mu_q^0 \ge C|q|$, $C>0$, 
and $c_{q,j}(t)$ is continuous in $t\in [0,1]$. Moreover,  
$$
c_{q,1}(t)=c'_{q,1} + c''_{1} \int_{\mathcal T}\, V_t\circ\pi_X d\tilde\mu\, , 
$$ 
$c'_{q,1}$ is  independent of $t$, and 
$$
c''_{1} =    2 \, \left(
\int_{\Lambda} A(\varrho)\,  d\mu \right)^{\, -1}\, = \, \frac{4}{\mbox{vol}({\mathcal T})}\, . 
$$
To construct the quasimodes we consider 
 for each  $j=0,\ldots,m-1$ a microlocal outgoing parametrix 
$\widetilde H_j:C^\infty(\Gamma)\rightarrow  C^\infty( \widetilde X )$ 
of the Dirichlet problem for $\Delta - \lambda^2 - V$  satisfying
$$
(\Delta-\lambda^2 - V_t)\widetilde H_j(\lambda) = O_M(|\lambda|^{-M })\ \mbox{in}\
\widetilde X\, . 
$$
We are looking for  $\widetilde H_j(\lambda)$  
of the form  $\widetilde H_j(\lambda) = H_j(\lambda) +
\lambda^{-1}H_{j,t}^0(\lambda) $, where $H_{j}(\lambda)$ is introduced in Sect. \ref{subsec:reduction} and  $H_{j,t}^0(\lambda)$ is a $\lambda$-FIO of order
$1/4$ having the same canonical relation as $H_j(\lambda)$ and such that
$$
(\Delta-\lambda^2 - V_t) H_{j,t}^0(\lambda) - V_t H_j(\lambda)
= O_M(|\lambda|^{-M})\ \mbox{in}\
\widetilde X\, . 
$$
The principal symbol $p_j^0(x,\xi)$ of $H_{j,t}^0(\lambda)$ satisfies the equation
$\{h,p_j^0\} = iV_t$ in $T^\ast X$.  Taking into account the boundary values
at  $U_j$ we
get 
$$
p_j^0(\varrho,s) = i\int_0^s V_t(\exp(uX_h)(\varrho))\, du\ ,\quad \varrho
\in U_j\, .
$$
As in Sect. \ref{subsec:reduction} we get
\begin{equation}
\imath_\Gamma^\ast \widetilde H_j(\lambda) = \Psi_j(\lambda) + 
\widetilde G_j(\lambda)
+ O_M(|\lambda|^{-M})
\, ,  
                     \label{boundary-trace2}
\end{equation}
where 
$$
\widetilde G_j(\lambda) := E(\lambda)(G_j^0(\lambda) + \lambda^{-1}
\widetilde G_j^0(\lambda))E(\lambda)^{-1}
$$ 
is a $\lambda$-FIO  
the canonical relation of which 
is just the graph of the restriction of the 
billiard ball map $B: U_j \to U_{j+1}$. 
Recall that   the principal symbol of $G_j^0(\lambda)$ is given by \eqref{principal-symbol-G}.
Moreover,  
the  principal symbol of  $\widetilde G_j^0(\lambda)$   is equal to 
\[
\exp\left(i\pi\vartheta_j/4\right)\exp (i\lambda A_j(\varrho))p_j^0(\varrho,T_j(\varrho)), \quad \varrho\in U_j,
\]  
where $T_j(\varrho)$ is defined in Sect. \ref{subsec:BNF} as the first return time (at $\Gamma$) for the flow starting at $U_j$. The corresponding parametrix has the form
\[
G(\lambda) = H_0(\lambda)\psi_0(\lambda) + 
\sum_{k=2}^{m} (-1)^{k-1}
H_{k-1}(\lambda)\left(\Pi_{j=0}^{k-2}\widetilde G_{j}(\lambda)\right)\psi_0(\lambda) \,  ,
\]
and we get the equation
\begin{equation}
(  \widetilde W(\lambda) -  \mbox{Id}\, )\psi_0(\lambda)v \ = \ 
O_{\mathcal B}(|\lambda| ^{-M-1})v\, ,
                               \label{equation2}
\end{equation}
where
\[
 \widetilde W(\lambda):= (-1)^{m-1}
\Pi_{j=0}^{m-1}  \left(\psi_{j+1}(\lambda) \widetilde G_j(\lambda)\right)\, .  
\] 
As in Proposition \ref{prop:monodromy2} we get
\[
  \widetilde W^0(\lambda):= E(\lambda)\widetilde W(\lambda)E(\lambda)^{-1}= (-1)^{m-1}
\psi_m(\lambda) \left(Q^0(\lambda) + \lambda^{-1}Q^1(\lambda)\right)
S(\lambda)\psi_0(\lambda)   + O_{\mathcal B}(\lambda^{-M-1}) \, .
\]
Here, $Q^0(\lambda)$ is a classical  $\lambda$-PDOs on $\Gamma$ with a $C^\infty$ symbol
independent of $V_t$ and  with  principal symbol
 $1$ in a neighborhood of $\Lambda$, and 
 $ Q^1\in {\rm PDO}_{l,2, M}(\Gamma;{\mathcal B};\lambda)$.  
By  Egorov's theorem (see Lemma \ref{Lemma:commutator}) the principal
symbol of $Q^1(\lambda)$ is 
$$
\sigma_0(Q^1)(\varrho) = 2i \sum_{j=0}^{m-1} p_j^0(B^{j}(\varrho),T_j(B^{j}(\varrho))) = 2i \int_{\mathcal T}\, V_0\circ\pi_X d\tilde\mu 
$$
in  $P(U_0)$. 
The operator  $  S(\lambda)$ (see Sect. \ref{subsec:reduction}) does not depend on $V_t$, 
and it is a classical $\lambda$-FIO 
of order 0 with a
large parameter $\lambda\in {\cal D}$. Its canonical relation is the graph of $P:U_0\to U_m$, and its principal symbol is given by 
$\exp\left(i\pi\vartheta/4\right)\exp (i\lambda A(\varrho))$, $\varrho\in U_0$. 
Arguing as in Sect. \ref{sec:construction} we complete the
construction of the quasimodes. The  invariants obtained in Theorem \ref{main-potential} 
are related to the first Birkhoff invariant of the monodromy operator.

In the same way one can deal with the spectral problem with Neumann (Robin) boundary conditions. 
By the same method one can obtain invariants of continuous deformations of the potential for the Schrödinger operator in closed manifolds.

\subsection{Further remarks}\label{subsec:FurtherRemarks}
 
{\bf 1. Are the Birkhoff coefficients of the monodromy operator $W(\lambda)$  isospectral invariants?}\quad 
Consider the Birkhoff invariants $p^0_{j,\alpha}$ of the QBNF  of the monodromy operator $W(\lambda)$ given by 
Proposition \ref{prop:commutator}. We are going to show that for almost any torus of a non-degenerate (in the sense of Kolmogorov)  completely integrable system and more generally of a KAM system, $p^0_{j,\alpha}$ are isospectral invariants. To this end we fix 
$s\in\N$ and  suppose that the continuous family $K_t\in C^\ell(\Gamma;\R)$, $t\in [0,1]$, satisfies the isospectral condition
(H$_1$)$_s$ - (H$_2$). 
 Denote by $M$ the smallest integer such that 
$M > \max\{2d,s\}$ and fix $\ell > M(\tau+2) + 2n + (n-1)/2$.  Consider the continuous family of quasimodes of order $M$ given  by Theorem \ref{quasimodes}, which is associated to the Kronecker torus $\Lambda$ and to the continuous family ${\mathcal B}:= \{K_t:\, t\in [0,1]\}\subset C^\ell(\Gamma;\R)$. 
Recall that the quasi-eigenvalues of the quasimode in Theorem \ref{quasimodes} associated with the Kronecker torus $\Lambda$ are given by
$$
\mu_q(t) = \mu_q^0 +c_{q,0} + c_{q,1}(t) (\mu_q^0)^{-1}  + c_{q,2}(t) (\mu_q^0)^{-2}+ \cdots + 
c_{q,M}(t) (\mu_q^0)^{-M}\, , 
$$
where $\mu_q^0$ and $c_{q,0}$ do not depend on $t$. The coefficients  $c_{q,j}(t)$, $1\le j\le M$, 
are obtained by \eqref{quasi-coefficients}.
Taking  $s=0$ we have shown that $c_{q,1}(t)$ does not depend on $t$. Then \eqref{quasi-coefficients} and \eqref{system1} imply that 
 $b_{q,1}(t)$ and $p^0_{0,0}$ do not depend on $t$ either.

Fix $1\le j\le s+1$ and consider the system \eqref{quasi-coefficients}. Recall from \eqref{W} that $W_{q,j}$ are polynomials of $c_{q, r}$ and $b_{q,r'}$,  
$0 \le r,r' \le j-1$,  and $V_{q,j}$  is given by \eqref{V}. 
 Suppose that the functions  $t\to p^0_{l,\alpha}$, $l+|\alpha|<j-2$,  $t\to b_{q,l}(t)$ and $t\to c_{q,l}(t)$, $l<j$,  do not depend on $t$. Using Lemma \ref{lemma:invariants} we would like to show that the functions  $t\to p^0_{l,\alpha}$, $l+|\alpha|=j-1$,  $t\to b_{q,j}(t)$ and $t\to c_{q,j}(t)$  do not depend on $t$ either. 
 At first glance this does not seem  possible  since there are too many unknown functions.  
 To get rid of the terms $p^0_{l,\alpha}$, where $l +|\alpha|=j-1$ and $|\alpha|>0$,  we suppose that either $P=B^m$ is   completely integrable and non-degenerate (in the sense of Kolmogorov)  or that it is a small perturbation of such a system (KAM system). 
 In both cases there is a large family of invariant tori $\Lambda(\omega)$ of $P$ with frequencies $\omega\in \Theta_\kappa^\tau$ satisfying a $(\kappa,\tau)$-Diophantine condition, where $\kappa>0$ and $\tau>n-1$ are fixed and $\kappa$ is sufficiently small. These tori are parameterized by their frequencies  $\omega\in \Theta_\kappa^\tau$,  and their union  is of a positive Lebesgue measure in the phase space. Moreover, the frequency  map $\omega\to I^0(\omega)$ assigning to each  $\omega\in \Theta_\kappa^\tau$ the corresponding 
 action $I^0(\omega)$ on $\Lambda(\omega)$ can be extended to a smooth map $\omega\to I^0(\omega)$ in a neighborhood of  $\Theta_\kappa^\tau$
 which is a diffeomorphism locally (cf. \cite{La, Poe, P2,P3}). In particular, one can parameterize the invariant tori by the corresponding action variables $I^0= I^0(\omega)$. 
 
  Consider the set $\tilde \Theta_\kappa^\tau$ of frequencies  $\omega\in \Theta_\kappa^\tau$ of positive Lebesgue density in $\Theta_\kappa^\tau$. By definition,  a frequency  $\omega\in  \Theta_\kappa^\tau$ belongs to 
 $\tilde \Theta_\kappa^\tau$ if  the intersection of any neighborhood of $\omega$ with  $\Theta_\kappa^\tau$  has a positive Lebesgue measure.  The advantage of working with $\tilde \Theta_\kappa^\tau$ is that if a smooth function is zero on $\tilde \Theta_\kappa^\tau$ then any partial derivative of that function is also zero on $\tilde \Theta_\kappa^\tau$. Moreover, the complement of  $\tilde \Theta_\kappa^\tau$ in  $\Theta_\kappa^\tau$ has Lebesgue measure zero. 
 On the other hand one can prove (using a suitable  KAM theorem and a {\em simultaneous} QBNF for all the invariant tori $\Lambda(\omega)$ with frequencies   $\omega\in \tilde \Theta_\kappa^\tau$) that the Birkhoff invariants $p^0_{j,0}(I^0)$ associated with $\Lambda(\omega)$ depend smoothly on $I^0 = I^0(\omega)$,  $\omega\in \tilde \Theta_\kappa^\tau$,  in a Whitney sense. In particular,  $p^0_{l,\alpha}(I^0)=\partial^\alpha p^0_{l,0}(I^0)$, for $l+|\alpha|\le j-1$ and any $I^0= I^0(\omega)$, where $\omega\in \tilde \Theta_\kappa^\tau$. Hence,  $p^0_{l,\alpha}(I^0)$, $l+|\alpha|=j-1$, $|\alpha|>0$, does not depend on $t$ by the inductive assumption. Now Lemma \ref{lemma:invariants},  \eqref{V} and \eqref{quasi-coefficients}   imply that  that the functions  $t\to p^0_{j-1,0}(I^0)$,   $t\to b_{q,j}(t)$ and $t\to c_{q,j}(t)$,  do not depend on $t$ for any $I^0= I^0(\omega)$, where $\omega\in \tilde \Theta_\kappa^\tau$.  
 
The Birkhoff invariants  $p^0_{j,0}(I^0)$,  
 $j\le s$,  would give  further isospectral invariants for the problem (\ref{thespectrum})   
involving integrals of polynomials of the derivatives of $K_t$. It can be  shown that  $p^0_{2,0}$ is of the form
\[
\sum_{j=0}^{m -1}
\int_{\Lambda_j}\left[
\left(\frac{K\circ \pi_\Gamma}{\sin \theta}\right)^2\, + D (K\circ \pi_\Gamma)\right]d\mu_j \, ,
\]
where $D$ is a differential operator of degree 2. This invariant could be used to remove one of the symmetries for the function $K$ in Corollary \ref{Coro:bouncing-ball} and Corollary \ref{Coro:liouville}. Using further invariants one may remove all the symmetries. 
Further invariants could be  obtained in the case of Theorem \ref{main-potential} as well.
The details will be given elsewhere. 

\vspace{0.3cm}
\noindent
{\bf 2. Deformations of the Riemannian metric. } 
Our method can be   applied as well in the case of  deformations of the Riemannian metric using a variant of Lemma \ref{lemma:invariants} and  
combining it with certain   results of  \cite{P2} and \cite{P3}. 
It could  be proved in this way that  KAM tori with Diophantine vectors of rotation in  $\tilde \Theta_\kappa^\tau$ 
are isospectral invariants (up to a symplectic conjugation) of   deformations of the Riemannian metric on compact manifolds with or without boundary. In particular, the Liouville classes of these tori (cf. \cite{P2}) will be  constant along the  deformation. Using  \cite{PT2} one could  generalize the results of Hezari and   Zelditch \cite{H-Z2} for Liouville billiard tables in dimensions 2 and 3. 

\vspace{0.3cm}
\noindent
{\bf 3. Elliptic tori.} Continuous families of quasimodes can be constructed as well for {\em elliptic} low-dimensional  tori under suitable Diophantine conditions, which will give further isospectral invariants. An elliptic closed geodesic can be considered  as an elliptic torus of dimension one. Elliptic tori of dimension $\ge 2$ usually appear in Cantor families and can be constructed via  the KAM theory \cite{Poe1}.

\vspace{0.3cm}
\noindent
{\bf 4. Resonances for exterior problems.} We  point out that the method we use can be applied 
whenever there exists a continuous
family of quasimodes of the spectral problem.  It can be used 
also for the Laplacian $\Delta_K$ in the exterior  $X=\R^n \setminus \Omega$ 
of a bounded domain in $\R^n$ with a $C^\infty$-smooth boundary with Robin boundary conditions on it. 
In this case an analogue of (H$_1$)-(H$_2$) can be formulated for the resonances  of $\Delta_K$ 
close to the real axis replacing the intervals in the definition of $\mathcal I$ by boxes in the complex upper
half plain of the form $[a_k,b_k] + i[0,d_k]$, $d_k>0$, $\lim d_k=0$. Given a Kronecker torus $\Lambda$ of $B$ we obtain  continuous in $K$ quasimodes of $\Delta_K$ associated to $\Lambda$. 
By a result of Tang and Zworski \cite{TZ} and Stefanov \cite{Ste} 
the corresponding  quasi-eigenvalues are close to resonances and one obtains an analogue of Lemma \ref{lemma:invariants} and of Theorem 1.1.

\appendix
\section{Appendix.}
\setcounter{equation}{0}

\renewcommand{\theequation}{A.\arabic{equation}}

\noindent
\subsection{ Parametrix of the Dirichlet problem for the Helmholtz equation at high frequencies.}\label{subsec:parametrix}
We are going to construct a  microlocal outgoing parametrix  
$H(\lambda):L^2(\Gamma)\rightarrow  C^\infty( \widetilde X )$ 
of the Dirichlet problem for the Helmholtz equation.
% with ``initial data'' on the boundary $\Gamma$ which are microlocally concentrated 
%in  an open subset $U$ of the domain of definition  $\widetilde B^\ast \Gamma$ of the billiard ball map. 
We consider $H(\lambda)$ as a Fourier integral operator with a large parameter $\lambda$ or equivalently as a semiclassical Fourier integral operator with a small parameter $ \hbar=1/\lambda$. 
It will be obtained  by means of  Maslov's canonical operator as it was presented by Duistermaat \cite{Du} (see also \cite{Alex}, \cite{C-P}, \cite{CV3}, \cite{E-Z}, \cite{Mein}, \cite{PP}). 
The operator $H(\lambda)$ will satisfy   the Helmholtz equation at high frequencies ($|\lambda|\to\infty$), i.e. 
\begin{equation}
\forall\,  N\in \N\, , \quad (\Delta-\lambda^2)H(\lambda)u = O_N(|\lambda|^{-N})u
\label{parametrix1}
\end{equation} 
in a neighborhood of $X$ in the smooth extension 
$\widetilde X$ of the Riemannian manifold $X$ with suitable ``initial data''  on $\Gamma$.

The corresponding canonical relation is an embedded Lagrangian  submanifold of  $T^\ast \widetilde X \times T^\ast\Gamma$  given by 
\begin{equation}\label{canonical-relation}
{\mathcal C} \ :=\ \left\{( 
\varrho', \varrho)\in T^\ast \widetilde X\times T^\ast  \Gamma\,
:\ \varrho'=  \exp(s X_ {\widetilde h})(\pi_\Sigma ^+(\varrho))\, ,\  \varrho\in U\, ,\ -\varepsilon <s < T(\varrho) +\varepsilon \right\}\, ,\ \varepsilon >0\, ,
\end{equation}
where $U$ is an open  subset of the domain of definition  $\widetilde B^\ast \Gamma$ of the billiard ball map. 
Recall that $\widetilde h$ is the Hamiltonian corresponding to the Riemannian metric $\widetilde g$ via the Legendre transform, $X_ {\widetilde h}$ is the corresponding Hamiltonian vector field, and the map $\pi_\Sigma ^+:B^\ast \Gamma\to \Sigma^+$ is defined by \eqref{outgoing-vector}. In particular,  $\exp(s X_ {\widetilde h})(\pi_\Sigma ^+(\varrho))$ lies on the cosphere bundle  
\begin{equation}
\widetilde \Sigma:= S^\ast \widetilde X= \{(x,\xi)\in T^\ast \widetilde X:\, \widetilde h(x,\xi)=1\} .
\label{cosphere}
\end{equation}
Moreover,  $T:U\to (0,+\infty)$ is the ``return time function'' which assigns  to each $\varrho\in U$ the time of the first impact at the boundary, i.e.  the 
smallest positive time $t=T(\varrho)$ such that 
\[
\exp(tX_{\widetilde h})(\pi^+_\Sigma(\varrho))\in \Sigma^-\, .
\]
We suppose that $U$ is a connected open subset of $\widetilde B^\ast \Gamma$ and that its closure $\overline U$  is still in  $\widetilde B^\ast \Gamma$. Then $T$ is a smooth function and its image is a finite  interval $(a,b)$, $0<a<b$. 
Denote by ${\mathcal C}'$ the corresponding Lagrangian submanifold of $T^\ast (X\times \Gamma)$ given by
\begin{equation}
(x,y,\xi,\eta)\in {\mathcal C}' \ \Leftrightarrow \ (x,\xi,y,-\eta)\in {\mathcal C}.
\label{lagrangian-manifold}
\end{equation}
${\mathcal C}$ can be  parameterized by the variables $(s,\varrho)$ in \eqref{canonical-relation} which gives also a parametrization of ${\mathcal C}'$. 
Notice that  ${\mathcal C}'$ is always an {\em exact Lagrangian submanifold} of  $T^\ast (X\times\Gamma)$ although $U$ may not be simply connected. This is important for the existence of globally defined  oscillatory integrals associated to ${\mathcal C}'$. 
\begin{Lemma}\label{Lemma:exact}
The pull-back $\alpha' = \alpha|_{{\mathcal C}'}$ of the canonical one-form $\alpha := \xi dx + \eta dy$ of $T^\ast(\widetilde X\times\Gamma)$ to  ${\mathcal C}'$ via the inclusion map is exact.
\end{Lemma}
{\em Proof}. 
To prove it let us identify the open co-ball bundle $B^\ast \Gamma$ with the symplectic manifold $\Sigma^+$ via the diffeomorphism $\pi^+_\Sigma$ and set $\tilde U := \pi^+_\Sigma(U)$. The symplectic two-form of $\Sigma^+$ is $- d\alpha_\Sigma$, where $\alpha_\Sigma$ is the pull-back of the fundamental one-form $\eta d y$ of $B^\ast \Gamma$ under $\pi_\Sigma$.
We have also $\alpha_\Sigma = \imath_\Sigma^\ast(\xi dx)$, where  $\imath_\Sigma:\Sigma^+\to T^\ast \widetilde X$ 
is the inclusion map. 
 Consider the  submanifold
\[
\widetilde{\mathcal C}:=\left\{((x,\xi),(y,\eta))\in T^\ast \widetilde X\times \Sigma^+ :\, 
(x,\xi)=\exp(s X_ {\widetilde h})(y,\eta),\  (y,\eta)\in \widetilde U\, ,\ -\varepsilon <s < T(\pi_\Sigma(y,\eta)) +\varepsilon \right\} 
\]
of the symplectic manifold $T^\ast \widetilde X\times \Sigma^+$ equipped with the symplectic two-form $- d \tilde\alpha$, where 
$\tilde\alpha:= \xi dx - \alpha_\Sigma$. By construction, $\widetilde{\mathcal C}$ is a Lagrangian submanifold of $T^\ast \widetilde X\times \Sigma^+$. On the other hand, any closed parameterized $C^1$  curve on $\widetilde{\mathcal C}$ can be written in the form 
$\gamma(t) = (\gamma_1(t);\gamma_0(t))$, $t\in [0,t_0]$, where 
$\gamma_1(t):=\exp(s(t) X_ {\widetilde h})(\gamma_0(t))$,   $s: [0,t_0]\to\R$ is  $C^1$, and $\gamma_0:[0,t_0]\to \Sigma^+$ is a closed $C^1$ curve , i.e. $\gamma_0(t_0)=\gamma_0(0)$. The cylinder 
\[
C:=\{\exp(s X_ {\widetilde h})(\gamma_0(t)):\, 0\le \pm s \le \pm s(t),\, 0\le t\le t_0\}
\]
 is an isotrope submanifold of $T^\ast \tilde X$ and by Stock's theorem  the integral of the fundamental one-form $\xi dx$ on its boundary is zero, which implies that the integral of 
 $\tilde \alpha$ on $\gamma$ is zero. Consequently the pull-back of the fundamental one-form  to ${\mathcal C}$ via the inclusion map is an exact one-form which means that ${\mathcal C}$ is an exact Lagrangian submanifold of $T^\ast \widetilde X\times T^\ast\Gamma$.
\finishproof

Hence, ${\mathcal C}'$ 
is an exact Lagrangian submanifold of $T^\ast (\widetilde X\times \Gamma)$ which allows us to   associate to  it a class of  oscillatory integrals defined by non-degenerate phase functions as in \cite{Du}. For sake of completeness we recall some basic facts about oscillatory integrals which we need below. 

Let   $M^d$  be a smooth paracompact  manifold of  dimension $d$ and  $\Lambda$ an embedded Lagrangian submanofold  of $T^\ast M^d$. 
Recall from   \cite{Du} and  \cite{Hor}  that a real valued phase function  $\Phi(x,\theta)$ defined in a neighborhood  of a point $(x^0,\theta^0)\in \R^d\times \R^N$ with  $d_\theta\Phi(x^0,\theta^0) = 0$ is  non-degenerate at $(x^0,\theta^0)$ if 
$$
{\rm rank}\, d_{(x,\theta)}d_\theta \Phi (x^0,\theta^0)= N,
$$
hence,  there is a  neighborhood $V$  of $(x^0,\theta^0)$ such that
$C_\Phi := \{(x,\theta)\in V: \,  d_\theta \Phi = 0\}$ is a smooth manifold of dimension $d$.  Moreover, the  map
\begin{equation}
\imath_{\Phi} : C_\Phi \ni (x,\theta) \longrightarrow (x,d_x\Phi(x,\theta))
\in \Lambda_\Phi:=\imath_{\Phi} (C_\Phi)
                                             \label{lagrangian1}
\end{equation}
is of rang $d$, and the non-degenerate function $\Phi(x,\theta)$, $(x,\theta)\in V$ is said to  define locally $\Lambda$ near $\nu^0\in \Lambda$  if the map \eqref{lagrangian1} 
is a  diffeomorphism onto an open neighborhood of $\nu^0$  of $\Lambda$. The collection $(\Lambda_\Phi, \imath_{\Phi}^{-1})$ provides $\Lambda$ with an atlas of local carts. 

Given $m\in\R$,  one can define a class of oscillatory integrals  associated to $\Lambda$ of order $m$ as follows.  An oscillatory integral of that class is a collection of smooth functions  $u(\cdot,\lambda)$, $\lambda\in {\cal D}$,  in $M^d$, where ${\cal D}$ is an unbounded subset of $\C$ lying in a band $|{\rm Im}\, \lambda|<C$,  such that  for any $\nu^0 = (x^0,\xi^0)\in \Lambda$ there
exist local coordinates $x$ in a neighborhood of $x^0$,  a non-degenerate real valued phase function $\Phi(x,\theta)$, and  a classical amplitude
$b(x,\theta,\lambda)$ of order $0$ such that  $\Phi(x,\theta)$ defines $\Lambda$ locally near $\nu^0$, and  the oscillatory integral
\begin{equation}
u^0(x,\lambda)\,  =\,  \left(\frac{\lambda}{2\pi}\right)^{m+\frac{d+2N}{4}}\, 
\int_{{\R}^N}\ e^{i\lambda\Phi(x,\theta)} b(x,\theta,\lambda)\, d\theta
                                          \label{oscillatory1}
\end{equation}
represents microlocally $u(x,\lambda)$ in a neighborhood of $\nu^0$ (cf. \cite{Alex}, \cite{Du}, \cite{Mein}, \cite{PP}, and also \cite{Hor}, Section 25.1). The latter  means that there is a neighborhood $U\subset T^\ast M^d$ of $\nu^0$ such that 
$$
{\rm WF} (u-u^0)\cap U=\emptyset,
$$
where ${\rm WF} $ stands for the frequency set  or semiclassical $\hbar$-wave-front with $\hbar=1/\lambda$ (cf. \cite{Alex}, \cite{D-S}, \cite{E-Z}). If $N=0$, then $u^0(x,\lambda)= \left(\frac{\lambda}{2\pi}\right)^{m+\frac{d}{4}}\, 
 e^{i\lambda\Phi(x)} b(x,\lambda)$. 
 
 By a {\em classical  symbol} of order $0$ we mean a collection of complex valued  functions 
$b(\cdot,\lambda)$ in $C_0^\infty(\R^d\times \R^{N})$,  $\lambda\in {\cal D}$,  such that the family $b(\cdot,\lambda)$, $\lambda\in {\mathcal D}$,  is uniformly compactly supported, and the symbol $b$ admits an asymptotic expansion 
\[
b(x,\xi,\lambda) \sim b_0(x,\xi) + b_1(x,\xi)\lambda^{-1} +\cdots \, , \quad \mbox{as}\ |\lambda| \to \infty\, ,
\]
which means that for any $k\ge 0$ and any multi-indices $\alpha,\beta\in \N^{n-1}$ there is a constant $C_{k,\alpha,\beta}>0$ such that
\begin{equation}
\big|\partial_x^\alpha \partial_\xi^\beta \big(b(x,\xi,\lambda) -\sum_{j=0}^k b_j(x,\xi)\lambda^{-j}\big)\big| \le C_{k,\alpha,\beta} (1+|\lambda|)^{-k-1}
\label{symbol}
\end{equation}
in $\R^d\times \R^{N}\times {\mathcal D}$. 
In particular, the  functions $b_j$, $j\ge 0$, belong to $C_0^\infty(\R^d\times \R^{N})$ and there is a compact  $V$ of $\R^d\times \R^{N}$ such that the support ${\rm supp}\, b_j$ of $b_j$ is a subset of $V$ for any $j\ge 0$. 
From now on we  consider $u(x,\lambda)$ as a half-density on $M^d$ and denote the class of these oscillatory integrals by $I^{m}(M^d, \Lambda; \Omega^{1/2}_{M^d})$.

The canonical relation  ${\mathcal C}$ defined by \eqref{canonical-relation} is quite special, since for any $s$ fixed it is a graph of a canonical transformation, and one can define  locally the Lagrangian manifold  ${\mathcal C}'$ by a non-degenerate
phase function as follows. 

Fix $\nu^0 = (x^0,y^0,\xi^0,\eta^0)\in {\mathcal C}' $ and choose a smooth submanifold  $M$ of $\tilde X$ of dimension $n-1$ passing through $x^0$ and  transversal to the geodesic starting from $y^0$ with codirection $(\eta^0)^+$. Let $x=(x',x_n)$ be normal coordinates to $M$. Then  $M$ is given locally by the equation $\{x_n=0\}$ and the Hamiltonian $\widetilde h$ is of the form  $\widetilde h(x,\xi) = \xi_n^2 + \widetilde h_0(x,\xi')$.   If $x^0\in\Gamma$, we take $M$ to be a neighborhood of $x^0$ in $\Gamma$, then $\widetilde h_0(x,\xi')= h_0(x,\xi')$. 
\begin{Lemma}\label{lamma:phase-function}
There exist  local coordinates  $y\in \R^{n-1}$ in a neighborhood of $y^0$ in $\Gamma$, and 
a non degenerate  phase function 
\begin{equation}
\Phi (x,y, \theta)= \phi (x,\theta) - \langle y,\theta\rangle 
                                             \label{phase-function}
\end{equation}
defined in a neighborhood $W$ of $(x^0,y^0, \theta^0)$, $\theta^0=-\eta^0$,  
in $\R^n\times\R^{n-1}\times\R^{n-1}$ such that locally ${\mathcal C}' _\Phi= {\mathcal C}'$. 
\end{Lemma}
The proof of the Lemma is the same as in   H\"ormander \cite{Hor}, Proposition 25.3.3. 
The non-degeneracy of the phase function means that
\begin{equation}
\det \frac{\partial^2\phi}{\partial x'\partial\theta} (x,\theta) \neq 0 
                                             \label{nondegenerate-phase-function}
\end{equation}
in a neighborhood of $(x^0,\theta^0)$. 
We take  $\phi(x',x_n,\theta) = \langle x',\theta\rangle + O(x_n)$ if the image of  $\nu^0$ by the involution in \eqref{lagrangian-manifold}   belongs to $\Sigma^+\times U$, which means that $\nu^0 = (y^0,y^0,-(\eta^0)^+,\eta^0)$. 

From now on we take  the set  ${\mathcal D}$ to be  as in  \eqref{strip}, and we consider  the corresponding oscillatory integrals as 1/2-densities by multiplying \eqref{oscillatory1} by the half-density $|dx|^{1/2}$ \cite{E-Z}.  In the same way, we consider the Helmholtz operator as an operator acting on $1/2$-densities $u\in C^\infty(\widetilde X, \Omega^{1/2}) $.

We are looking for a parametrix $H(\lambda)$ of  \eqref{parametrix1}  in the class $I^{-1/4}(\tilde X, \Gamma, {\mathcal C}; \Omega^{1/2})$ which means that its  Schwartz kernel $K_H(x,y,\lambda)$  is in  
$I^{-1/4}(\tilde X\times \Gamma, {\mathcal C}'; \Omega^{1/2}_{\tilde X\times\Gamma})$. 
The order $-1/4$ comes from the 
``initial data'' at the boundary which will be given  by a classical $\lambda$-PDO $\Psi(\lambda)$ of order zero (see Lemma \ref{Lemma:restriction} below). 
Moreover, as we shall see below,  operators of that class are uniformly bounded  as operators from $L^2(\Gamma)$ to $L^2(\widetilde X)$. 
In any 
local chart  of $\Gamma$ the Schwartz kernel of a $\lambda$-PDO $\Psi(\lambda)$ is a $1/2$-density $K_\Psi|dx|^{1/2}|dy|^{1/2}$, where $K_\Psi$ is  an oscillatory integral of the form 
\begin{equation}
K_\Psi(x,y,\lambda):=\left(\frac{\lambda}{2\pi}\right)^{n-1} \int_{\R^{n-1}} \, e^{i\lambda \langle x-y,\xi\rangle} q(x,\xi,\lambda)\, d\xi \, ,\quad \lambda\in {\mathcal D}\, ,
\label{classical-pdo}
\end{equation} 
and $q\sim q_0 + q_1 \lambda^{-1} + \cdots $ is a  classical amplitude of order zero.

To any non-degenerate  phase function $\Phi(x,\theta)$ of the form \eqref{phase-function} generating ${\mathcal C}'$ near a point $\nu^0$ there is a  classical amplitude $b(x,\theta,\lambda)\sim b_0(x,\theta) + b_1(x,\theta)\lambda^{-1}+ \cdots$  such that  
the Schwartz kernel of $H(\lambda)$ can be written microlocally near  $\nu^0$ as a $1/2$-density $I_\Phi|dx|^{1/2}|dy|^{1/2}$, where 
\begin{equation}
I_\Phi(x,y,\lambda) = \left(\frac{\lambda}{2\pi} \right)^{n-1}
\int_{{\R}^{n-1}}\ e^{i\lambda \Phi(x,y,\theta)} b(x,\theta,\lambda)  d\theta 
                                          \label{oscillatory}
\end{equation}
(see  \eqref{oscillatory1}). 
 In particular, $H(\lambda)u$ is a $C^\infty$ smooth function for any fixed $\lambda$ and $u\in L^2(\Gamma)$ (the essential support of  $b$ with respect to $\theta$ is a compact set)
 but it highly oscillates at  its frequency set.  We are going to describe the principal symbol of $H(\lambda)$. 

 By Lemma \ref{Lemma:exact}  the Lagrangian manifold  ${\mathcal C}'$ is exact.  
Then the principal symbol
of $H(\lambda)$ is of the form $(\lambda/2\pi)^{-1/4}\sigma(H)$, where  
\begin{equation}
\sigma(H) =  e^{i\lambda f}\,  \sigma_1 \otimes \sigma_2,
\label{principal-symbol1}
\end{equation}
$\sigma_1$ and $\sigma_2$ are smooth sections of the half-density bundle
$\Omega_{1/2}({\mathcal C}')$ and the Maslov bundle $M({\mathcal C}')$, respectively, 
and $f$  is a smooth function representing the Liouville factor 
(cf. \cite{Du}). 
 The function $f$ is defined locally by means of the phase functions $\Phi$ in \eqref{oscillatory} as follows
\begin{equation}
f(\nu) = \widetilde \Phi(\nu):=  \Phi(\imath^{-1}_\Phi(\nu)),\
\nu\in {\mathcal C}'_\Phi \subset {\mathcal C}' ,
\label{liouville-factor}
\end{equation}
where $\imath_\Phi$ is given by \eqref{lagrangian1}.
Such a function $f(\nu)$ exists globally on ${\mathcal C}'$ since the latter is exact.
Indeed, we have locally
$$
d \widetilde \Phi = \left(\imath^{-1}_\Phi\right)^\ast \left(d\left(\Phi_{|C_\Phi}\right)\right)
=\left(\imath^{-1}_\Phi\right)^\ast \left(\frac{\partial\Phi}{\partial x}dx + \frac{\partial\Phi}{\partial y}dy +
 \frac{\partial\Phi}{\partial \theta}d\theta\right)
_{|C_\Phi}$$
$$
= \left(\imath^{-1}_\Phi\right)^\ast \left(\frac{\partial\Phi}{\partial x}dx + \frac{\partial\Phi}{\partial y}dy\right)
_{|C_\Phi} =  \alpha'.
$$
If $I_{\Phi_1}$ is  another microlocal  representation of $H$ near $\nu^0$ given by (\ref{oscillatory}) with
a phase function $\Phi_1$, then 
$$
d\widetilde\Phi = \alpha' = d\widetilde\Phi_1 \ \mbox{and}\ \widetilde\Phi(\nu^0) =
\widetilde\Phi_1(\nu^0)
$$
which implies $\widetilde\Phi = \widetilde\Phi_1$ in a neighborhood  of $\nu^0$. Since
${\mathcal C}'$ is exact, there exists a globally defined exponent $f$ of
the Liouville factor. Note that the function $f$ is uniquely determined
on ${\mathcal C}'$ modulo a constant, hence, it is enough to know $f$ at a
single point in order to recover it on the whole manifold ${\mathcal C}'$. We have $f=0$ for any $\nu^0 = (y^0,y^0,-(\eta^0)^+,\eta^0)$, where
 $(x^0,-\eta^0)\in U$. Parameterizing ${\mathcal C}'$ by the variables $(s,\varrho)$ we obtain that $f(s,\varrho) = 2s$ is just the action along the corresponding bicharacteristic arc.

Consider the flow $S^t:{\mathcal C}\to {\mathcal C}$ defined by 
\begin{equation}
S^t(x,\xi, y,\eta)=S^t(\exp(sX_{\widetilde h})(y,\eta),y,\eta)= (\exp((s+t)X_{\widetilde h})(y,\eta),y,\eta), 
\label{flow}
\end{equation}
and denote by $Y=(X_{\widetilde h},0)$ the corresponding vector field. We use the same notations for the corresponding flow and vector field on 
${\mathcal C}'$. We are going to choose a half-density invariant with respect to $S^t$. 

Consider the volume form $\beta_0$ on ${\mathcal C}'$  given by 
$ds\wedge (dy_1 \wedge d\eta_1)\wedge \cdots  \wedge (dy_{n-1} \wedge d\eta_{n-1})$ 
in the coordinates $(s,y,\eta)$.
It is obviously invariant with respect to the flow $S^t$ or equivalently, the Lie derivative ${\mathcal L}_Y \beta_0$ vanishes. 
 Then the $1/2$-density $\sigma_0:= |\beta|^{1/2}$ on ${\mathcal C}'$ is invariant with respect to $S^t$, and we set 
\begin{equation}
\sigma_1= \widetilde b_0 \sigma_0, \quad \widetilde b_0\in C_0^\infty({\mathcal C}'). 
\label{principal-symbol2}
\end{equation}
There is another way to describe the half-density $\sigma_0$. The cosphere bundle 
$\widetilde \Sigma$ in \eqref{cosphere} 
is a contact manifold with a contact one-form $\alpha_0 = \imath^\ast (\xi dx)$ given by  the pull-back of the fundamental one-form $\xi dx$ of $T^\ast \widetilde X$ to $\widetilde \Sigma$   via the embedding  $\imath:\widetilde\Sigma \to T^\ast \widetilde X $. Consider the map 
\[
\psi: (s,y,\eta) \mapsto (x(s,y,\eta),\xi(s,y,\eta)):=\exp(sX_{\widetilde h})(\pi_\Sigma^+(y,\eta))\in \widetilde \Sigma .
\]
Then we have 
\[
\begin{array}{lcrr}
\displaystyle \psi^\ast(\alpha_0) = \xi(s,y,\eta)\frac{\partial x}{\partial s}(s,y,\eta) ds + (\pi_\Sigma^+)^\ast \left( \exp(sX_{\widetilde h}) \right)^\ast \alpha_0 \\[0.3cm]
\displaystyle = \xi \frac{\partial \widetilde h}{\partial \xi} ds + (\pi_\Sigma^+)^\ast \alpha_0  + \alpha_1= 
2ds + \eta dy +  \alpha_1, 
\end{array}
\]
where $\alpha_1 = a\, dy + b\,  d\eta$ is a closed form, 
which implies that
\begin{equation}
\sigma_0 = \psi^\ast\left(\left(\frac{|\alpha_0\wedge (d\alpha_0) ^{n-1}|}{2 (n-1)! }\right)^{1/2}\right). 
\label{density}
\end{equation}
Consider as above  the 
 local coordinates $(x',x_n)$ near $x^0$ which are normal to the submanifold $M=\{~x_n~=~0~\}\subset \widetilde X$ and recall that  $\tilde h(x,\xi) = \xi_n^2 + \widetilde h_0(x,\xi')$. Setting $\xi_n= \pm \sqrt{1- \widetilde h_0(x,\xi')}$ and parameterizing $\widetilde\Sigma$ near $(x^0,\xi^0)$  by $(x,\xi')$ 
we easily obtain 
$$
\begin{array}{lcrr}
\beta_0:=\displaystyle \frac{\alpha_0\wedge (d\alpha_0) ^{n-1}}{(n-1)!} = \left(\xi_n - \sum_{j=1}^{n-1} \xi_j\frac{\partial \xi_n}{\partial \xi_j} \right)dx_n\wedge (dx_1 \wedge d\xi_1)\wedge \cdots  \wedge (dx_{n-1} \wedge d\xi_{n-1})\\[0.5cm] 
\displaystyle = \frac {dx_n\wedge (dx_1 \wedge d\xi_1)\wedge \cdots  \wedge (dx_{n-1} \wedge d\xi_{n-1})}{\xi_n}
\end{array}
$$
in a neighborhood of $(x^0,\xi^0)\in S^\ast \widetilde X$. Then \eqref{principal-symbol2} and \eqref{density} imply
\begin{equation}
\label{leray-density}
\sigma_0= \frac{|dx\wedge d\xi'|^{1/2}}{\sqrt{2|\xi_n|}}\quad \mbox{and}\quad   
\sigma_1= \frac{{\widetilde b_0}}{\sqrt{2|\xi_n|}}\,  |dx\wedge d\xi'|^{1/2} 
\end{equation}
in the coordinates $(x,\xi')$, where $ |dx\wedge d\xi'| :=|dx_1\wedge\cdots\wedge dx_n\wedge d\xi_1\wedge\cdots\wedge d\xi_{n-1}|$. 
Notice also that $ \beta_0$ is a Leray form of $\widetilde\Sigma$, i. e. 
$\pm \beta_0= \imath^\ast \beta_1$, where $\imath:\widetilde\Sigma \to T^\ast \widetilde X $ is the corresponding embedding  and $\beta_1$ is  a $2n-1$ - form on $T^\ast \widetilde X $ such that 
$\beta_1\wedge d\widetilde h = dx\wedge d\xi.$ 

We are going to relate the leading term $b_0$ of the amplitude of the oscillatory integral \eqref{oscillatory} to the function $\widetilde b_0$ in \eqref{principal-symbol2}. 
There are two natural choices of $\sigma_1$ depending on the parametrization of ${\mathcal C}'$. 
Recall that  
$${\mathcal C}'=\{(x,d_\theta\phi(x,\theta),d_{x}\phi(x,\theta),-\theta)\}$$ 
in a neighborhood of $\nu^0$, where $\phi$ is given by \eqref{phase-function}. 
Firstly, parameterizing ${\mathcal C}'$ by the variables $(x,\eta)$, $\eta=-\theta$,  and using the definition of the principal symbol as in   H\"ormander \cite{Hor}, Sect. 25.3,  we get locally near $\nu^0$ that 
$$
\sigma_1  = b_0(x,-\eta)  |dx\wedge d\eta'|^{1/2}. 
$$ 
On the other hand, parameterizing  ${\mathcal C}'$  by the variables $(x,\xi')$ we get
$$
\left|dx\wedge d\eta'\right|=  \left|\det \phi''_{x'\theta}(x,-\eta(x,\xi'))\right|^{-1}\left|dx\wedge d\xi'\right|,
$$ 
where $\eta=\eta(x,\xi')$ is obtained from $\xi'=\phi'_{x'}(x,-\eta)$ by the implicit function theorem.  Hence, 
$$
\sigma_1  = b_0(x,-\eta(x,\xi'))  \left|\det \phi''_{x'\theta}(x,-\eta(x,\xi'))\right|^{-1/2} \left|dx\wedge  d\xi'\right|^{1/2}. 
$$ 
and taking into account \eqref{leray-density} we obtain
\begin{equation}
b_0(x,-\eta(x,\xi')) = \frac{\widetilde b_0'(x,\xi')}{\sqrt{2|\xi_n|}}\left|\det \phi''_{x'\theta}(x,-\eta(x,\xi')) \right|^{1/2} 
\label{principal-symbol3}
\end{equation} 
in a neighborhood of $(x^0,\xi^0)$, where 
\[
\xi_n = \pm \sqrt{1-\widetilde h_0(x,\xi')}\, ,\quad \widetilde b_0'(x,\xi') = \widetilde b_0(\pi_1^{-1}(x,\xi',\xi_n)),
\]
and 
$\pi_1: {\mathcal C}' \to \widetilde \Sigma \subset T^\ast \widetilde X$ stands for the local projection  $\pi_1(x,y,\xi,\eta)=(x,\xi)$.  

The Keller-Maslov bundle $M({\mathcal C}')$ of ${\mathcal C}'$ is trivial. We recall from  H\"ormander \cite{Hor1} that a section of the line bundle 
$M({\mathcal C}')$ is given by a family of functions $f_\Phi:{\mathcal C}'_\Phi\to \C$, where $\Phi$ is a non-degenerate generating function of ${\mathcal C}'$ and ${\mathcal C}'_\Phi$ is defined by \eqref{lagrangian1} such that
$f_{\widetilde \Phi}= \exp\left\{ i\frac{\pi}{4}\left({\rm sgn}\,  \Phi''_{\theta\theta}- {\rm sgn}\,  \widetilde\Phi''_{\widetilde \theta\widetilde \theta}\right)\right\}f_{\Phi}$ on ${\mathcal C}'_\Phi\cap {\mathcal C}'_{\widetilde \Phi}$. Moreover, ${\rm sgn}\,  \Phi''_{\theta\theta}- {\rm sgn}\,  \widetilde\Phi''_{\widetilde \theta\widetilde \theta}\in 2\Z$ and it is constant on ${\mathcal C}'_\Phi\cap {\mathcal C}'_{\widetilde \Phi}$. In our case, we can trivialize $M({\mathcal C}')$ in a band $0<s<\delta$ using the phase functions $\Phi$ given by Lemma \ref{lamma:phase-function}. Indeed, ${\mathcal C}'$ is generated in a neighborhood of any point $\nu^0 = (y^0,y^0,-(\eta^0)^+,\eta^0)$ by a phase function $\Phi(x,y,\theta) = \phi(x,\theta)-\langle y,\theta\rangle$, which is a local solution of the Hamilton-Jacobi equation  $\phi'_{x_n}(x,\theta) = \sqrt{1-h_0(x,\phi'_{x'}(x,\theta))}$ with Cauchy data $\phi(x',0, \theta) = \langle x',\theta\rangle$. In particular, 
\[
\Phi(x,y,\theta) = \langle x'-y,\theta\rangle + x_n q(x',\theta) + O(x_n^2),  
\]
where $q(x',\theta):=\sqrt{1-h_0(x',0,\theta)}$. 
Moreover, $q''_{\theta\theta}(x',\theta)$ is a negative definite matrix since $\theta\to h_0(x',0,\theta)$ is a positive definite quadratic form, and we obtain that ${\rm sgn}\,  \Phi_{\theta\theta}(x,y,\theta) = -n+1$ for any $0<x_n<\delta$ and $(y,\theta)$ in a neighborhood of $(y^0,\theta^0)$. This yields a trivialization of the Keller-Maslov bundle in a band ${\mathcal C}'\cap \{0<s<\delta\}$ for some $\delta>0$. Moreover, it  gives a global section $\sigma_2$ of $M({\mathcal C}')$, which is locally constant. In particular, the Lie derivative ${\mathcal L}_Y \sigma_2$ vanishes.

The oscillatory integral $(\Delta_x-\lambda^2)K_H(x,y,\lambda)$ belongs to $I^{3/4}(\tilde X\times \Gamma, {\mathcal C}'; \Omega^{1/2}_{\tilde X\times\Gamma})$, and its principal symbol is just the Lie derivative ${\mathcal L}_{Y}\sigma(H)$ of \eqref{principal-symbol1} multiplied by $(\lambda/2\pi)^{3/4}$ (cf. \cite{Du}, Sect. 1.3 as well as \cite{Hor}, Theorem 25.2.4) since the subprincipal symbol of the Laplace-Beltrami operator is $0$. 
Moreover,  the Lie derivative with respect to $Y$ of the sections $\sigma_0= |ds\wedge dy\wedge d\eta|^{1/2}$ and $\sigma_2$ and of the function $f$ vanishes, hence,  the transport  equation  ${\mathcal L}_{Y}\sigma(H)=0$ becomes 
\begin{equation}
(S^t)^\ast \tilde b_0 = \tilde b_0, 
\label{transport-equation}
\end{equation}
where $\tilde b_0$ is given by \eqref{principal-symbol2}. Multiplying  $\tilde b_0$ with a suitable cut-off function, which  equals $1$ in a neighborhood of 
${\mathcal C}'\cap T^\ast(X\times \Gamma)$,  we can suppose that  $\tilde b_0$ has a compact support with respect to $(s,y,\eta)$.  In this way we obtain 
an operator $H_0(\lambda)\in I^{-1/4}(\widetilde X, \Gamma, {\mathcal C}; \Omega^{1/2})$ such that 
$(\Delta_x-\lambda^2)H_0(\lambda)\in I^{-1/4}( X, \Gamma, {\mathcal C}; \Omega^{1/2})$. Repeating this procedure we get an operator $H_1(\lambda)$ such that $H_0(\lambda)+H_1(\lambda)$ solves (\ref{parametrix1}) modulo  a $\lambda$-FIO of order $-5/4$ and so on. 
The initial data $\tilde b_0|_{s=0}$  will be determined by Lemma \ref{Lemma:restriction}, where taking  
$G_0(\lambda):=\Psi(\lambda)$  as initial data at $\Gamma$ for $s=0$   we determine  $\widetilde b_0|_{s=0}$, and so on.

Denote by   
$\imath_\Gamma^\ast: C^\infty(\widetilde X) \to  C^\infty(\Gamma)$ the operator of restriction  $\imath_\Gamma^\ast(u)=u_{|\Gamma}$.  We would like to represent $\imath_\Gamma^\ast$ microlocally as a $\lambda$-FIO. To this end, denote by ${\mathcal N}$ the conormal bundle of the graph of the inclusion map $\imath_\Gamma: \Gamma \to \widetilde X$ and by ${\mathcal R}= {\mathcal N}^{-1}$ the corresponding inverse canonical relation. In other words,
\[
{\mathcal R}:= \{(x,\xi;x,\widetilde \xi)\in T^\ast \Gamma\times T^\ast \widetilde X:\, x\in\Gamma, \xi=\widetilde \xi|_{T_x\Gamma}\}. 
\]
The operator $\imath_\Gamma^\ast $ can be considered microlocally  as a  $\lambda$-FIO of the class
$I^{1/4}(\Gamma, \tilde X,  {\mathcal R}; \Omega^{1/2})$ (the composition $\imath_\Gamma^\ast \circ A(\lambda)$ belongs to that class for any classical $\lambda$-PDO $A(\lambda)$ of order $0$).  Moreover, its   principal symbol can be identified with $(\lambda/2\pi)^{1/4}$ modulo the corresponding $1/2$-density. Indeed, 
fix  $(x^0,\xi^0)\in T^\ast \widetilde X$ such that $x^0\in \Gamma$ and introduce  as above  normal coordinates $x=(  x',  x_n)$ to $\Gamma$ in a neighborhood of $x^0$. Then $\Gamma = \{  x_n= 0\}$ and $\widetilde h(  x,   \xi) =  \xi_n^2 + h_0(  x,   \xi')$ in the local  coordinates $(  x,   \xi)$ in $T^\ast \widetilde X$.  The Schwartz kernel of  $\imath_\Gamma^\ast$ is represented microlocally in a neighborhood of $(x^0,\xi^0|_{T_x\Gamma},x^0,-\xi^0)\in T^\ast(\Gamma\times\widetilde X)$ by the half-density 
$ K |dx'|^{1/2}|d z|^{1/2}$, where
\begin{equation}
K(x',z,\lambda):=\left(\frac{\lambda}{2\pi}\right)^{n} \int_{\R^{n}} \, e^{i\lambda (\langle x'-z',\xi'\rangle -z_n\xi_n)} \kappa(x,\xi)\, d\xi \, ,\quad \lambda\in {\mathcal D}\, ,
\label{restriction map}
\end{equation} 
and $\kappa = 1$ in a neighborhood of $(x^0,\xi^0)$. In particular, this gives $m=1/4$ in \eqref{oscillatory1}.  

In what follows, we shall investigate the composition $\imath_\Gamma^\ast I(\lambda)$ of $\lambda$-FIOs, where $I(\lambda)\in I^{m}(\tilde X, \Gamma, {\mathcal C}; \Omega^{1/2})$. Firstly, notice that the composition ${\cal R}\circ {\mathcal C}$ of the corresponding canonical relations is transversal. Indeed, introduce 
\[
Z:= T^\ast \Gamma \times \Delta \left(T^\ast \widetilde X \times T^\ast \widetilde X \right) \times T^\ast \Gamma,
\]
where $\Delta$ stands for the diagonal, and consider $Z$, ${\mathcal R}\times {\mathcal C}$ and $({\mathcal R}\times {\mathcal C})\cap Z$ as submanifolds of 
$Y:=T^\ast \Gamma \times T^\ast \widetilde X \times T^\ast \widetilde X  \times T^\ast \Gamma$. We have ${\rm codim}\, (Z)= 2n$, ${\rm codim}\, ({\mathcal R}\times {\mathcal C})= 4n-2$, and  ${\rm codim}\, (({\mathcal R}\times {\mathcal C})\cap Z)= 6n-2$, which implies that the intersection of $Z$ and ${\mathcal R}\times {\mathcal C}$ is transversal along $({\mathcal R}\times {\mathcal C})\cap Z$ in $Y$. Denote by $\sigma(I) =   (\lambda/2\pi)^m e^{i\lambda f} \widetilde b_0\,  \sigma_0\otimes \sigma_2$ the principal symbol of $I(\lambda)$, where $\widetilde b_0 \in C^\infty_0({\mathcal C}')$ and the half-density $\sigma_0$ is given by \eqref{density}. 
Recall that $\pi_1: {\mathcal C} \to T^\ast X$ and $\pi_2: {\mathcal C} \to T^\ast \Gamma$ are the projections 
$\pi_1(x,y,\xi,\eta)= (x,\xi)$ and $\pi_2(x,y,\xi,\eta)= (y,\eta)$. Denote by $dv(\varrho):= dy\wedge d\eta$  the symplectic volume form on $T^\ast\Gamma$, and  recall that $\nu(x)\in T_x\widetilde X$ is the unit inward normal to $\Gamma$, and that $\pi_\Sigma^\pm (x,\xi)= (x,\xi^\pm)\in \Sigma^{\pm} $ (see Sect. \ref{subsec:billiard-ball}).  Now we have 
\begin{Lemma}\label{Lemma:restriction}
The composition of canonical relations 
${\cal R}\circ {\mathcal C}$ is transversal    and it is a disjoint union $\Delta_0 \sqcup {\mathcal C}_0$ of
the diagonal $\Delta_0$ in $U\times U$ (for $s=0$) and  the graph ${\mathcal C}_0$ of the billiard ball map
$B:U\to B(U)$ (for $s=T$). 
 Moreover, 
\begin{equation}
\imath_\Gamma^\ast I(\lambda) = P(\lambda) + 
G(\lambda)
+ O_M(|\lambda|^{-M})
\, ,  
                     \label{boundary-trace1}
\end{equation}
where  $P(\lambda)$ is a classical $\lambda$-PDO on $\Gamma$ of order $m+1/4$ and $G(\lambda)\in I^{m+1/4}(\Gamma, \Gamma,  {\mathcal C}_0; \Omega^{1/2})$.   The principal symbol of the operator 
$P(\lambda)$ can be identified by 
\[
\left(\frac{\lambda}{2\pi}\right)^{m+1/4}\widetilde b_0( \pi_2^{-1}( \varrho))(2| \langle \xi^+,\nu\rangle(\varrho) |)^{-1/2}|dv(\varrho)|^{1/2}, \quad \varrho \in U.  
\] 
The principal symbol of  $G(\lambda)$  can be identified with  
\[
\left(\frac{\lambda}{2\pi}\right)^{m+1/4}\widetilde b_0( \pi_1^{-1}( \pi_\Sigma^-(\varrho)))|2\langle\xi^-,\nu\rangle(\varrho)|^{-1/2}e^{i\lambda A(B^{-1}(\varrho))}e^{i \pi m/2}|dv(\varrho)|^{1/2}, \quad \varrho \in B(U), 
\]
where 
$A(\varrho)=
\int_{\gamma(\varrho)}\xi dx $ is the action along the integral curve $\gamma(\varrho)$ of the Hamiltonian vector
field $X_{\widetilde h}$  starting from  $\pi^+_\Sigma(\varrho)$ with endpoint  $\pi^{-}_\Sigma(B(\varrho))$ and  $m\in \Z$ is a Maslov's index. 
\end{Lemma}
{\em Proof}. Let $(\varrho',\varrho)\in {\cal R}\circ {\mathcal C}$. Then there is $s\ge 0$ such that 
$\exp(sH_{\widetilde h}(\pi^+_\Sigma(\varrho)))\in \Sigma$, hence, $\varrho'=\varrho$  if $s=0$ and $J\circ \pi^+_\Sigma(\varrho) = \pi^-_\Sigma(\varrho')$ if $s=T(\varrho)$ and we get $\varrho'= B(\varrho)$ in the second case (see Sect. \ref{subsec:billiard-ball}). 
Hence, ${\cal R}\circ {\mathcal C} = \Delta_0 \sqcup {\mathcal C}_0$, where  $\Delta_0= \{(y,\eta,y,\eta):\, (y,\eta)\in U\}$  and  
${\mathcal C}_0= \{(x,\xi,y,\eta):\, (x,\xi)=B(y,\eta),\, (y,\eta)\in U\}$.

The statement about the principal symbol follows  from the global calculus of $\lambda$-FIOs. 
It can be proved easily as well using the calculus in local coordinates.  For pedagogical reason we give the complete proof for readers who are not well familiar with the global calculus of $\lambda$-FIOs. 
Take $\nu^0=(x^0,y^0,\xi^0,\eta^0)\in {\mathcal C}_0'$ so that $(x^0,\xi^0)= B(y^0,\eta^0)$. 
 Choose as above normal to $\Gamma$  coordinates $(x',x_n)$  in a neighborhood of $x^0$ such that $x_n<0$ in the interior of $X$. Using Lemma \ref{lamma:phase-function}, choose   local coordinates $y'$ in a neighborhood of $y^0$ in $\Gamma$ and a non-degenerate phase function $\Phi$ of the form $\Phi (x,y', \theta)= \phi (x,\theta) - \langle y',\theta\rangle$  which defines locally ${\mathcal C}'$ near $(\nu^0)^- = (x^0,y^0, (\xi^0)^-,\eta^0)$, where $(x^0,(\xi^0)^-)= \pi_\Sigma^-(x^0,\xi^0)$, 	and then take normal to $\Gamma$ coordinates $y=(y',y_n)$ such that $y_n>0$ in the interior of $X$. Using \eqref{oscillatory1} we write microlocally the Schwartz kernel of $I(\lambda)$ as an oscillatory integral with a phase function $\Phi(x,y',\theta)$ and amplitude 
 $(\lambda/2\pi)^{m+n-3/4}b(x,\theta,\lambda)|dx|^{1/2}|d y'|^{1/2}$, where $b=b_0 + \lambda^{-1}b_1 + \cdots $ is a classical amplitude. 
Then  the Schwartz kernel of the composition can be written microlocally near $\nu^0$ as a half-density $ \widetilde K_1 := K_1 |dx'|^{1/2}|d y'|^{1/2}$, where $K_1$ is given by the oscillatory integral 
\[
K_1(x',y',\lambda):= \left(\frac{\lambda}{2\pi}\right)^{m+2n-3/4} \int_{\R^{3n-1}} \, e^{i\lambda\Phi_1(x',y',z,\xi,\theta )} a(x', z,\xi,\theta,\lambda)\, dz d\xi d\theta  \, ,
\] 
with  phase function  $\Phi_1(x',y',z,\xi,\theta )= \langle x'-z',\xi'\rangle -z_n\xi_n + \phi(z,\theta)-\langle y',\theta\rangle$ and amplitude $a$,    the leading term of which is 
$ a_0(x', z,\xi,\theta) = \kappa(x',\xi)b_0(z,\theta)$. The stationary points with respect to $(z,\xi)$ are given by $z'=x'$, $z_n=0$, and $\xi = \phi'_{z}(z,\theta)$,  and they are non-degenerate. 
Applying the stationary phase method one gets 
\begin{equation}\label{stationary-phase}
 K_1(x',y',\lambda):=\left(\frac{\lambda}{2\pi}\right)^{m+1/4+(n-1)} \int_{\R^{n-1}} \, e^{i\lambda \Phi_0(x',y',\theta )} w(x',\theta,\lambda)\,  d\theta \, ,
\end{equation} 
where $\Phi_0(x',y',\theta )=\phi(x',0,\theta )-\langle y',\theta\rangle$, and $w(x',\theta,\lambda)$ is a classical symbol of order $0$ with a leading term
\[
w_0(x',\theta) = b_0(x',0,\theta)\kappa(x',0, \phi'_{x}(x',\theta))  .
\]
The phase function $\Phi_0(x',y', \theta)$ generates the Lagrangian submanifold 
${\mathcal C}_0'$ 
of $T^\ast \Gamma$ in a neighborhood of $\nu^0$,  since  $J\circ \pi^+_\Sigma(\phi'_{\theta}(x',0,  \theta),\theta) = (x', 0, \phi'_{x}(x',0,\theta)) = \pi^-_\Sigma(x',  \phi'_{x'}(x',0,\theta))$. 
The half-density part of the principal symbol of $ \widetilde K_1$ can be identified  in a neighborhood of $\nu^0$ in  ${\mathcal C}'_0$ with 
\[
\sigma:= b_0(x,0,\theta)  |dx'\wedge d\theta|^{1/2} = 
b_0(x',0,\theta(x',\xi'))  \left|\left(\det \phi''_{x'\theta}(x',0,\theta(x',\xi'))\right)\right|^{-1/2} |dx'\wedge d\xi'|^{1/2},
\]
where $\theta=\theta(x',\xi')$ is the solution of $\xi'=\phi'_{x'}(x',0, \theta)$. 
Then setting $(  x',0,  \xi^-):= \pi^-_\Sigma(x',\xi')$, where $\xi^-=(\xi',\xi_n)$ and $\xi_n = -\sqrt{1-h_0(x',0,\xi')} = \langle \xi^-,\nu\rangle(x',\xi')$,    and taking into account \eqref{principal-symbol3},  we obtain
\[
\sigma =  
\widetilde b_0( \pi_1^{-1}( x',0,  \xi^-)) (2|  \langle \xi^-,\nu\rangle(x',\xi')|)^{-1/2}|d  x'\wedge d    \xi'|^{1/2} . 
\]
On the other hand, the trivialization of the Keller-Maslov bundle of ${\mathcal C}'$ by $\sigma_2$ induces a trivialization of the Keller-Maslov bundle of ${\mathcal C}_0'$ and we can identify the principal symbol of $G(\lambda)$ with 
\[
\left(\frac{\lambda}{2\pi}\right)^{m+1/4} e^{i\lambda A(B^{-1}(\varrho))}e^{i \pi m/2}\sigma , 
\]
where $m\in\Z$.

Consider now  the case when $s=0$, taking the phase function $\Phi$ in \eqref{phase-function}  so that $\phi(x',0,\eta') = \langle x',\eta'\rangle$. The corresponding phase function of \eqref{stationary-phase} is  $\Phi_0(x', y',\eta')= \langle x'-y',\eta'\rangle$ and \eqref{stationary-phase} becomes 
\[
 K_0(x',y',\lambda):=\left(\frac{\lambda}{2\pi}\right)^{m+ 1/4+(n-1)} \int_{\R^{n-1}} \, e^{i\lambda \langle x'-y',\eta'\rangle} q(x',\eta',\lambda)\,  d\eta' \, ,
\]
where $q(x',\eta',\lambda)$ is a classical symbol of order $0$ with a principal term  $q_0(x',\eta')$. Setting  $(  y',0,  \eta^+):= \pi^+_\Sigma(y',\eta')$, where 
$\eta^+=(\eta',\eta_n)$ and  $\eta_n= \sqrt{1-h_0(y',0,\eta')}= \langle \eta^+,\nu\rangle(y',\eta')$ 
 we obtain as above
\[
q_0(x',\eta') = \widetilde b_0( \pi_1^{-1}( y',0,\eta^+)(2| \langle \eta^+,\nu\rangle(y',\eta') |)^{-1/2}. 
\] 
%Choosing a section of the Keller-Maslov bundle of  ${\mathcal C}_1'$ we obtain the Maslov factor $e^{i\pi m/2}$. 
\finishproof

Applying Lemma \ref{Lemma:restriction} to the $\lambda$-FIO $H_0(\lambda)$, which is of order $m=-1/4$,  we shall complete the construction of $H(\lambda)$. 
Let  $\Psi(\lambda)$ be a classical $\lambda$-PDO of order $0$ with frequency set in $ U$ and principal symbol $\Psi_0(\varrho)$, $\varrho\in U$. 
 We take $\Psi(\lambda)$ as initial data of $H_0(\lambda)$ as $s=0$, setting $P(\lambda)=\Psi(\lambda)$ in Lemma \ref{Lemma:restriction}. Recall that $\widetilde b_0$ satisfies \eqref{transport-equation}. 
On the other hand $(x',\xi')=B(y',\eta')$ if and only if $\pi_1^{-1}( x',0,\xi^-)=   S^{T( y',\eta')}(\pi_2^{-1}( y',\eta'))$, and 
 \eqref{transport-equation} implies
\[
\widetilde b_0( \pi_1^{-1}( x',\xi',\xi_n)) = \widetilde b_0( (\pi_2^{-1}( y',\eta'))) = q_0(y',\eta')(2|\langle \eta^+,\nu\rangle(y',\eta')|)^{1/2}.
\]
Then  parameterizing ${\mathcal C}_0'$ by the variables $(y',\eta')\in U$ we obtain
 \begin{equation}
\sigma(H_0(\lambda)) =  
\Psi_0(y',\eta') \frac{|\langle\eta^+,\nu\rangle(y',\eta')|^{1/2}}{|\langle \xi^-,\nu\rangle(B(y',\eta'))|^{1/2}} \, e^{i\lambda A(y',\eta')}e^{i \pi m/2} | d  y'\wedge d \eta'  |^{1/2} . 
                                                   \label{principal-symbol}
\end{equation}
In the same way, using Lemma \ref{Lemma:restriction} we  determine the initial conditions of $H_1(\lambda)$ and so on. In this way  we obtain an operator 
$H(\lambda)=H_0(\lambda) + H_1(\lambda) + \cdots $ in $ I^{-1/4}(\tilde X, \Gamma, {\mathcal C}; \Omega^{1/2})$ satisfying \eqref{parametrix1} and such that $P(\lambda)=\Psi(\lambda)$. From now on, to simplify the notations we drop the corresponding half-density.
Denote by $E(\lambda)$ a classical $\lambda$-PDO of order $0$ on $\Gamma$ with a principal symbol $E_0\in C_0^\infty(\widetilde B^\ast \Gamma)$ such that  
\begin{equation}
E_0(\varrho) = |\langle \xi^+,\nu\rangle(\varrho)|^{1/2} = |\langle \xi^-,\nu\rangle(\varrho)|^{1/2}
                    \label{boundary-trace5}
\end{equation}
in a compact neighborhood of $\overline U$ in $\widetilde B^\ast \Gamma$.  Then using Egorov's theorem,  \eqref{principal-symbol} and \eqref{boundary-trace5} we obtain
\begin{equation}
\imath_\Gamma^\ast H(\lambda) = \Psi(\lambda) + 
E(\lambda)^{-1}G^0(\lambda)E(\lambda)
+ O_M(|\lambda|^{-M})
\, ,  
                     \label{boundary-trace3}
\end{equation}
where the principal symbol of $G^0(\lambda)$ can be identified with 
\begin{equation}
\Psi_0(\varrho)e^{i\lambda A(\varrho)}e^{i \pi m/2}|dv(\varrho)|^{1/2}, \quad \varrho \in U.  
                      \label{boundary-trace4}
\end{equation}
In particular,  the frequency 
set $WF'$  of
$G(\lambda)= E(\lambda)^{-1}G^0(\lambda)E(\lambda)$  
is contained in 
$B(U)\times U$.  

Consider the adjoint operator $H(\lambda)^\ast$ of $H(\lambda)$ in $L^2$, which is well-defined for any $\lambda\in {\mathcal D}$ fixed as an operator
from  $L^2(\widetilde X)$ to $L^2(\Gamma)$. Moreover, it can be considered as a $\lambda$-FIO of the class $ I^{-1/4}(\Gamma, \tilde X,  {\mathcal C}^{-1}; \Omega^{1/2})$. 
\begin{Prop}\label{Lemma:continuity}
The operator $C(\lambda):=H(\lambda)^\ast H(\lambda): L^2(\Gamma) \to L^2(\Gamma)$ is a classical $\lambda$-PDO of order $0$. Its principal symbol can be identified with 
$$ C_0(y,\eta):=\int_\R |\widetilde b_0 (s,y,\eta)|^2 ds \, ,\quad (y,\eta)\in U.$$
Moreover, 
$$ C_0(y,\eta)\ge 2 A(y,\eta) | \langle \xi^+,\nu\rangle(y,\eta) |\Psi_0(y,\eta)|^2\, ,\quad (y,\eta)\in U. $$
\end{Prop}
Indeed, the composition ${\mathcal C}^{-1}\circ {\mathcal C}$ of the canonical relations ${\mathcal C}^{-1}$ and ${\mathcal C}$ is clean with excess $1$, which means that the ``diagonal'' $M_1:=T^\ast \Gamma \times \Delta(T^\ast \widetilde X \times T^\ast \widetilde X)\times T^\ast \Gamma$ intersect cleanly $M_2:={\mathcal C}^{-1}\times {\mathcal C}$  along $M_1\cap M_2$ in $T^\ast \Gamma \times T^\ast \widetilde X \times T^\ast \widetilde X\times T^\ast \Gamma$ with excess $1$. 
The fiber over  any $(\varrho,\varrho)\in {\mathcal C}^{-1}\circ {\mathcal C}$ under the corresponding  projection  $M_1\cap M_2\to {\mathcal C}^{-1}\circ {\mathcal C}$ can be identified with the Hamiltonian trajectory  starting from  $\pi^+_\Sigma(\varrho)$. Moreover, 
$\widetilde b_0(s,y,\eta)= (2| \langle \xi^+,\nu\rangle(y,\eta) |)^{1/2}\Psi_0(y,\eta)$, for $s\in [0, T(y,\eta)]$ and $T(y,\eta)= 2A(y,\eta)$, which proves the claim.  One can prove the Proposition  using only the local theory of $\lambda$-FIOs representing microlocally $H(\lambda)$ by an oscillatory integral with a phase function given by \eqref{phase-function}. 
\finishproof

\subsection*{A.2. Symbols of finite smoothness.}\label{subsec:finite-smoothness}
We start this section by recalling a result of  Boulkhemair  \cite{Boulk} about the $L^2$  boundedness of  a class of  $\lambda$-FIOs. 
Consider  a $\lambda$-FIO $A_\lambda$, $\lambda\in {\mathcal D}$,  acting on $C^\infty_0(\R^{n-1})$ by
\begin{equation}
A_\lambda u(x)=(\lambda/2\pi)^{n-1} \int e^{i\lambda (\langle x,\xi\rangle+\psi(x,\xi))} q(x,\xi,\lambda)\, \hat{u}(\lambda\xi)\, d\xi \, ,
\label{FIO1}
\end{equation}
where ${\mathcal D}$ is given by \eqref{strip} and $\hat u$ is the Fourier transform of $u$. 
The distribution kernel of $A_\lambda$ is 
\begin{equation}
K_{A_\lambda}(x,y)=(\lambda/2\pi)^{n-1} \int e^{i\lambda (\langle x-y,\xi\rangle+\psi(x,\xi))} q(x,\xi,\lambda)\, d\xi \, . 
\label{FIO}
\end{equation}
Suppose that the amplitude $q$ satisfies the following conditions
\begin{enumerate}
\item[(i)]
$q_\lambda:=q(\cdot,\cdot,\lambda) \in C^n(\R^{n-1}\times \R^{n-1})$ for any $\lambda \in \cal D$, 
\item[(ii)]  $\displaystyle\sup_{\lambda\in \cal D} \|q_\lambda\|_{C^n(\R^{2n-2})}  <\infty $, 
\item [(iii)]
there is compact $F\subset T^\ast R^{n-1}$ such that for any $\lambda\in  {\cal D}$
the  support of the function $q_\lambda$  is contained in $F$. 
\end{enumerate}
Suppose also that the phase function $S(x,\xi) = \langle x,\xi\rangle+\psi(x,\xi)$ is $C^\infty$-smooth 
in a neighborhood $U$ of $ F$ and that 
 $|\det (\partial_x\partial_\xi S)| \ge \delta >0$ in $U$.
 
Using a result of  Boulkhemair  \cite{Boulk},   we are going to show that
\begin{equation}
\|A_\lambda\|_{{\cal L}(L^2)}\ \le\  C\,  \sup_{\lambda\in {\mathcal D}} \|q_\lambda\|_{C^n}  \, ,
                              \label{L2estimates}
\end{equation}
where $C=C(S,F)>0$ does not depend on $q_\lambda$ and ${\cal L}(L^2)$ stands for the space of linear continuous operators in $L^2(\R^{n-1})$. To this end we   extend the phase function $S$ in $T^\ast\R^{n-1}$ so that $|\det (\partial_x\partial_\xi S)| \ge \delta/2$ in the whole space. 
Indeed,  let  $\varrho^0\in F$ and 
$$
F\subset B_\varepsilon (\varrho^0):=\{\varrho\in T^\ast\R^{n-1}:\, |\varrho-\varrho^0|<\varepsilon\}
\subset B_{2\varepsilon} (\varrho^0)\subset U.
$$ 
If $\varepsilon>0$ is sufficiently small we can
extend    $S$ to a globally defined smooth function $\widetilde S$ in $T^\ast\R^{n-1}$
which coincides with $S$ in $B_\varepsilon (\varrho^0)$ and equals the Taylor polynomial of 
degree $2$ of $S(x,\xi)$ at
$\varrho^0$ outside $B_{2\varepsilon}(\varrho^0)$ and such that $|\det  \partial_x\partial_\xi\widetilde S| \ge \delta/2$ in $T^\ast \R^{n-1}$. Now applying 
 \cite{Boulk}, Corollary 1, to the oscillatory integral with phase function $\widetilde S$ and
amplitude $q_\lambda$ we obtain 
\[
\|A_\lambda u\|_{L^2}\ \le\  C\, (|\lambda|/2\pi)^{(n-1)/2}  \sup_{\lambda\in {\mathcal D}} \|q_\lambda\|_{C^n}\|\hat u(\lambda \, \cdot)\|_{L^2} =  C\,  \sup_{\lambda\in {\mathcal D}} \|q_\lambda\|_{C^n}\|u\|_{L^2}\, .
\]
In the general case we use a suitable partition of
the unity on $F$.

 Consider now the oscillating integral 
\begin{equation}
q_\lambda(x,\xi) = \lambda ^{n-1} \int_{\R^{2n-2}} \, 
e^{-i\lambda \langle v,u\rangle}
a_\lambda(x,\xi,u,v)\,  dv du\ , \quad  
\lambda\in {\cal D}\, ,
\label{oscillating-integral}
\end{equation}
where $a_\lambda =a(\cdot,\lambda)\in C^{2n}(\R^{n-1}\times \R^{n-1}\times\R^{n-1}\times \R^{n-1})$, and 
\begin{enumerate}
\item[(a)] $\sup_{\lambda\in {\mathcal D}} \|a_\lambda\|_{C^{2n}}<\infty$, 
\item[(b)]
 $a_\lambda$ is  uniformly compactly supported with respect to $v\in \R^{n-1}$, 
which means that there is a  
compact subset $F$  of $\R^{n-1}$ such that the support of the function 
$v\to a_\lambda(x,\xi,u, v)$ is contained in $F$ for any
$(x,\xi,u,\lambda)\in \R^{n-1}\times \R^{n-1}\times \R^{n-1}\times {\cal D}$.  
\end{enumerate}
\begin{Lemma}\label{Lemma:estimate1}
Let $a_\lambda\in C^{2n} $  satisfy (a) and (b). 
Then 
\[
 \|q_\lambda\|_{C^{n}}\,  \le\,  C 
 \|a_\lambda\|_{C^{2n}}\,   , 
\]
where $C=C(n, F)>0$. Moreover, if $a_\lambda$ is uniformly compactly supported with respect to $(x,\xi,v)$, then the operator \eqref{FIO1} is bounded in $L^2$ and 
$$\|A_\lambda\|_{{\cal L}(L^2)}\ \le\  C\,  \sup_{\lambda\in {\mathcal D}} \|q_\lambda\|_{C^n} \le   C' \sup_{\lambda\in {\mathcal D}} \|a_\lambda\|_{C^{2n}}. $$ 
\end{Lemma}
{\em Proof}. 
%Choose $\psi\in C_0^\infty(\R^{n-1})$ such that $\sum_{k\in\Z^{n-1}}\psi(u+2\pi k) = 1$ for any $u\in\R^{n-1}$, and then set
%\[
%q_{k,\lambda}(\varphi,I) = \lambda ^{n-1} \int_{\R^{2n-2}} \, 
%e^{-i\lambda \langle v,u\rangle}
%\tilde a(\varphi,u,I,v,\lambda)\psi(u+2\pi k)\,  dv du\, .
%\]
Set $z=(x, \xi)$. For any $\alpha\in \N^{2n-2}$ with $|\alpha|\le n$ we get
\[
\displaystyle \partial^\alpha_z q_{\lambda}(x,\xi) =  \int_{\R^{n-1}} \, 
b_\lambda(x,\xi, u/\lambda,u)\,  du\, ,
\]
where 
\[
b_\lambda(x,\xi,w, u) = \, \int_{\R^{n-1}} \, 
e^{-i\langle v,u\rangle}
\partial^\alpha_z  a_\lambda(x,\xi,w, v)\,  dv \, .
\]
Integrating  by parts with respect to $v$ in the oscillating integral above we get
\[
(1+|u|)^n |b_\lambda(x,\xi,w, u)| \le C \|a_\lambda\|_{C^{2n}}, 
\]
where $C=C(n, F)>0$. The last assertion of the Lemma follows from \eqref{L2estimates}.  
\finishproof

We consider families of $\lambda$-PDOs with symbols of finite smoothness which depend 
continuously on  $K\in C^\ell(\Gamma)$.
Let $Y$ be a smooth paracompact manifold of dimension $n-1$, $n\ge 2$.
 Let ${\mathcal B}= {\mathcal B}^\ell$ be  a bounded subset of $C^{\ell}(\Gamma,\R)$.
We denote by $O_{\mathcal B}(|\lambda|^{-N})\, :\ L^2(Y)
\to  L^2(Y)$, $\lambda\in {\cal D}$,  any family of linear continuous operators depending on $K\in {\cal B}$,
the norm of which is {\em uniformly bounded} by $C_{\mathcal B}(1 +|\lambda|)^{-N}$ for some positive  constant $C_{\mathcal B}$. 
\begin{Def}\label{Def:finite-smoothness}
Let  $\widetilde l \ge 0$,  $N\ge 1$ and 
$m\ge 2$  be such that $\widetilde l \ge mN + 2n$.   
We say  that a family of operators $Q(\lambda, K) \, :\ C^\infty(Y) \to  C^\infty (Y)$, $\lambda\in {\cal D}$,  depending on $K\in {\cal B}^\ell\subset C^{\ell}(\Gamma,\R)$  belongs to
${\rm PDO}_{\widetilde l,m,N}(Y; {\mathcal B}^\ell; \lambda)$ if in any local coordinates
$Q(\lambda,K)$ can be written in the form ${\rm OP}_\lambda (q) + O_{\mathcal B}(|\lambda|^{-N})$, where
the Schwartz  kernel of ${\rm OP}_\lambda (q)$ is  
\begin{equation}
{\rm OP}_\lambda(q)(x,y):=(\lambda/2\pi)^{n-1} \int e^{i\lambda \langle x-y,\xi\rangle} q(x,\xi,\lambda)\, d\xi \, ,
\label{pdo}
\end{equation} 
 with amplitude 
\begin{equation}
q(x,\xi,\lambda) = \sum_{k=0}^{N-1} q_k(x,\xi)\lambda^{-k}\, ,
                                      \label{thesymbol}
\end{equation} 
and $q_k\in C^{\widetilde l-mk}_0(T^\ast R^{n-1})$, $0\le k \le N-1$,  depends continuously in $K\in  C^\ell(\Gamma,\R)$ in the sense
that  the support of $q_k$ is contained in a fixed compact set independent of $K$ 
and the map 
\[
C^\ell(\Gamma,\R)\ni K  \ \to \ q_k\in C^{\widetilde l-mk}(T^\ast R^{n-1})
\]
is continuous. 
We denote the class of symbols $q$ by $S_{\widetilde l,m,N}(T^\ast \R^{n-1}; {\mathcal B}^\ell; \lambda)$. 
\end{Def}
Throughout the paper we take $Y$ to be either $\Gamma$ or  $\T^{n-1}$ and $\R^{n-1}$. To simplify the notations we  write ${\mathcal B}$ and $Q(\lambda)$ instead of ${\mathcal B}^\ell$ and $Q(\lambda,K)$ respectively. 
\begin{Remark}\label{Rem:reminder}
Let $q_{\lambda,K}(\cdot,\cdot)\in  C^n(T^\ast \R^{n-1})$,  $\lambda\in {\cal D}$, be  a family of symbols depending on  $K\in {\cal B}$.  Suppose that 
$\displaystyle\sup_{\lambda\in {\mathcal D},\, K\in {\mathcal B}} \|q_{\lambda,K}\|_{C^n(T^\ast\R^{n-1})}  < C_{\mathcal B} $  and that 
there is a compact $F\subset T^\ast R^{n-1}$ such that the  support of the function $q_{\lambda,K}$  is contained in $F$ for any $\lambda\in  {\cal D}$ and $K\in {\cal B}$. Then by \eqref{L2estimates}  the family of operators ${\rm OP}_\lambda (\lambda^{-N}q_{\lambda,K})$ is $O_{\mathcal B}(|\lambda|^{-N})\, :\ L^2( \R^{n-1})
\to  L^2(\R^{n-1})$. 
\end{Remark}
In particular,  we obtain (see also the $L^2$-continuity theorem, \cite{Hor}, Theorem 18.1.11$'$) that 
\begin{Lemma}\label{Lemma:L2-continuity}
Let  $Y$ be  compact and $\widetilde l \ge m(N-1) + n$. Than any  family of operators $Q(\lambda,K)$, $\lambda\in {\cal D}$,   in ${\rm PDO}_{\widetilde l,m,N}(Y; {\mathcal B}; \lambda)$ is  uniformly bounded in $L^2$ with respect to $K\in {\cal B}$ and $\lambda\in {\cal D}$, i.e. there exists $C_{\mathcal B}>0$ such that $\|Q(\lambda,K)\|_{{\cal L}(L^2)} \le  C_{\mathcal B}$ for any $K\in {\cal B}$ and any $\lambda\in {\cal D}$. 
\end{Lemma}
We shall see in Remark \ref{Rem:commutator} that  the
class $ {\rm PDO}_{\widetilde l,m,N}(Y; {\mathcal B}; \lambda)$, $\widetilde l \ge mN + 2n$,  is closed under multiplication and 
transposition   and that  it 
does not depend on the choice of the local coordinates.

The frequency  set  $\mbox{WF}'(Q)$  (modulo $O(|\lambda|^{-N})$) of a $\lambda$-PDO $Q(\lambda)$ 
with symbol $q$ locally given by (\ref{thesymbol}) is 
\[
\mbox{WF}'\, (Q(\lambda)) := \cup_{j=0}^{N-1} \, \mbox{supp}\, (q_j)
\]
in each local cart (the essential support of the symbol).

Using  Lemma \ref{Lemma:estimate1} we are going to investigate the composition of    $\lambda$-PDOs in 
${\rm PDO}_{\widetilde l,s,N}( Y,{\mathcal B};\lambda)$
 with  classical 
$\lambda$-FIOs $A(\lambda)$ associated to a smooth canonical transformation 
$\kappa:T^\ast  Y \to T^\ast  Y$ and having  $C_0^\infty$ amplitudes in each local cart. 
\begin{Lemma}\label{Lemma:commutator}
If  $Q(\lambda) \in   {\rm PDO}_{\widetilde l,m,N}( Y; {\mathcal B}; \lambda)$, 
$\widetilde l \ge mN+ 2n$, and the $\lambda$-FIO  $A(\lambda)$ is elliptic on $\mbox{WF}'(Q)$, then 
there exists $Q'(\lambda) \in   {\rm PDO}_{\widetilde l,m,N}( Y; {\mathcal B}; \lambda)$ such that 
\begin{equation}
Q(\lambda) A(\lambda) - A(\lambda)Q'(\lambda) =  O_{\mathcal B}(|\lambda|^{-N})\, :\ L^2( Y)
\longrightarrow  L^2( Y)
\label{commutation}
\end{equation}
 and wise versa, if $Q'(\lambda) \in   {\rm PDO}_{\widetilde l,m,N}( Y; {\mathcal B}; \lambda)$ and $A(\lambda)$ is elliptic on $\mbox{WF}'(Q')$, then 
there exists $Q(\lambda) \in   {\rm PDO}_{\widetilde l,m,N}( Y; {\mathcal B}; \lambda)$ such that \eqref{commutation} holds.
Moreover, the  principal symbols of $Q(\lambda)$ and $Q'(\lambda)$  are related  by the Egorov's 
theorem, $\sigma(Q') = \sigma(Q)\circ \kappa$.
\end{Lemma}
{\em Proof}. We define $Q' = BQA$, where
$\mbox{WF}'(AB -I)\cap \mbox{WF}'(Q) = \emptyset$. To prove that $Q'(\lambda) \in   {\rm PDO}_{\widetilde l,m,N}( Y; {\mathcal B}; \lambda)$, we 
choose  local coordinates $x$ in $ Y$ and write the distribution kernel of 
 $Q(\lambda)$ in the form (\ref{pdo}) with symbol 
 $q\in S_{\widetilde l,m,N}(T^\ast \R^{n-1}; {\mathcal B}; \lambda)$.  
 We can suppose that the distribution kernel of 
 $A(\lambda)$  is given microlocally  by (\ref{FIO}), where $\langle x,\xi\rangle + \psi(x,\xi)$ is a   generating function of the symplectic transformation $\kappa$ and  $a$ is a smooth compactly supported amplitude. 
 
 More generally, we suppose that
 $a\in S_{\widetilde l,m,N}(T^\ast \R^{n-1}; {\mathcal B}; \lambda)$.  Then the 
 distribution kernel of $Q(\lambda)A(\lambda)$ 
 is given modulo $O_{\cal B}(|\lambda|^{-N})$ by the oscillatory integral (\ref{FIO}) with
 amplitude given by the oscillatory integral
\[
 \displaystyle K_{QA}(x,\xi,\lambda)=
\left(\frac{\lambda}{2\pi}\right)^{n-1}\, \int_{{\R}^{2n-2}} 
\, e^{i\lambda (\langle x-z,\eta-\xi\rangle + \psi(z,\xi) -\psi(x,\xi))} \, 
q(x,\eta,\lambda)
a(z,\xi,\lambda)\, d\eta dz \, .
\] 
The amplitude of the  oscillatory integral $K_{QA}$ is uniformly compactly supported with respect to  
$(x,z,\xi,\eta)$.
We are going to write it as a symbol in $S_{\widetilde l,m,N}(T^\ast \R^{n-1}; {\mathcal B}; \lambda)$. Set  
\[ 
\psi_1(x,z,\xi) = \int_0^1 \nabla_x \psi(x +\tau z,\xi) d \tau\, .
\]
Changing the variables and using Lemma \ref{Lemma:estimate1} we get 
\[
K_{QA}(x,\xi,\lambda) = \sum_{j=0}^{N-1} \sum_{r + s = j}\lambda^{-j} \left(\frac{\lambda}{2\pi}\right)^{n-1}\,  \int_{{\R}^{2n-2}} 
\, e^{-i\lambda \langle z,\eta\rangle } A_{r,s}(x,\xi,z,\eta)\,  \, d\eta dz    + O_{\cal B}(|\lambda|^{-N})\, ,
\] 
where the amplitude 
\[
A_{r,s}(x,\xi,z,\eta) = 
 q_r(x,\eta +\xi + \psi_1(x,z,\xi))
 a_{s}(z +x, \xi)   
\]
has a compact support. 
We expand  $A_{r,s}$ by Taylor formula with respect to $\eta$ 
at $\eta=0$  up to order $O(|\eta|^{N-j})$ 
\[
\begin{array}{lcrr}
\displaystyle A_{r,s}(x,\xi,z,\eta)\,  = \, \sum_{|\beta|<N-j}\frac{\eta^\beta}{\beta!}\, 
(\partial_\eta^\beta A_{r,s})(x,\xi,z,0) \\[0.3cm]
\displaystyle  +\,  (N-j)\int_{0}^{1} (1-t)^{N-j-1}\sum_{|\beta|=N-j}\frac{\eta^\beta}{\beta!}\, (\partial_\eta^\beta A_{r,s})(x,\xi,z,t\eta) dt. 
\end{array}
\]
Integrating $\beta$ times by parts with respect to $z$ in the corresponding oscillatory integrals  
we obtain  
$$
K_{QA}(x,\xi,\lambda) = 
\sum_{j=0}^{N-1} F_{j}(x,\xi)\lambda^{-j} + F_{N}(x,\xi,\lambda)\lambda^{-N}\, ,
$$ 
where 
\begin{equation}
\begin{array}{lcrr}
\displaystyle F_{j}(x,\xi)\ =\ 
\displaystyle \sum_{r+s +|\beta|= j}\, 
\frac{1}{\beta!}\, \left[ \left(-i\partial_z\right)^\beta  \left(\partial^\beta_\eta\,
q_r(x,\eta +\xi + \psi_1(x,z,\xi)\, 
a_s(z+x,\xi )\right)\right]_{|z=0, \eta =0}\,  
                    \label{symbolQA} 
\end{array}
\end{equation}
for $j\le N-1$. Hence,  $F_j\in C_0^{\tilde l -mj}(T^\ast \R^{n-1})$, $0\le j\le N-1$,  and it depends continuously on $K$. 

Consider now the reminder term $F_N$. Notice that for $|\beta|\le N-j\le N-r$ we have $\partial_\eta^\beta q_{r}\in C^{\tilde l -mr-|\beta|}(T^\ast \R^{n-1})$, 
where  
\begin{equation}
\tilde l -mr-|\beta| = |\beta| + \tilde l -mr-2|\beta| \ge   |\beta|+\tilde l -m(r + |\beta|) \ge  
|\beta|+ \tilde l -mN \ge |\beta|+ 2n.
\label{order}
\end{equation}
Then $\partial_z^{\gamma}\partial_\eta^{\beta} q_{r}\in C^{2n}(T^\ast \R^{n-1})$ if $|\gamma|\le |\beta|\le N-j$. In the same way we get $\partial_z^{\gamma} a_{s}\in C^{2n}(T^\ast \R^{n-1})$ 
if $|\gamma| \le |\beta|\le N-j$. 
Now using  Lemma \ref{Lemma:estimate1} (with $v=z$ and $u=\eta$)  and the definition of the class of symbols $S_{\widetilde l,m,N}(T^\ast \R^{n-1}; {\mathcal B}; \lambda)$, we estimate the reminder term by 
$$
\begin{array}{rcl}
\displaystyle\sup_\lambda \|F_N\|_{C^n} 
&\le& C \sum_{r+s\le N-1} \|q_r\|_{C^{2n}}\,  
\|a_s\|_{C^{2n}}\\[0.3cm]
\displaystyle &\le& C \sum_{r\le N-1} (\|q_r\|^2_{C^{\tilde l-mr}}
+ \|a_r\|^2_{C^{\tilde l-mr}})  \le\ C_{\mathcal B}. 
\end{array}
$$

In the same way,  we write 
 $A(\lambda)Q'(\lambda)$ modulo $O_{\cal B}(|\lambda|^{-N})$ as a $\lambda$-FIO with distribution
 kernel  (\ref{FIO}) with amplitude 
given by the oscillatory integral 
\[
K_{AQ'}(x,\xi,\lambda) = \sum_{j=0}^{N-1}\sum_{s+r=j}\lambda^{-j}\left(\frac{\lambda}{2\pi}\right)^{n-1}\,   
 \int_{{\R}^{2n-2}} 
\, e^{-i\lambda \langle z,\eta\rangle} \, 
a_s(x,\eta+\xi)q_r'(z+x +\psi_2(x,\xi,\eta),\xi)
 d\eta dz \, ,
\] 
where $ \psi_2(x,\xi,\eta) = \int_0^1 \nabla_\xi \psi(x, \xi +\tau \eta) d \tau$. 
We get as above
\[
K_{AQ'}(x,\xi,\lambda) =
\sum_{j=0}^{N-1} H_{j}(x,\xi)\lambda^{-j} + H_{N}(x,\xi,\lambda)\lambda^{-N}
\] 
where 
$$
\|H_N\|_{C^n} 
\le C \sum_{r\le N-1}( \|a_r\|^2_{C^{\tilde l-mr}} +  
 \|q_r'\|^2_{C^{\tilde l - mr}}) \le\ C_{\mathcal B}
$$
and  
\begin{equation}
\begin{array}{lcrr}
\displaystyle H_{j}(x,\xi)\ =\ 
\displaystyle \sum_{r+s +|\beta|= j}\, 
\frac{1}{\beta!}\, \left[  \left(-i\partial_\eta\right)^\beta  \left(
a_s(x, \eta +\xi )
\partial^\beta_z\, q_r'(z +x + \psi_2(x,\xi,\eta),\xi)\, \right)\right]_{|\eta =0,z=0}\,  
                    \label{symbolAQ} 
\end{array}
\end{equation}
for $0\le j\le N-1$. Hence,  $H_j\in C_0^{\tilde l -mj}(T^\ast \R^{n-1})$, $0\le j\le N-1$,  and it depends continuously on $K$. 
Note that $\psi_1(x,0,\xi)=\nabla_x\psi(x,\xi)$, $\psi_2(x,\xi,0)=\nabla_\xi\psi(x,\xi)$, and
that locally 
$${\rm graph}\, \kappa = 
\{(x,\xi +\nabla_x\psi(x,\xi), x + \nabla_\xi\psi(x,\xi),\xi)\}.
$$
In particular, this relation implies that $\sigma(Q') = \sigma(Q)\circ \kappa$. 
 Since $A(\lambda)$ is elliptic 
on $\mbox{WF}'(Q)$ we obtain that $a_0(x,\xi)\neq 0$ on the support of the functions
$(x,\xi)\to q_r(x,\xi + \nabla_x\psi(x,\xi))$, and we 
 determine $q_j'$ by recurrence from the equations 
\[
H_{j}(x,\xi) \ = \ F_{j}(x,\xi)\, ,\ j=0,\ldots, N-1. 
\]
It is easy to see by recurrence 
that $q_j'\in C_0^{\tilde l - mj}(T^\ast R^{n-1})$ and that it depends  continuously with respect 
to $K\in C^\ell(\Gamma)$. If $A(\lambda)$ is elliptic on $\mbox{WF}'(Q')$, we obtain $Q(\lambda)$ in the same way.
\finishproof 
\begin{Remark}\label{Rem:commutator} We have proved in particular that if $Q(\lambda)$, $\lambda\in {\mathcal D}$,  is a  $\lambda$-PDO of order $0$ in $\R ^{n-1}$
defined  by (\ref{pdo}) with symbol 
 $q\in S_{\widetilde l,m,N}(T^\ast \R^{n-1}; {\mathcal B}; \lambda)$ and  if the
 distribution kernel of $A(\lambda)$  is  given by (\ref{FIO}) with amplitude 
 $a\in S_{\widetilde l,m,N}(T^\ast \R^{n-1}; {\mathcal B}; \lambda)$, 
 then  $Q(\lambda)A(\lambda)$ and $A(\lambda)Q(\lambda)$ are
 $\lambda$-FIOs in $\R^{n-1}$ with distribution kernels (\ref{FIO}) and amplitudes 
 in $ S_{\widetilde l,m,N}(T^\ast \R^{n-1}; {\mathcal B}; \lambda)$. Their principal symbols are given by the product of the corresponding principal symbols. 
By the same argument, the
class of zero order $\lambda$-PDOs  with symbols in $S_{\widetilde l,m,N}(T^\ast \R^{n-1}; {\mathcal B}; \lambda)$ is closed under multiplication and 
transposition   and  it 
does not depend on the choice of the local coordinates (modulo $O_{\mathcal B}(|\lambda|^{-N})$). 
\end{Remark}

\subsection{Proof of Proposition \ref{prop:commutator}.}
\label{subsec:commutator}

Given $f\in C^N(\T^{n-1}\times D)$ we denote by $T_N f$ its Taylor polynomial
with respect to $I$ at $I=I^0$  
 \[
 T_Nf(\varphi, I) = 
\sum_{0\le |\alpha|\le N}(I-I^0)^\alpha \, f_\alpha(\varphi)\, , 
 \] 
 where $f_\alpha(\varphi)=\partial^\alpha_I f (\varphi, I^0)/\alpha !$ are the corresponding 
 Taylor coefficients. Recall that $W_1(\lambda)$ is a $\lambda$-FIO operator of the form \eqref{operator} with symbol 
 $\sigma(W_1)(x,I,\lambda) = w_0(x,I) + \lambda^{-1}w^0 (x,I,\lambda)$, where $w_0=1$ in 
 $\T^{n-1}\times D^0$ and $w^0$ satisfies \eqref{w-symbol} with given  $\ell$ and $M$ such that $\ell\ge 2M + 2n$. 
 In this notations, $w_{j,\alpha}^0(\varphi)= \partial^\alpha_I w_{j}^0(\varphi, I^0)/\alpha !$ and  \eqref{w-symbol} implies
 \begin{equation}
 \label{w-symbol1}
 \left\{
 \begin{array}{lcrr}
 w_{j,\alpha}^0\in C^{\ell -2j-|\alpha|}(\T^{n-1})\ \mbox{for }\ 0\le j\le M-1 \ \mbox{and}\ |\alpha|\le \ell - 2j\, ,\ \mbox{and} 
\\[0.3cm]
 \mbox{the map}\ K \to w_{j,\alpha}^0\in C^{\ell -2j-|\alpha|}(\T^{n-1})\ \mbox{is continuous in}\  K\in {\mathcal B}\subset C^\ell(\Gamma,\R). 
 \end{array}
 \right.
 \end{equation} 
%Fix $l\ge M(\tau +2) + 2n$ and $\ell> l + (n-1)/2$. 
We need  the following 
\begin{Lemma}\label{lemma:commutator2}
Let  $A(\lambda)$ be a $\lambda$-PDO on $\T^{n-1}$ with a symbol $\sigma(A)(\varphi,I, \lambda) = a_0(I) + \lambda^{-1}a^0(\varphi,I,\lambda)$, and let 
 $W^0(\lambda)$ be a $\lambda$-FIO operator of the form \eqref{operator} with a  symbol 
$\sigma(W^0)(\varphi,I, \lambda)=p_0(I) + \lambda^{-1} p^{0}(I,\lambda)$. Suppose that 
$a_0(I)=p_0(I)=1$  in a neighborhood $D^0$ of $I^0$, 
  and that the symbols 
\[
a^0(\varphi,I,\lambda)=\ \psi(I)\,  \sum_{j+|\alpha| \le M-1}\, 
\lambda^{-j}(I-I^0)^\alpha a^0_{j,\alpha}(\varphi)
\]
 and $p^0(I,\lambda) =   \sum_{j=0}^{M-1} 
\lambda^{-j}  p^0_{j}(I)$ satisfy (2) and (3) in Proposition \ref{prop:commutator} with $l \ge M(\tau +2) +2n$ and $\ell>l+(n-1)/2$.  Set
\[
 R(\lambda):= W_1(\lambda)A(\lambda) - A(\lambda)W^0(\lambda) . 
\]
Then
\[
R(\lambda) = \lambda^{-1}R_1(\lambda) +  R^0(\lambda) +  O_{\mathcal B}(|\lambda|^{-M-1}) \, ,
\]
 where 
$R_1(\lambda)$ and  $R^0(\lambda)$  
are $\lambda$-FIOs  of order $0$ of the form (\ref{operator}) and such that 
\begin{itemize}
\item[(1)]
the  symbol of $R^0(\lambda)$   satisfies (\ref{reminders}), 
\item[(2)]
the symbol of $R_1(\lambda)$ has the form
\[
R_1(\varphi, I, \lambda)=\sum_{j=0}^{M-1}R_{1j}(\varphi, I)\lambda^{-j},
\]
where 
\item[(3)]
for any $0\le j\le M-1$ we have 
\begin{equation}
R_{1j}(\varphi, I) = 
a_{j}^0(\varphi-2\pi\omega,I)-a_{j}^0(\varphi,I)\, +\, T_{M-j-1}w_{j}^0(\varphi,I)\, -\, p^0_{j}(I)\,
+\, h^0_{j}(\varphi,I)\, , 
\label{reminder}
\end{equation}
where
$h_j^0 =  f_j^0 - g_j^0$, $0\le j\le M-1$,  are  such that
\begin{itemize}
\item[(i)]
 the Taylor coefficients $f_{j,\alpha}^0(\varphi):= \partial_I^\alpha f^0_j(\varphi, I^0)/\alpha !$,  $|\alpha|\le M-j-1$,  are linear combinations, with coefficients independent of $K$,  of terms of the form 
%\begin{equation}
%\label{taylor1}
%\left\{
%\begin{array}{lll}
\begin{itemize}
\item[(a)]
$\partial_\varphi^{\beta} a_{j,\gamma}^0(\varphi-2\pi\omega)$, where $|\beta|+|\gamma|\le|\alpha|$  and $|\beta|\ge 1$, 
\item[(b)]
$\partial_\varphi^{\beta} a_{s,\gamma}^0(\varphi-2\pi\omega)$, where $s \le j-1$, $|\gamma|\le |\alpha|$,   and $|\beta|+|\gamma|\le 2(j-s)+|\alpha|, $  
\item[(c)]
 $w_{r,\delta}^0(\varphi)\partial_\varphi^{\beta} a_{s,\gamma}^0(\varphi-2\pi\omega)$, where
$r+s\le j-1$, $|\gamma|\le |\alpha|$ and $|\beta| + |\delta| + |\gamma| \le 2(j-1-r-s) + |\alpha|, $ 
\end{itemize}
%\end{array}
%\right.
%\end{equation}
\item[(ii)]  $g^0_{0}= 0$ and 
 the Taylor coefficients 
 \[
 g^0_{j,\alpha}(\varphi):= \partial_I^\alpha g_j(\varphi, I^0)/\alpha !\, ,\ 1\le j\le N-1\, ,\ |\alpha|\le M-j-1\, ,
 \]
 are linear combinations, with coefficients independent of $K$,   of functions 
\begin{equation}
p^0_{k,\beta}\, a^0_{j-k-1,\gamma}(\varphi)\, ,\ \mbox{where}  \ 0\le k\le j-1\,  ,\  \mbox{and}  \ \beta +\gamma =\alpha\, .
\label{taylor2}
\end{equation}
\end{itemize}
\end{itemize}
\end{Lemma}
\begin{Remark}\label{Rem:BNF}
$h^0_{j}$ does not depend on  $a^0_{r}$ for $r > j$. Moreover, $h^0_{j,\alpha}$ does not depend on  $a^0_{j,\gamma}$
for $|\gamma|\ge |\alpha|$. 
\end{Remark}
The proof of the lemma is given below. 
We proceed with the proof of Proposition \ref{prop:commutator}.
 Our goal is first to solve the system of equations $R_{1j}=0$, $0\le j\le N-1$,  with respect to $p^0_{j,\alpha}$ and $a^0_{j,\alpha}$ in the corresponding classes of functions and then to apply Lemma \ref{lemma:commutator2} in order to prove that  the  reminder term $R^0(\lambda)$ satisfies \eqref{reminders}. More precisely, we are going to find  by recurrence  
\[
p^0_{j,\alpha}\in \C\  \mbox{and}\  a^0_{j,\alpha}\in {\mathcal  A}^{l-\tau-j(\tau+2)-\tau|\alpha|}(\T^{n-1})\, ,\ 
 0\le j \le M-1\, ,\ |\alpha| \le  M-j-1 \, , 
 \] 
so that $R_{1j}=0$, where $R_{1j}$ is given by \eqref{reminder}. Moreover, we shall 
prove    that 
 the maps 
\begin{equation}
K \mapsto p^0_{j,\alpha}\in \C \, ,\ 
K\mapsto a^0_{j,\alpha}\in {\mathcal  A}^{l-\tau-j(\tau+2)-\tau|\alpha|}(\T^{n-1}) 
                                                            \label{maps2}
\end{equation}
are continuous with respect to
$K\in C^\ell(\Gamma,\R)$. In order to obtain uniqueness of the solutions we normalize $a^0_{j,\alpha}$ by
$\int_{\T^{n-1}}  a^0_{j,\alpha}(\varphi)\, d \varphi = 0$. 

For $j=0$ we have $h_0^0 = f_0^0$,   where $f_{0,\alpha}^0(\varphi)$ is a linear combination of 
$\partial_\varphi^{\beta} a_{0,\gamma}^0(\varphi-2\pi\omega)$, where $|\beta|+|\gamma|\le|\alpha|$ and  $|\gamma|<|\alpha|, $
and  we put
\[
p_{0,\alpha}^0\ =\ \frac{1}{(2\pi)^{n-1}}\int_{\T^{n-1}}  
w^0_{0,\alpha}(\varphi) \, d\varphi\, , \quad |\alpha|\le M-1\, .
\]
Setting $u_\alpha:=a_{0,\alpha}$ and $f_\alpha:=
p_{0,\alpha}^0 -w^0_{j,\alpha}-f_{0,\alpha}^0$  
we obtain from (\ref{reminder})  the  homological 
equations (\ref{homological}) for any $\alpha$ with $|\alpha|\le N-1$, which we solve by recurrence with respect to $|\alpha|$. By Remark \ref{Rem:BNF}    $h_{0,\alpha}^0=f_{0,\alpha}^0$ does not depend on $a_{0,\gamma}$, $|\gamma|\ge |\alpha|$, hence,  the  homological equation has a unique solution 
$a_{0,\alpha}$, $|\alpha|\le N-1$,  normalized by $\int_{\T^{n-1}}  a_{0,\alpha}(\varphi)\, d \varphi = 0$. 

Recall from \eqref{w-symbol1} and \eqref{bernstein} that for  $2j + |\alpha|\le l$ the map 
\begin{equation}
 C^\ell(\Gamma,\R)\ni K \ \to\  w^0_{j,\alpha} \in C^{\ell-2j-|\alpha|}(\T^{n-1}) \hookrightarrow  {\mathcal A} ^{l-2j-|\alpha|}(\T^{n-1})
\label{maps1}
\end{equation}
is continuous since $\ell>l+(n-1)/2$.
%\begin{equation}
%\frac{1}{i}\, {\cal L}_\omega \, 
%u(\varphi)\, = \, v(\varphi)\, , \quad \int_{\T^{n-1}}  v(\varphi)\, d \varphi = 0\, .
%\label{homological}
%\end{equation}
Now using Lemma \ref{lemma:homological} and (\ref{maps1}) for $j=0$ and $|\alpha|\le M-1$  
we obtain that $p_{0,\alpha}^0\in \C$ and 
$a_{0,\alpha}\in {\mathcal  A}^{l-\tau-\tau|\alpha|}(\T^{n-1})  $  and we prove that the corresponding 
maps (\ref{maps2}) are continuous. Moreover,
\[
p^0_0(I^0) = \frac{1}{(2\pi)^{n-1}}\, \int_{\T^{n-1}} \, w^0_0(\varphi,I^0) d\varphi\, .
\]

Fix $1\le j \le M-1$ and suppose that 
the inductive assumption holds for all indices $k\le j-1$. Moreover, take $0\le \alpha_0\le M-j-1$ and suppose as well that the inductive assumption holds for any pair $(j,\alpha)$ with $0\le |\alpha|< \alpha_0$ if $\alpha_0>0$. We shall prove that the inductive 	assumption holds for any pair $(j,\alpha)$ with $|\alpha|= \alpha_0$. 
Firstly, using   Lemma \ref{lemma:commutator2} (i) and  (ii) we are going to show that
the maps
\[
K\mapsto h_{j,\alpha}\in {\mathcal  A}^{l-j(\tau + 2) -\tau|\alpha|}(\T^{n-1})\, ,\ 
|\alpha| = \alpha_0 \le M-j-1\, ,
\]
are continuous with respect to
$K\in C^\ell(\Gamma,\R)$. 

Let $|\alpha|= \alpha_0$. 
If $\alpha_0>0$
the function $\varphi \to \partial_\varphi^{\beta} a_{j,\gamma}^0(\varphi-2\pi\omega)$  in Lemma \ref{lemma:commutator2} (3)-(i)-(a), where $|\beta|+|\gamma|\le |\alpha|=\alpha_0$  and $|\beta|\ge 1$
(then $|\gamma|< |\alpha|=\alpha_0$),   belongs to ${\mathcal  A}^{p}(\T^{n-1})$, where 
$$
\begin{array}{rcll}
p&:=&l-j(\tau + 2) -\tau|\gamma| -|\beta| > l-j(\tau + 2) -\tau(|\beta|+|\gamma|)\\[0.3cm]
&\ge& l-j(\tau + 2) -\tau|\alpha|.
\end{array}
$$

The term $\partial_\varphi^{\beta} a_{s,\gamma}^0$ in (b) belongs to ${\mathcal  A}^p(\T^{n-1})$, where
$$
\begin{array}{rcll}
p&:=&l - \tau-s(\tau+2)-\tau|\gamma|-|\beta |\\[0.3cm]
 &=& l - \tau -s(\tau+2)  -(\tau-1)|\gamma| -(|\beta|+|\gamma|)\\[0.3cm]
&\ge& l - \tau -s(\tau+2)  -(\tau-1)|\alpha| -(2(j-s) +|\alpha|)\\[0.3cm]
&=& l  -(s+1)\tau -2j -\tau|\alpha|
\ge l - j(\tau +2) -\tau|\alpha|. 
\end{array}
$$

Consider now the terms in (c). By \eqref{maps1} we have  
\[
w_{r,\delta}^0\in {\mathcal  A}^{l-2r -|\delta|}(\T^{n-1}) \subset {\mathcal  A}^{l-j(\tau + 2) -\tau|\alpha|}(\T^{n-1})
\]
and   it depends continuously on $K\in C^\ell(\Gamma,\R)$ in these spaces, since $2r +|\delta|\le 2(j-1)+ |\alpha|$.
 Moreover, $\partial_\varphi^\beta a_{s,\gamma}(\varphi-2\pi \omega)$ belongs to ${\mathcal A}^{p}(\T^{n-1})$, where
$$
\begin{array}{lcrr}
p:=l-\tau -s(\tau+2)-\tau|\gamma| - |\beta|= l-(s+1)\tau -(\tau-1)|\gamma|-(2s+|\beta+\gamma|)\\[0.3cm]
\ge  l-j(\tau+2)-\tau|\alpha|\, , 
\end{array}
$$
hence, $\partial_\varphi^{\beta} a_{s,\gamma}^0(\varphi-2\pi\omega)$ belongs to 
${\mathcal A}^{l-j(\tau+2)-\tau|\alpha|}(\T^{n-1})$ and  it depends continuously on $K\in C^\ell(\Gamma,\R)$. 
We have proved that $f^0_{j,\alpha}\in {\mathcal A}^{l-j(\tau+2)-\tau|\alpha|}(\T^{n-1})$ for any multi-index $\alpha$ with length $|\alpha|= \alpha_0$ and that it depends continuously on $K\in C^\ell(\Gamma)$. In the same way, using (\ref{taylor2}), we obtain that  $g^0_{j,\alpha}\in {\mathcal A}^{l-j(\tau+2)-\tau|\alpha|}(\T^{n-1})$ and that it depends continuously on $K\in C^\ell(\Gamma,\R)$.
Moreover,   
$$w^0_{j,\alpha}\in {\mathcal A}^{l-j(\tau+2)-\tau|\alpha|}(\T^{n-1})$$ 
and it  depends continuously on $K\in C^\ell(\Gamma,\R)$.

We set as above 
\[
p_{j,\alpha}^0\ =\ \frac{1}{(2\pi)^{n-1}}\int_{\T^{n-1}}  
(w^0_{j,\alpha}(\varphi) -h_{j,\alpha}(\varphi))\, d
\varphi\, .
\] 
Obviously it depends  continuously on   $K\in C^\ell(\Gamma,\R)$. 
Setting $u=a_{j,\alpha}$ and $f=
p_{j,\alpha}^0 -w^0_{j,\alpha} + h_{j,\alpha}\in {\mathcal  A}^{l-j(\tau + 2) -\tau|\alpha|}(\T^{n-1})$, $|\alpha| = \alpha_0 \le M-j-1$, 
we solve (\ref{homological}) by Lemma \ref{lemma:homological},  first for $\alpha_0 =0$, then for $\alpha_0=1$, and so on,  and we prove as above that the maps (\ref{maps2}) 
are continuous. 
In this way we obtain symbols $p ^0$ and $a^0$ satisfying (2) and (3) in Proposition \ref{prop:commutator} with $l \ge M(\tau +2) +2n$ and $\ell>l+(n-1)/2$ and such that
$R_{1j}=0$ for $1\le j\le M-1$. 
Now Lemma \ref{lemma:commutator2} implies that $R(\lambda) = R^0(\lambda) + O_{\mathcal B}(|\lambda|^{-M-1})$, where
$R^0(\lambda)$ satisfies (\ref{reminders}). This completes the proof of Proposition \ref{prop:commutator}.
 \finishproof

\noindent
{\em Proof of Lemma \ref{lemma:commutator2}.} The proof of the Lemma is similar to that of Lemma \ref{Lemma:commutator}. 
First we write the  Schwartz kernel the  operator $\widetilde{W_1(\lambda)A(\lambda)}$ as an oscillatory  integral 
\[
\left(\frac{\lambda}{2\pi}\right)^{n-1}\, \int_{{\R}^{n-1}} 
\, e^{i\lambda (\langle x-y,I\rangle + \Phi(x,I))} \, 
F(x,I,\lambda)\, d x d I \, ,
\]
where the  amplitude $F$ is given by the oscillatory integral
\[
F(x,I,\lambda) = 
\left(\frac{\lambda}{2\pi}\right)^{n-1}\, \int_{{\R}^{2n-2}} 
\, e^{i\lambda (\langle x-z,\xi-I\rangle + \Phi(x,\xi) -\Phi(x,I))} \, 
w(x,\xi,\lambda)
a(z,I,\lambda)J(x,\xi)\, d\xi dz 
\] 
modulo $O_{\mathcal B}(|\lambda|^{-M-1})$.  Recall that the function $\Phi(x,I)=L(I) + \Phi^0(x,I)$ is $C^\infty$-smooth,  and 
\begin{equation}
\label{phase}
\nabla L(I^0) = 2\pi\omega\ \mbox{and}\ \Phi^0(x,I)= O_N(|I-I^0|^N)\ \mbox{for any $N\in\N$}. 
\end{equation}
Set
\[ 
\begin{array}{lcrr} 
\Psi(x,I,\eta) =  \int_0^1 \nabla_I \Phi (x,I +\tau \eta) d \tau= L_0(I,\eta) + H_0(x,I,\eta),\ \mbox{where}  \\[0.3cm]
L_0(I,\eta)=\int_0^1 \nabla_I L(I +\tau \eta) d \tau\ \mbox{and}\
H_0(x,I,\eta) = \int_0^1 \nabla_I \Phi^0(x,I +\tau \eta) d \tau\, .
\end{array}
\]
Changing the variables $z= v+x+ \Psi(x,I,\eta)$, $\xi=I +\eta$,  and using (\ref{newsymbols}) and Lemma \ref{Lemma:estimate1} we obtain 
\begin{equation}
F(x,I,\lambda)\, =\, \tilde F_0(x,I,\lambda) + \sum_{j=0}^{M-1}\lambda^{-j-1} \tilde F_{j}^0(x,I,\lambda) +
 \widetilde F(x,I,\lambda)\, ,
\label{amplitude}
\end{equation}
modulo $O_{\mathcal B}(|\lambda|^{-M-1})$, where $\widetilde F$ belongs to the residual class $\widetilde R_{M+1}(\T^{n-1}\times D;{\mathcal B};\lambda)$, and
\[
\tilde F_0(x,I,\lambda) =\left(\frac{\lambda}{2\pi}\right)^{n-1}\, \int_{{\R}^{2n-2}} 
\, e^{-i\lambda \langle v,\eta\rangle } \, c_0(x,I,v,\eta)\, d\eta dv \, ,
\]
\[
\tilde F_{j}^0(x,I,\lambda) =\left(\frac{\lambda}{2\pi}\right)^{n-1}\, \int_{{\R}^{2n-2}} 
\, e^{-i\lambda \langle v,\eta\rangle } c_j^0(x,I,v,\eta)\, d\eta dv\, .
\]
Moreover, 
\[
\begin{array}{rcll} 
c_0(x,I,v,\eta) &=&  w_0(x,I + \eta)J(x,I + \eta)a_0(I), \\[0.3cm]
 c_j^0(x,I,v,\eta) &=& 
w_j^0(x,I + \eta)J(x,I + \eta)a_0(I) \\[0.3cm]
&+& \displaystyle \sum_{|\gamma| \le M-j-1} \,c_{j,\gamma}^0(x,I,v,\eta)J(x,I + \eta)  \,  (I-I^0)^\gamma\,
\end{array}
\]
and
\[
\begin{array}{rcll} 
\displaystyle
c_{j,\gamma}^0(x,I,v,\eta)   &=&  \psi(I)w_0(x,I + \eta)\,  
a^0_{j,\gamma}(v +x  + \Psi(x,I,\eta) )  \\  [0.3cm]
 &+&  \displaystyle \sum_{r + s = j-1}  \psi(I)\,w_r^0(x,I + \eta)
 a^0_{s,\gamma}(v +x  + \Psi(x,I,\eta)) \, ,
\end{array}
\]
where $\psi\in C_0^\infty(\R^{n-1})$ and $\psi=1$ in a neighborhood $D^0$ of $I=I^0$. 
Hereafter, to simplify the notations we denote the class of $x\in\R^{n-1}$ in $\T^{n-1}$ by $x$ as well. 
Recall that the function  $J(x,\xi)-1$ vanishes up to any order at $\xi=I^0$, hence, 
\[
J(x,\xi) = 1 + O(|\xi-I^0|^{M+1}) , 
\] 
and as in the proof of Lemma \ref{Lemma:commutator} we obtain  
\[
\tilde F_0(x,I,\lambda)=w_0(x,I)a_0(I)J(x,I) = 1 + O(|I-I^0|^{M+1}) \, .
\]
In the same way we write
\[
\tilde F_{j}^0(x,I,\lambda) =w_j^0(x,I)
 + \sum_{|\gamma| \le M-j-1} \,F_{j,\gamma}^0(x,I,\lambda)  \,  (I-I^0)^\gamma\, +  O(|I-I^0|^{M+1})\, ,
\]
where the amplitude of $F^0_{j,\gamma}(x,I,\lambda)$ is  $c_{j,\gamma}^0(x,I,v,\eta)J(x,I + \eta)$. 
Consider the composition
\begin{equation}
\label{map}
\T^{n-1}\times D_1\times T^{n-1}\times D_2\ni(x,I,v,\eta) \to a^0_{j,\gamma}(v+x+\Psi(x,I,\eta))\, , 
\end{equation}
where $D_1$ and $D_2$ are suitable neighborhoods  of $I=I^0$  and $\eta=0$  respectively such that $\Psi\in C^\infty(\T^{n-1}\times D_1\times D_2)$. 
One should be  careful when changing variables on the torus since the class ${\mathcal A}^{q}(\T^{n-1})$ is not conserved except in the case when the corresponding diffeomorphism is affine linear. 
By assumption $a^0$ satisfies \eqref{symbols}, where $l\ge M(\tau+2)+2n$. Then the second relation of \eqref{bernstein} implies 
\[
a^0_{j,\gamma}\in {\mathcal A}^{q}(\T^{n-1})\subset C^q(\T^{n-1})\ \mbox{for}\  j +|\gamma|\le M-1,
\]
and $a^0_{j,\gamma}$  depends continuously on  $K\in C^\ell(\Gamma,\R)$ in these spaces, where 
\begin{equation}
\label{q}
q:=l-\tau -j(\tau +2) -\tau|\gamma| \ge M(\tau +2) - (j + |\gamma| +1)\tau - 2j + 2n 
\ge 2(M-j)+ 2n. 
\end{equation}
In particular,  the map \eqref{map} belongs to $C^q(\T^{n-1}\times D_1\times T^{n-1}\times D_2)$ and it depends continuously on  $K\in C^\ell(\Gamma,\R)$. 

We proceed as in the proof of Lemma \ref{Lemma:commutator}. We  develop $a^0_{j,\gamma}(v+x+\Psi)$ in Taylor series with respect to $v$ 
at $v=0$  up to order $O(|v|^{M-j-|\gamma|})$ using the 
Taylor formula with an integral reminder term. 
For any $|\beta|\le M-j$ and $j+|\gamma|\le M-1$
we have  
$\partial_x^\beta a^0_{j,\gamma}\in {\mathcal A}^p(\T^{n-1})\subset C^p(\T^{n-1})$,  
where by \eqref{q} we get
\begin{equation}
p = l-\tau -j(\tau +2) -\tau|\gamma| -|\beta|   \ge 
2(M-j) -|\beta| +2n \ge  |\beta| + 2n\, .
\label{order1}
\end{equation}
Moreover,  $\partial_x^\beta a^0_{j,\gamma}\in {\mathcal A}^{|\beta|+ 2n}(\T^{n-1}) \subset C^{|\beta|+ 2n}(\T^{n-1}) $ 
depends continuously on $K\in {\cal B}$  if $|\beta|\le M-j$ and  $j+|\gamma|\le M-1$. 
Integrating  $\beta$ times by parts with respect to $\eta$ in the corresponding oscillatory integrals,    we get  in view of (\ref{order1}) and Lemma \ref{Lemma:estimate1} that
\[
\begin{array}{lcrr}
\displaystyle F^0_{j,\gamma}(x, I,\lambda) 
 \ =\  
\sum_{ |\beta| \le  M-j-|\gamma|}\, 
\frac{\lambda^{-|\beta|}}{\beta!}\, \left[ (-i\partial_\eta)^\beta    (\partial^\beta_x\,
a^0_{j,\gamma})(x - \Psi(x,I,\eta) )\right]_{|\eta =0}  \\ [0.7cm]
\displaystyle + \sum_{r+s = j-1}\, \sum_{|\beta|\le M-j-1-|\gamma|}\, 
\frac{\lambda^{-|\beta|}}{\beta!}\, \left[ (-i\partial_\eta)^\beta  \left(w_r^0(x,I + \eta)\, 
(\partial^\beta_x\,
a_{s,\gamma}^0)(x - \Psi(x, I,\eta) )\right)\right]_{|\eta =0} \, \\ [0.7cm]
+\,  \lambda^{-M+j+|\gamma|}\widetilde F_{j,\gamma}(x,I,\lambda)\, ,
\end{array}
\]
where $\sup_\lambda\|\widetilde F_{j,\gamma}\|_{C^n}\le C_{\cal B}$, hence $ \lambda^{-M+j+|\gamma|}\widetilde F_{j,\gamma}$ belongs to the residual class. 
It follows from \eqref{phase}  that  
all the derivatives of $H_0$ 
vanish at  $(\eta, I)=(0,I^0)$, and we have $\partial_\eta^\gamma H_0(x,I,0) = O(|I-I^0|^{M+1})$
for any $\gamma$, hence, 
$$
F(x,I,\lambda) = 
F_0(x,I) +   \sum_{j=0}^{M-1} F^0_{j}(x,I)\lambda^{-j-1} +  F^1(x,I,\lambda) \, ,
$$ 
where $F_0=1$ in $\T^{n-1}\times D^0$,  
\begin{equation}
\label{F^0}
F^0_{j}(x,I)\ =\ a^0_j(x - 2\pi\omega ,I)  + w^0_j(x,I) + f^0_j(x,I)\,  , \ 0\le j\le M-1\, ,
\end{equation}
$f^0_{0}= a^0_0(x - \nabla L(I) ,I) - a^0_0(x - 2\pi\omega ,I)$, and  for $j\ge 1$ we have
\begin{equation}
\begin{array}{lcrr}
\displaystyle f^0_{j}(x,I)\ =\  a^0_j(x - \nabla L(I) ,I) - a^0_j(x - 2\pi\omega ,I)\\ [0.5cm]
\displaystyle +\, \sum_{s=0}^{j-1}\, 
\sum_{ |\beta| = j-s}\, \sum_{|\gamma|\le M-j-1}\, 
\frac{(I-I^0)^\gamma}{\beta!}\, \left[ (-i\partial_\eta)^\beta   \partial^\beta_x\,
a^0_{s,\gamma}(x - L_0(I,\eta) )\right]_{|\eta =0}  \\ [0.5cm]
\displaystyle +  \sum_{r+s +|\beta|= j-1}\, \sum_{|\gamma|\le M-j-1}\, 
\frac{(I-I^0)^\gamma}{\beta!}\, \left[ (-i\partial_\eta)^\beta  \left(w_r^0(x,I + \eta)\, 
\partial^\beta_x\,
a_{s,\gamma}^0(x - L_0(I,\eta) )\right)\right]_{|\eta =0} .
                    \label{f} 
\end{array}
\end{equation}
We have also $F^1\in \widetilde R_{M+1}(\T^{n-1}\times D;{\mathcal B}, \lambda)$ 
in view of (\ref{order1}) and Lemma \ref{Lemma:estimate1}. 
Expanding  the right hand side of (\ref{f}) in Taylor series with respect to $I$ at $I^0$ up to order $M-j-1$, we obtain 
 $$
 f^0_{j}(x,I) = \sum_{|\alpha|\le M-j-1}f^0_{j,\alpha}(x) (I-I^0)^\alpha + 
\sum_{|\alpha|= M-j}\tilde f^0_{j,\alpha}(x,I)(I-I^0)^\alpha,
$$ 
where $\tilde f^0_{j,\alpha}\in C^{2n}(\T^{n-1}\times D)$ is bounded in $K\in {\mathcal B}$ and it contributes to the residual term.
We write $f^0_{j,\alpha}= f^{0,1}_{j,\alpha} + f^{0,2}_{j,\alpha}+f^{0,3}_{j,\alpha}$, where $f^{0,1}_{j,\alpha}$ comes from the first line  in the right hand side of  (\ref{f}), $f^{0,2}_{j,\alpha}$ comes from the second one, and  $f^{0,3}_{j,\alpha}$ from the third line.

\noindent
{\em  Case 1.} 
Since $\nabla L(I^0) = 2\pi \omega$, 
expanding $a^0_j(x - \nabla L(I) ,I)$ in Taylor series with respect to $I$ at $I^0$ up to order $M-j$, we get
\[
a^0_j(x - \nabla L(I) ,I) = a^0_j(x - 2\pi \omega ,I) + b^0_j(x - 2\pi \omega   ,I)
 + c_j(x, I), 
\]
where $b^0_j(x -  2\pi \omega ,I)$ is a linear combination of terms 
\[
\partial_x^{\beta} a_{j,\gamma}^0(x-2\pi\omega)(I-I^0)^\alpha,  \ \mbox{where} \ j+|\alpha|\le N-1, \ |\beta|+|\gamma|\le |\alpha|\ \mbox{and}\ |\beta|\ge 1, 
\]
while 
$c_j(x,I)= \sum_{|\alpha|=N-j}c_{j,\alpha}(x,I)(I-I^0)^\alpha$ and $c_{j,\alpha}\in C^{2n}(\T^{n-1}\times D)$ is bounded in $K\in {\mathcal B}$. Hence,  $\partial_x^{\beta} a_{j,\gamma}^0(x-2\pi\omega)$  is in the set described by (a) and $\lambda^{-j-1}c_j$ is a residual symbol.   \\

\noindent
{\em  Case 2.} 
Consider the term $f^{0,2}_{j,\alpha}$, $j+|\alpha|\le M-1$. The derivative 
\[
\partial_\eta^\beta  
\partial^\beta_x\,
a_{s,\gamma}^0(x - L_0(I,\eta))|_{\eta=0}
\]
is a linear combination of terms
$\partial^{\beta+\beta'}_x\,
a_{s,\gamma}^0(x - L(I))\, L_{\beta'}(I,0)$, 
where $|\beta'|\le |\beta|$ and $L_{\beta'}(I,\eta)$ is polynomial of derivatives of $L_0$. Expanding the function 
\[
\partial^{\beta+\beta'}_x\,
a_{s,\gamma}^0(x - L(I))\, L_{\beta'}(I,0)(I-I^0)^\gamma\, ,
\]
in Taylor series at $I=I^0$ up to order $O(|I-I^0|^{N-j})$ and taking the coefficients of $(I-I^0)^{\alpha}$ we obtain that 
 $f^{0,2}_{j,\alpha}$ is a linear combination of terms
\[
\partial_x^{\beta+\beta'+\delta} a_{s,\gamma}^0(x-2\pi\omega)\, ,   
\]
where $\ s\le j-1\, , s+|\beta|= j\, ,\beta'\le \beta\, ,\ |\gamma|+|\delta|\le |\alpha|$.   
Hence $f^{0,2}_{j,\alpha}$ can be written as a linear combination of 
\[
\partial_x^{\beta} a_{s,\gamma}^0(x-2\pi\omega)\, ,\ \mbox{where}\ s \le j-1\, ,\ |\gamma|\le |\alpha|\, ,\ \mbox{and}\ |\beta|+|\gamma|\le 2(j-s)+|\alpha|\, ,    
\]
and we obtain the terms in (b). \\

\noindent
{\em  Case 3.} 
Consider the term $f^{0,3}_{j,\alpha}$. The function 
\[
\partial_\eta^\beta  \left(w_r^0(x,I + \eta)\, 
\partial^\beta_x\,
a_{s,\gamma}^0(x - L_0(I,\eta) )\right)
\]
is a linear combination of terms
\[
(\partial_\eta^{\beta'} w_r^0)(x,I + \eta)\, 
(\partial^{\beta+\beta''}_x\,
a_{s,\gamma}^0)(x - L_0(I,\eta))\, L_{\beta',\beta''}(I,\eta)
\]
where $|\beta'+\beta''|\le |\beta|$ and $L_{\beta',\beta''}(I,\eta)$ are polynomials of derivatives of $L_0$. Taking $\eta=0$ and using the Taylor expansion at $I=I_0$, we obtain as above that 
  $f^{0,3}_{j,\alpha}$ is a linear combination of terms
\[
 w_{r,\delta'+\beta'}^0(x)\, \partial_x^{\beta+\beta''+\delta'' } a_{s,\gamma}^0(x-2\pi\omega)\, , 
\]
where 
\[
r+s+|\beta|=j-1\, ,\ |\beta'+\beta''|\le |\beta|\, ,\ |\delta'+\delta''+\gamma| \le |\alpha|\, . 
\]
Hence, $f^{0,3}_{j,\alpha}$ can be written as a linear combination of
\[
 w_{r,\delta}^0(x)\partial_x^{\beta} a_{s,\gamma}^0(x-2\pi\omega)\, ,
\]
where 
\[
r+s\le j-1\, ,\  |\gamma|\le |\alpha|\, ,\   |\beta| + |\delta| + |\gamma| \le 2(j-1-r-s) + |\alpha|.   
\]
This gives (c) and  completes the proof of  (i). \\

In the same way we write 
 $A(\lambda)W_0(\lambda)$ in 
the form (\ref{operator}) with amplitude $G(x,I,\lambda)$ 
given by the oscillatory integral 
\[
\left(\frac{\lambda}{2\pi}\right)^{n-1}\,   
(p_0 + \lambda^{-1}p^0)(I,\lambda)\, \int_{{\R}^{2n-2}} 
\, e^{i\lambda (\langle x-z,\xi-I\rangle + \Phi(z,I) -\Phi(x,I))} \, 
a(x,\xi,\lambda)
 d\xi dz \, .
\] 
Changing the variables we obtain 
$G = a(p_0 +\lambda^{-1}p^0) + \widetilde G$, where 
$\widetilde G$ is given by 
\[
\left(\frac{\lambda}{2\pi}\right)^{n-1}\,  (p_0 + \lambda^{-1}p^0)(I,\lambda)\,  
\int_{{\R}^{2n-2}} 
\, e^{-i\lambda \langle v,\eta\rangle } 
[a(x, \eta  + I + H_1(x,v,I), \lambda) - a(x, \eta + I, \lambda)]\, 
d\eta dv \ , 
\] 
and 
\[
H_1(x,v,I) = \int_0^1 \nabla_x \Phi^0(x+\tau v,I) d \tau = O_N(|I-I^0|^N)
\]
for any $N\in \N$. 
Using Taylor formula and  Lemma \ref{Lemma:estimate1} 
we obtain  that $\widetilde G$ is a residual symbol satisfying  (\ref{reminders}),  and we get 
$$
G(x,I,\lambda) = G_0(x,I) +   \sum_{j=0}^{M-1} G^0_{j}(x,I)\lambda^{-j-1} +  
G^1(x,I,\lambda)\, ,
 $$ 
where
$G^1\in \widetilde R_{M+1}(\T^{n-1}\times D;{\mathcal B}, \lambda)$. Moreover, 
$G_0=1$ in $\T^{n-1}\times D^0$, and 
\[
 G^0_{j}(x,I)  \ =\ a^0_j(x ,I) + p^0_{j}(I) + g^0_j(x,I)\, , 
\]
where, $g^0_0 = 0$ and 
\[
g^0_j(x,I) = 
\sum_{k=0}^{j-1}a^0_{j-k-1}(x, I) p^0_{k}(I)
\]
for $j\ge 1$. Now it is easy to see that $g^0_{j,\alpha}(x)$, $j+|\alpha|\le N-1$,  are linear combinations of terms given by (\ref{taylor2}). 
Finally, 
\[
R_1(x,I, \lambda)=\sum_{j=0}^{M-1}T_{M-j-1}(F^0_j - G^0_j)(x,I)\lambda^{-j} \, ,
\]
and $R_1(\lambda)$ is the corresponding $\lambda$-FIO.
\finishproof

\subsection*{A.4.  Proof of Proposition \ref{prop:spectral-decomposition}.}\label{subsec:spectral-decomposition}

We obtain as above 
\[
\begin{array}{lcrr}
\widetilde{W^0(\lambda)e_k}(x)\ =\  \widetilde{e_k}(x)\, 
e^{i\lambda \Phi(x,\xi_k)} \\ [0.3cm]
\displaystyle \times \, 
\left(\frac{\lambda}{2\pi}\right)^{n-1}\, \int_{{\R}^{2n-2}} 
\, e^{i\lambda\langle x-y +  
\Phi_0(x,\xi_k, \eta_k),\eta_k \rangle} \, 
(p_0 + \lambda^{-1}p^0)(I,\lambda)\, dI \,  dy \, ,
\end{array}
\]
where $\Phi_0(x,\xi,\eta)=\int_0^1\nabla_\xi\Phi(x,\xi+\tau \eta)d\tau$, 
$\xi_k =  (k + \vartheta_0/4)/\lambda$ and $\eta_k=I-(k + \vartheta_0/4)/\lambda$. 
Deforming the contour of integration we obtain
\[
\begin{array}{lcrr}
W^0(\lambda)e_k(\varphi)\ =\  e_k(x)\, 
e^{i\lambda \Phi(\varphi,(k + \vartheta_0/4)/\lambda)} \\ [0.3cm]
\displaystyle \times \, 
\left(\frac{\lambda}{2\pi}\right)^{n-1}\, \int_{{\R}^{2n-2}} 
\, e^{-i\lambda\langle u,v\rangle } \, 
(p_0 + \lambda^{-1}p^0)(v + (k + \vartheta_0/4)/\lambda,\lambda)\, du\,  dv \,  + \, 
O_{\cal B}(|\lambda|^{-M-1}) \,  ,
\end{array}
\]
which implies  (\ref{newoperator}).

To prove (\ref{newrest}) we write $\widetilde{R^0(\lambda)e_k(x)}$ as an oscillatory integral as above, and then we change the contour of integration with respect to $y$ 
by 
\[
y \to v = y-x-\Phi_0(x,(k + \vartheta_0/4)/\lambda, I - (k + \vartheta_0/4)/\lambda)\, .
\]
This implies 
\[
\begin{array}{lcrr}
R^0(\lambda)e_k(\varphi)\ =\  e_k(\varphi)\, 
e^{i\lambda \Phi(\varphi,(k + \vartheta_0/4)/\lambda)}  \\ [0.3cm]
\displaystyle \times \, \left(\frac{\lambda}{2\pi}\right)^{n-1}\,
 \int_{{\R}^{2n-2}} 
\, e^{-i\lambda\langle v, I -(k + \vartheta_0/4)/\lambda\rangle }\, 
(r_0 + \lambda^{-1}r^0)(\varphi, I,\lambda)\, dI \, dv \, 
\end{array}
\]
modulo $O_{\cal B}(|\lambda|^{-M-1})$. We have 
\[
r_0(\varphi, I) = \sum_{|\alpha|=M+1}\, (I-I^0)^\alpha r_{0,\alpha}(\varphi, I)\, ,
\]
where $r_{0,\alpha}\in C^\infty_0(\T^{n-1}\times D)$ does not depend on $K$, and we
 write now $r^0$ in the form (\ref{newreminders}).   
Integrating  by parts with respect to $v$  in the corresponding oscillatory integral 
with amplitude $r_{j,\alpha}^0(\varphi,I)(I-I^0)^\alpha$, $|\alpha| = M-j$, 
we replace $(I-I^0)^\alpha$ by 
$( (k + \vartheta_0/4)/\lambda)-I^0)^\alpha$ .
Hence, 
\[
\begin{array}{lcrr}
R^0(\lambda)e_k(\varphi)\ =\   e_k(\varphi)
e^{i\lambda \Phi(\varphi,(k + \vartheta_0/4)/\lambda)} \\ [0.3cm] 
\displaystyle \times \, \left(\frac{\lambda}{2\pi}\right)^{n-1}\,\int_{{\R}^{2n-2}} 
\, e^{-i\lambda\langle v, I -(k + \vartheta_0/4)/\lambda\rangle }\, 
f_k(\varphi, I,\lambda)\, dI \, dv \, + \, O_{\cal B}(|\lambda|^{-M-1})\, ,
\end{array}
\]
where 
\[
\begin{array}{lcrr}
f_k(\varphi, I,\lambda)= \sum_{|\alpha|=M+1}\, ((k + \vartheta_0/4)/\lambda -I^0)^{\alpha} 
r_{0,\alpha}(\varphi, I)\\
[0.3cm]
 \displaystyle +\,  \sum_{j=0}^{M}\, \sum_{|\alpha|=M-j}\, \lambda^{-j}\, 
 ((k + \vartheta_0/4)/\lambda -I^0)^\alpha \, r^0_{j,\alpha}(\varphi, I)\, .
 \end{array}
\]
Since $r^0_{j,\alpha}\in C_0^{n}(\T^{n-1}\times D)$   is continuous 
with respect to $K\in {\cal B}$,  
integrating $n$ times by parts with respect to $I$ in the last integral 
we gain $O_{\cal B}((1 + |\lambda v|)^{-n})$, 
and we obtain  (\ref{newrest}). \finishproof

\vspace{0.5cm} 
\noindent 
G. P.: 
Universit\'e de Nantes,  \\
Laboratoire de mathématiques Jean Leray,
CNRS: UMR 6629,\\
2, rue de la Houssini\`ere, \\
BP 92208,  44072 Nantes 
Cedex 03, France \\
%e-mail: georgi.popov@univ-nantes.fr
 
\vspace{0.5cm} 
\noindent 
P.T.: Department of Mathematics,\\
 Northeastern University,\\ 
360 Huntington Avenue, Boston, MA 02115 

\end{document}